\documentclass[12pt]{article}
\usepackage{amssymb}
\usepackage{latexsym}
\usepackage{amscd}

\usepackage{amsmath}
\usepackage[english]{babel}
\usepackage{a4}
\usepackage{epsf}
\usepackage{curves}
\usepackage{graphics}
\usepackage{rotating}
\input xypic

\usepackage{graphics} 
\usepackage{array} 
\usepackage{enumerate}
\input xy
\xyoption{all}

\def\a{\alpha}
\def\b{\beta}
\def\g{\gamma}
\def\G{\Gamma}
\def\d{\delta}
\def\e{\eta}

\def\ape{\varepsilon}
\def\la{\lambda}

\def\pat{\partial}
\def\o{\omega}
\def\k{\kappa}
\def\O{\Omega}
\def\na{\nabla}
\def\tr{\triangle}

\def\r{\rho}
\def\phi{\varphi}
\def\ta{\theta}
\newcommand{\osc}{\operatorname{osc}}

\newcommand{\sbs}{\subset}
\newcommand{\qed}{\hfill\rule{5pt}{5pt}}
\title{Novikov-Morse Theory for Dynamical Systems }
\author{\large  HuiJun Fan \qquad J\"urgen Jost\\  \\
 Max-Planck-Institut f\"ur Mathematik\\ in den Naturwissenschaften, Inselstr. 22-26,\\ D-04103 Leipzig,
 Germany}
\date{}
\begin{document}
\maketitle
\vskip 20pt
\begin{center}
\begin{minipage}{120mm}
{
{\bf Abstract}$\;\;$ The present paper contains an interpretation and generalization of Novikov's theory for Morse type inequalities for closed 1-forms in terms of concepts from Conley's theory for dynamical systems. We introduce the concept of a flow carrying a cocycle $\a$, (generalized) $\a$-flow for short, where $\a$ is a cocycle in bounded Alexander-Spanier cohomology theory. Gradient-like flows can then be characterized as flows carrying a trivial cocycle. We also define $\a$-Morse-Smale flows that allow the existence of ``cycles'' in contrast to the usual Morse-Smale flows. $\a$-flows without fixed points carry not only a cocycle, but a cohomology class, in the sense of [8], and we shall deduce a vanishing theorem for generalized Novikov numbers in that situation. By passing to a suitable cover of the underlying compact polyhedron adapted to the cocycle $\a$, we construct a so-called $\pi$-Morse decomposition for an $\a$-flow. On this basis, we can use the Conley index to derive generalized Novikov-Morse inequalitites, extending those of M. Farber [13]. In particular, these inequalities include both the classical Morse type inequalities (corresponding to the case when $\a$ is a coboundary) as well as the Novikov type inequalities ( when $\a$ is a nontrivial cocycle).
}   
\end{minipage}
\end{center}
\vskip 20pt
\noindent
{\bf\large {1 Introduction}}\\
  \\
\indent
Let $X$ be an $m$-dimensional compact manifold and let $f$ be a Morse function on $X$. Define $S(f)$ to be the union of the sets $S_i(f)$ which contain the critical points of $f$ with index $i$. Define 
\begin{align*}
&c_i:=\# S_i(f)\\
&b_i:=\dim H_i(X,{\Bbb R})
\end{align*}
Then the classical Morse inequalities are
\begin{align*}
c_i\ge & b_i\\
\sum^i_{j=0}(-1)^{i-j}c_j\ge& \sum^i_{j=0}(-1)^{i-j}b_j\tag{$1.1$}
\end{align*}
for $i=0,1\cdots,m$.
\indent
The Morse inequalities provide an important connection between the analytic information and the topological information of $X$ (see [17]). The method to obtain this result is to consider the change of the homotopy type of the level set
$$
X^a=\{x\in X|f(x)\le a\}
$$
when $a$ crosses the critical value. If $c$ is a critical value, then there is the homotopy decomposition 
$$
X^{a+\ape}\simeq X^{a-\ape}\cup_{p_j\in f^{-1}(c)\cap S(f)} e(p_j)
$$
where $e(p_j)$ is a $j$-cell corresponding to the critical point $p_j$ of index $j$. Hence, $X$ has a CW-decomposition
$$
X=\cup_{p_j\in S(f)}e(p_j)
$$
\indent
However, this decomposition only provides information about the dimension of the basis. It is not clear how each cell is attached. The successive work by R. Thom, S. Smale and J. W. Milnor in the 50's and 60's lead to the construction of the Morse complex which includes also the boundary operators.\\
\indent
Let us describe the Morse complex simply. Consider the negative gradient flow $v$ of $f$:
$$
\dot{v}(t)=-\na f(v(t)).
$$
Then for each critical point $p\in S(f)$, there are the stable and unstable manifolds
\begin{align*}
&W^u(p)=\{x\in X|\;x\cdot t\to p, \;\hbox{as}\;t\to -\infty\}\\
&W^s(p)=\{x\in X|\;x\cdot t\to p, \;\hbox{as}\;t\to \infty\}
\end{align*}
R. Thom gave a disjoint decomposition of $X$ by unstable manifolds
$$
X=\cup_{p\in S(f)}W^u(p)
$$
In [26], S. Smale obtained an important result for the gradient flow on a compact manifold. The result is that generically $W^u(q)$ intersects transversally with $W^s(p)$ for each point pair $(p,q)\in S(f)\times S(f)$. Therefore we can choose a generic Riemannian metric such that for each critical point pair $(q,p)\in S_i\times S_{i-1}$, $W^u(q)\cap W^s(p)$ is 1-dimensional and is transversal to any level set $f^{-1}(a)$ between them. This condition for the flow is called the Morse-Smale condition.\\
\indent
Since $f^{-1}(a)$ is compact, $W^u(q)\cap W^s(p)\cap f^{-1}(a)$ is a finite set. Assign arbitrary orientations to the unstable manifolds of all the critical points. Now we can assign a number $\ape(\g)$ to each trajectory $\g\in W^u(q)\cap W^s(p)$ which takes its value in $\{\pm 1\}$. Let $x\in \g\cap f^{-1}(a)$. Fix the direction $\dot{\g}(t)$, then $W^u_x(q), W^s_x(p)$ and $T_x X$ will give the induced orientation to $T_x(f^{-1}(a)\cap W^u(q)), T_x(f^{-1}(a)\cap W^s(p))$ and $T_x(f^{-1}(a))$. If the orientation of $(T_x(f^{-1}(a)\cap W^u(q)), T_x(f^{-1}(a)\cap W^s(p)))$ is the same as the one of $T_x(f^{-1}(a))$, then we let $\ape(\g)=1$, otherwise, let $\ape(\g)=-1$.\\
\indent
Define the incidence coefficient $n(q,p)$ to be 
\begin{align*}
n(q,p)=\sum_{\g\in W^u(q)\cap W^s(p)}\ape(\g).
\end{align*}
Now we can define the boundary operator $\pat_i:S_i(f)\to S_{i-1}(f)$ as follows. For each critical point $q\in S_i(f)$,
$$
\pat_i(q)=\sum_{p\in S_{i-1}(f)}n(q,p)\cdot p
$$
It can be proved that $\pat_{i-1}\pat_i=0$. Therefore $(S_i(f),\pat_i)$ becomes a complex. This is the Morse complex. The important thing is that the homology groups of the Morse complex are isomorphic to the standard homology groups of the underlying compact manifold. The reader can consult [24] for more information about the Morse homology theory.\\
\indent
In fact, the Morse complex was first used in [30] to prove the Morse inequalities (1.1). The Witten deformation technique used there provides an analytic way to prove the Morse type inequalities (see [4], [3], etc.), even the Novikov inequalities introduced below (see [6]).\\
\indent
Let us also remark that Morse theory has been generalized to the infinite dimensional case. In that situation, one needs to impose a compactness condition on the funtional, the so-called Palais-Smale condition. See e. g. the books by K. C. Chang [7] and M. Struwe [29].\\
\indent
In the 80's, S. P. Novikov considered in [18] and [19] Morse type inequalities which are used to count the zeros of a closed Morse 1-form $\o$ on $X$. By the Poincar\'e lemma, a closed 1-form is locally exact. For instance, we can assume that $\o=df_p$ near a zero point $p$ of $\o$. So closed Morse 1-form means that each zero point $p$ of $\o$ is a non-degenerate critical point of $f_p$ and has a Morse index from $f_p$. These inequalities are similar to (1.1) formally. The difference is that the ordinary Betti number $b_i$ at the right hand side of (1.1) is replaced by the Novikov number $b_i([\o])$ which is determined by $X$ and the cohomology class $[\o]$. These Morse type inequalities involving Novikov numbers are called Novikov inequalities.\\
\indent
Novikov inequalities are obtained by constructing the Novikov complex $(C_*,\pat)$ with respect to the closed Morse 1-form $\o$. We can give a simple description of the Novikov complex for $\hbox{rank}\;\o=1$. In this case, $\o$ can be lifted to the covering space $\bar{X}$ with the deck transformation group being the infinite cyclic group ${\Bbb Z}=\langle t \rangle$ such that the pull-back form by the projection map is an exact form $d\bar{f}$ where $\bar{f}$ is a Morse function on $\bar{X}$. Introduce the Novikov ring ${\Bbb Z}((t)):={\Bbb Z}[[t]][t^{-1}]$. Now the module $C_i$ is defined as the free ${\Bbb Z}((t))$-module generated by $S_i(\o)$, i.e., the set of all the zeros of $\o$ with index $i$. As for the Morse complex, we can consider the flow generated by the dual vector field of $\o$ satisfying the Morse-Smale condition and choose an orientation for each unstable manifold for each zero in $S(\o)$. Consequently, we are able to define the boundary operator of the Novikov complex. Let $(q,p)\in S_i(\o)\times S_{i-1}(\o)$ and choose any lifting points $\bar{q}$ and $\bar{p}$. Let $T_k(\bar{q},\bar{p})$ be the finite set of all the trajectories joining $\bar{q}$ and $\bar{p}\cdot t^k$. Define
\begin{align*}
&n_k(q,p)=\sum_{\g\in T_k(\bar{q},\bar{p})}\ape(\g)\\
&n(q,p)=\sum_k n_k(q,p)\cdot t^k
\end{align*}
Note that $n(q,p)\in {\Bbb Z}((t))$. $n(q,p)$ is called the Novikov incidence coefficient.\\
\indent
Finally we can define the boundary operator $\pat_i:S_i(\o)\to S_{i-1}(\o)$ to be
$$
\pat_i(q)=\sum_{p\in S_{i-1}(\o)}n(q,p)p
$$
There are some discussions about the properties of the Novikov incidence coefficients, for instance, the rationality and the growth properties of $n_k(q,p)$ (see [21], [22] and the references there).\\
\indent
The above construction of the Novikov complex stems from the analysis of the gradient flow generated by $\bar{f}$ in the noncompact manifold $\bar{X}$. However, we can define the Novikov complex by triangulating the underlying manifold $X$ and then studying the action of the boundary operators in the ${\Bbb Z}$-CW complex in the covering space $\bar{X}$. In more detail, we can let $C_*(X)$ be the simplicial complex of $X$. After lifting this triangulation to the covering space, we get $C_*(\bar{X})$ which can be viewed as a free ${\Bbb Z}$-module chain complex. Tensor $C_*{\bar{X}}$ with the Novikov ring ${\Bbb Z}((t))$, to get a complex ${\Bbb Z}((t))\otimes_{{\Bbb Z}[t,t^{-1}]}C_*(\bar{X})$, because the boundary operator in $C_*(\bar{X})$ can be viewed as a matrix operator with each entry being a Laurent polynomial. This is also called the Novikov complex. It was proved in [23] that there is a closed Morse 1-form $\o$ such that the second kind Novikov complex can be realized by the first Novikov complex given by $\o$. The construction of the Novikov complex with respect to a Morse closed 1-form implies the Novikov inequalities.\\
\indent
Since the Novikov ring is a principle ideal domain (see [16]), hence the homology module of the Novikov complex can be decomposed into two parts, the free module part and the torsion part. The Novikov numbers are refered to be the ranks of the homology modules, and the number of the minimal generators of the torsion part.\\
\indent
The Novikov numbers can also be obtained from another chain complex with coefficient ring different from the Novikov ring. One way is to consider the ring of rational functions: $\mathcal{ R}:=(1+t{\Bbb Z}[t])^{-1}{\Bbb Z}[t,t^{-1}]$. In [11], M. Farber has proved that the homology of $\mathcal{ R} \otimes_{{\Bbb Z}[t,t^{-1}]}C_*(\bar{X})$ is isomorphic to $H_*((C_*,\pat))$. We will use this kind of definition of the Novikov numbers in this paper.\\
\indent
Recently, M. Farber ([14]) has used a noncommutative localization method to construct a universal complex, which can induce many kinds of Novikov complexes.\\
\indent
The Morse complex provides a way to construct the underlying manifold by attaching handles, and we know not only how many and  which kinds of handles should be attached but also how to attach handles. Actually, we can get more refined Morse inequalities than (1.1) that contain the information from the torsion part of the homology group which is determined by the boundary operators. The construction of the Morse complex is based on the analysis of the gradient flow of the Morse function. From the point of view of dynamical systems, the information contained in the right hand side of (1.1) is nothing but the local information (Morse index) of the fixed points of the gradient flow of the Morse function $f$. Fixed points are only special invariant sets in a dynamical system. If we are given an arbitrary flow in a compact manifold and its invariant sets, can we get Morse-type inequalities that contain the local information of the invariant sets at the right hand side and the topological information of the underlying manifold at the left hand side? Of course, since a general flow is very complicated, it is very difficult to do so. However, C. C. Conley successfully generalized the concepts of nondegenerate critical points and the corresponding Morse index to the isolated invariant sets and the Conley index (see [9]). The isolated invariant sets contain many interesting invariant sets, e.g., the Bott type nondegenerate critical manifolds of the gradient flow of a smooth functions (see [4][5]), and the hyperbolic periodic orbits, or more general the hyperbolic sets. The Conley index characterizes the ``local'' topological information of the isolated invariant sets. C. C. Conley and E. Zehnder also obtained a generalization of (1.1) (see [10]). They consider the flow restricted to an isolated invariant set $S$ and assume that $S$ has a Morse decomposition $\{M_n,\cdots,M_1\}$. It can be proved that each Morse set $M_i$ is isolated and hence has a Conley index. Now the Morse type inequalities provide a connection between the Conl
ey index of $M_i$ and the Conley index of $S$.\\   
\indent
However, in order to prove the Morse type inequalities, the flow must have a ``global'' gradient-like structure with respect to all the invariant sets in $X$. ``Global'' here means that if we collapse all the invariant sets into points, then the quotient flow is gradient-like. Only with a ``global'' gradient-like flow , we can get the Morse decomposition.\\
\indent
Now we return to the Novikov inequalities. We can construct a closed Morse 1-form such that its flow has no attractor and no repeller (see [19]). Therefore the flow generated by this 1-form is not a gradient-like flow in $X$ though its lifting flow on $\bar{X}$ is. Therefore, it is not possible to give a Morse decomposition for the whole manifold $X$.\\
\indent
After an introduction to the Morse theory for a gradient-like flow and the Novikov theory for a Morse closed 1-form, we reach the goal in the present paper. In this paper, we will introduce the concept, ``flow carrying a cocycle $\a$ ($\a$-flow and generalized $\a$-flow)'', where $\a$ is a 1-dimensional cocycle in the bounded Alexander-Spanier cohomology which will be defined in section 3.1. When $\a$ is a coboundary, then this flow is a gradient-like flow. Conversely, if a flow is a gradient-like flow ,then this flow carries a coboundary $\d g$ where $g$ is a Lyapunov function. Moreover, in theorems 3.4.2 and 3.4.3, we give a sufficient and necessary conditon for a flow to be a non gradient-flow if $\a$ is a nontrivial cocycle. The notion of $\a$-flow contains many interesting flows. The gradient flow of a Morse function $f$ is a $df$-flow. The flow generated by a closed 1-form $\o$ is an $\o$-flow and the pertubation flows are also $\o$-flows. If an $\a$-flow has no fixed point, then this flow is a flow ``carrying a cohomology class'' which is introduced in [8]. A generalized $\a$-flow has a ``global'' $\a$-flow structure with respect to the $\pi$-stable nonwandering sets. The $\pi$-stable nonwandering sets are those sets that when lifted to the covering space related to $\a$, are also bounded sets. An example of a generalized $\a$-flow is the $\a$-Morse-Smale flow. In the definition of the Morse-Smale flow, one requires that such a flow has no ``cycles'' connecting the different hyperbolic fixed points or hyperbolic periodic orbits. The $\a$-Morse-Smale flow is a Morse-Smale flow when lifted to the covering space with respect to the cocyle $\a$. However, in the underlying manifold, it is possible that there exist ``cycles'' connecting different nonwandering sets or there exist  $(\mathcal{ U}(s), T)$-chain cycles for arbitrary small $\ape>0$ connecting different nonwandering sets. The existence of such ``cycles'' or such chain cycles makes the flow a nongradient-like flow by theorem 3.4.2.\\
\indent
A flow carrying a cocycle has a $\pi$-Morse decomposition (Theorem 3.5.2). Hence there are the Novikov-Morse type inequalities, Theorem 4.5.1, if $\mathcal{ IN}_\pi=\mathcal{ N}_\pi$. This theorem uses the Conley index to get the ``local'' topological information of the isolated invariant sets.\\
\indent
The inequalities (4.5.1) are a generalization of the Novikov inequalities for a closed Morse 1-form. Starting from these inequalities, we can recover many Novikov type inequalities found before. For example, we can consider the $\o$-flow generated by a smooth 1-form $\o$ having Bott type nondegenerate zero sets, then (4.5.1) induces (4.6.2), which was essentially given in [13]. If we consider an $\a$-Morse Smale flow, then we have Theorem 4.6.6 which provides new Novikov inequalities. Also, the existence of a flow carrying a cohomology class $\a$ implies the vanishing of the Novikov numbers $b_i([\o])$. In addition, if we consider the trivial cocycle $\a$, then Theorem 4.5.1 provides the Morse type inequalities, because all the representations involving $\a$ are trivial.\\   
\indent
Since the isolated invariant sets are more complicated than a hyperbolic fixed point, hence it is not clear what the connecting orbits in the $\pi$-Morse decomposition are. The appendix gives another proof of the classical Novikov inequalities starting from the equality (4.5.1). These Novikov inequalities also contain the information about the torsion parts of the homology groups of the Novikov complex.\\
  \\
\noindent
{\bf Structure of this paper}\\
 \\
1 Introduction\\
\noindent
2 Basic results about flows\\
\indent
2.1 Basic concepts and structures of flows\\
\indent
2.2 Local structure of gradient-like flow\\
\noindent
3 Flows carrying a cocycle $\a$\\
\indent
3.1 Bounded Alexander-Spanier cohomology\\
\indent
3.2 Integration of cocylces along chains\\
\indent
3.3 Integration of cocycles in a gradient-like flow\\
\indent
3.4 Flows carrying a cocycle $\a$\\
\indent
3.5 $\pi$-Morse decomposition of $\a$-flows\\
\noindent
4 Novikov-Morse theory\\
\indent
4.1 Deformation of complexes relative to $\a$-flows\\
\indent
4.2 Monodromy representations and Novikov numbers\\
\indent
4.3 Morse type inequalities for a filtration\\
\indent
4.4 Isolated invariant sets and Conley index\\
\indent
4.5 Novikov-Morse type inequalities for flows carrying a cocycle\\
\indent
4.6 Applications to special flows carrying a cocycle\\
\noindent
Acknowledgement\\
\noindent
References\\
\noindent
Appendix\\
  \\
\noindent
{\bf \large 2 Basic results about flows}\\
\ \\
{\bf 2.1 Basic concepts and structures of flows}\\
\par
Let $X$ be a topological space. The continuous function from $X\times{\Bbb R}\rightarrow X$ given by $(x, t)\rightarrow x\cdot t$ is called a flow on $X$ if for all $x\in X$, and $s, t\in {\Bbb R}$ it satisfies (1) $x\cdot 0=x$ and (2) $(x\cdot s)\cdot t=x\cdot(s+t)$. For $Y\subset X$ and ${\Bbb R}^{\prime}\subset {\Bbb R}$, let $Y\cdot{\Bbb R}^{\prime}$ denote the set of points $x\cdot t$ such that $y\in Y$ and $t\in {\Bbb R}^{\prime}$.\\
\indent
A subset $I\sbs X$ is called invariant if $I\cdot{\Bbb R}=I$. If $N\sbs X$ is a subset, we denote by $I(N)$ the invariant set in $N:I(N)=\{x\in N|x\cdot{\Bbb R}\sbs N\}$. Clearly, $I(N)$ is invariant and it is closed if $N$ is closed.\\
\indent
For $Y\subset X$, the $\o-$limit set $\o(Y)$ of $Y$ is defined to be the maximal invariant set in the closure of $Y\cdot[0, \infty]$. If $\overline{Y\cdot[0, \infty]}$ denotes the closure in $X$, then $\o(Y)=I(\overline{Y\cdot[0, \infty])}$. Similarly, the $\o^*-$limit set $\o^*(Y)$ of $Y$ is defined to be the maximal invariant set in the closure of $Y\cdot[0, -\infty]$.\\
\indent
The $\o-$limit set $\o(Y)$ has some properties:\\
\qquad A. $\o(Y)$ is a closed set.\\
\qquad B. $\o(Y)=\cap \{\overline{ Y\cdot[t, \infty)}|t\ge 0\}.$\\
\qquad C. Let $S$ be a closed compact Hausdorff invariant set in $X$. If $Y\subset S$, then $\o(Y)$ is a nonempty compact invariant set in $S$ and is connected if $Y$ is.\\
\indent
$\o^*(Y)$ has analogous properties, since $\o^*(Y)$ is the $\o$-limit set of the inverse flow. The above conclusions can be easily obtained from the definition of the $\o$-limit set.\\
\indent
Let $S$ be an invariant set in $X$. The nonwandering set $N_\o(S)$ consists of the points $x$ in $S$ with the property that if $U$ is a neighborhood of $x$ in $S$, then $x\in\o(U)$. The nonwandering set is invariant and closed in $S$. Since the definition of $N_\o(S)$ depends on the invariant set $S$, it is not generally true that $N_\o(S)=N_\o(N_\o(S))$. For $x\in S, \o(x)\subset N_\o(S)$.\\
\noindent
{\bf Attractor-repeller pairs}\\
\indent
Let $S$ be a compact Hausdorff invariant set in $X$. A subset $A\subset S$ is called an attractor (relative to $S$) if there is a neighborhood $U$ of $A$ in $S$ such that $\o(U)=A$. Similarly a set which is the $\o^*$-limit set of one of its neighborhoods is called a repeller. Let $A$ be an attractor in $S$, and let $A^*:=\{x\in S| \o(x)\not\in A\}$. Then $A^*$ is a repeller in $S$ which is called the complementary repeller in $S$. $(A, A^*)$ is called an attractor-repeller pair in $S$. $C(A, A^*)=S\backslash A\cup A^*$ is called the set of connecting orbits to the pair, since the point $x\in C(A, A^*)$ if and only if $\o(x)\in A$ and $\o^*(x)\in A^*$.\\
\indent
There are two criteria to judge if a set is an attractor.\\
\indent
(1) Suppose $U\subset S$ and, for some $t>0,\,(\overline{U})\cdot t\subset\mbox{int}U$. Then $\o(U)$ is an attractor contained in the interior of $U$.\\
(2) Let $U$ be a compact set in $S$. If all its boundary points leave $S$ in backward time, then the maximal invariant set in $U$ is an attractor.\\
\indent
If $(A,A^*)$ is an attractor-repeller pair in $S$ and $S^\prime$ is an invariant subset in $S$, then $(A\cap S',\,A^*\cap S')$ is an attractor-repeller pair in $S'$.\\
\noindent
{\bf Chain recurrent set}\\
\indent
Let $\mathcal{U}$ be a cover of $S$ and let $x,\,x^\prime\in S$, A $(\mathcal{U},\,t)$-chain from $x$ to $x^\prime$ is a sequence $\{x=x_1,\,\cdots,\,x_{n+1}=x^\prime|t_1,\cdots,t_n\}$ such that $t_i\ge t$ and each pair $(x_i\cdot t_i, x_{i+1})(i=1,\cdots,n)$ belongs to the same open set in $\mathcal{U}$.\\
\indent
Let $Q(S)\subset S\times S$ be the set of pairs $(x,x')$ where for any covering $\mathcal{U}$ and any $t>0$ there is a $(\mathcal{U},t)$-chain from $x$ to $x'$.\\
\indent
For $Y\subset S$, define 
$$
\aligned
\O(Y):= & \{x^\prime|\hbox{there is }x\in Y \hbox{s.t.}(x,x^\prime)\in Q(S)\}\\
\O^*(Y):= & \{x^\prime|\hbox{there is }x\in Y \hbox{s.t.}(x^\prime,x)\in Q(S)\}
\endaligned
$$
\indent
$Q(S)$ is a closed, transitive relation on $S$, and if $(x,x^\prime)\in Q(S)$, then $(x\cdot I,x^\prime\cdot I)\in Q(S)$, where $I$ is any interval in the real line. For example, we prove the last conclusion.\\
  \\
Proof $\;\;$ Let $\mathcal{U}$ be any open covering of $S$ and take any positive $t_0,t$ and $t^\prime$. Suppose that $(x,x^\prime)\in Q(S)$, we want to prove that $(x\cdot t, x^\prime\cdot t^\prime)\in Q(S)$.\\
\indent
Let $\bar{t}=\max{t_0+|t|,t_0+|t^\prime|}$, and $\{x=x_1,\,\cdots,\,x_{n+1}=x^\prime|t_1,\cdots,t_n\}$ be a $(\mathcal{U},\bar{t})$-chain from $x$ to $x^\prime$, then $\{x\cdot t=x_1,\cdots,x_n,x_{n+1}=x'\cdot t'|t_1-t,\cdots,t_{n-1},t_n+t'\}$ is a $(\mathcal{U},t_0)$-chain from $x\cdot t$ to $x'\cdot t'$.\qed\\
  \\
\indent
$\O(Y)$ and $\O^*(Y)$ have the following two properties:
\\
\indent
(1) If $Y$ is closed, $\O(Y)$ and $\O^*(Y)$ are closed invariant sets containing $\o(Y)$ and $\o^*(Y)$, resp.\\
\indent
(2) If $Y$ is closed, $\O(Y)$ is the intersection of the attractors in $S$ which contain $\o(Y)$; similarly, $\O^*(Y)$ is the intersection of the repellers containing $\o^*(Y)$.\\
\indent
The chain recurrent set $R(S)$ is defined to be the set of $x$ such that $(x,x)\in Q(S)$, or equivalently $x\in\O(x)$. The relation of chain recurrent sets with the attractor-repeller pairs and the nonwandering sets are as follows:
\\
\indent
(1) $R(S)=\cap \{A\cap A^*|A \mbox{ is an attractor in }S\}$\\
\indent
(2) $R(S)\supset N_\o(S)$
\\
  \\
{\bf Definition} $\;\;$ A invariant set $S$ is called chain recurrent if $R(S)=S$. A flow in $X$ is called  a chain recurrent flow if its chain recurrent set is the whole space.\\
  \\
\noindent
{\bf Morse decomposition}\\
\indent
Let $S$ be a compact Hausdorff invariant set in $X$. The intersection of an attractor and a repeller in $S$ is called a Morse set (rel $S$).\\
\indent
Let $\emptyset=A_{n+1}\subset A_n\subset\cdots\subset A_0=S$ be a decreasing sequence of attractors in $S$ and let $M_j=A_j\cap A^*_{j+1}(j=1,2,\cdots,n)$. Then the ordered set $D=\{M_n,M_{n-1},\cdots,M_1\}$ is called a Morse decomposition of $S$. Clearly, the Morse sets in $D$ are disjoint and for any $x\in S$, either $x$ is in a Morse set, or there is $i,j$ satisfying $i<j$ such that $\o(x)\subset M_j$ and $\o^*(x)\subset M_i$. This means that each point in a connection orbit of $S$ goes from a Morse set with low index to a Morse set with high index.\\
\indent
The Morse decomposition is not unique. For example, let $1\le i<j\le n$, consider the decreasing sequence of attractors: $\emptyset=A_{n+1}\subset A_n\subset\cdots\subset A_0=S$, then $D_{ji}=\{M_0,\cdots,M_j,M_{ji},M_{i-1},\cdots,M_1\}$ is also a decomposition of $S$. Here $M_{ji}=A_i\cap A^*_j$.\\
\indent
Conversely, let $\{M_n,\cdots,M_1\}$ be any collection of disjoint invariant sets in $S$. If for any $x\in S$, there is $i<j$ such that $\o^*(x)\subset M_i$ and $\o(x)\subset M_j$, then it can be proved that $D=\{M_n,\cdots,M_1\}$ is a Morse decompostion.\\
\indent
Let $M(D)=M_n\cup\cdots\cup M_1$. Then for any decomposition, we have $M(D)=\cap^n_{i=1}\{A_i\cup A_i^*\}$. Hence $M(D)\supset R(S)$.\\
\indent
Suppose that the nonwandering set $N_\o(S)$ has finitely many components, then these can be ordered and become a Morse decomposition with $M(D)=N_\o(S)$. But generally $M(D)\supset R(S)\supset N_\o(S)$, hence in this case $M(D)=R(S)=N_\o(S)$.\\
\noindent
{\bf Gradient-like flows and chain recurrent flows}\\
\indent
A flow is called gradient-like, if there is a continuous real valued function which is strictly decreasing on non constant solutions. Such a function is called a Lyapunouv function. The concept of a gradient-like flow is a generalization of one of a gradient flow, i.e., a flow of the form:
$$
\dot{v}=\text{grad}f(v)
$$
where $f$ is a $C^1$-function defined on $X$. Clearly, $-f$ decreases along the nonconstant solutions as $t$ increases.\\
\indent
A gradient-like flow need not be a gradient flow.\\
\indent
A flow is called a strongly gradient-like flow if its chain recurrent set is totally disconnected. A gradient-like flow with finitely many fixed points is a strongly gradient-like flow.\\
\indent
We have the basic structure theorem of flows (see [9]).\\
 \\
\noindent
{\sl {\bf Theorem 2.1}$\;\;$ Every flow on a compact space is uniquely represented as the extension of a chain recurrent flow by a strongly gradient-like flow.}\\
\par
This theorem means that if we collapse every connected component of the chain recurrent set to one point, then we get a quotient flow on the collapsed space which is a strongly gradient-like flow. The flow when restricted to every component of the chain recurrent set becomes a chain recurrent subflow.\\
 \\
\noindent
{\bf 2.2 Local structure of gradient-like flows}\\
\par
This section will provide a detailed description of the local structure near a fixed point of a gradient-like flow. The conclusion obtained here will be applied in sections 3.3 and 3.4.\\
\indent
Let $v$ be a gradient-like flow on a compact metric space $X$ having $n$ fixed points $p_i,\;i=1,\cdots, n$. 
Then there exists an associated Lyapunov function $g(x)$ on $X$ such that for a small $r>0$, the sets 
$g(B_r(p_i)), i=1,\cdots,n$ do not intersect mutually.\\
\indent
Let $p$ be a fixed point of the flow $v$, and let $B_{\frac{r}{2}}(p)$ be a closed ball centered at $p$ with radius $\frac{r}{2}$. 
For $0<s<\frac{r}{2}$, we define two sets on the sphere $\pat B_r(p)$, 
$$
\aligned
&B^+_{r,s}(p)=\{x\in\partial B_r(p)|[x,x\cdot t]\cap \pat B_s(p)\neq\emptyset\;\hbox{for}\; t>0\}\\
&B^-_{r,s}(p)=\{x\in\partial B_r(p)|[x,x\cdot t]\cap \pat B_s(p)\neq\emptyset\;\hbox{for}\; t<0\}
\endaligned
$$
\indent
For any point $x\in B^+_{r,s}(p)$, let $t_x$ denote the arrival time of the trajectory $[x, x\cdot t](t>0)$ 
to the sphere $\pat B_{\frac{r}{2}}(p)$. We have the following lemma.\\
  \\
{\sl {\bf Lemma 2.2.1}$\;\;$ (1) $B^\pm_{r,s}(p)$ are closed sets on $\pat B_r(p)$.\\
(2) $t_x$ is a lower semicontinuous function on $B^+_{r,s}(p)$.\\
}  
 \\
Proof $\;\;$ To prove (1), we need only to prove that the set $\pat B_r(p)-B^+_{r,s}(p)$ is open in 
$\pat B_r(p)$. Let $x\in \pat B_r(p)-B^+_{r,s}(p)$, then $x\cdot t(t>0)$ does not intersect
 with $\pat B_s(p)$ and will finally drop into another fixed point $q$. In terms of the choice of $g(x)$,
  we can choose $a<b$ satisfying 
\begin{equation}
B_r(p)\subset g^{-1}([b,+\infty));\; B_r(q)\subset g^{-1}((-\infty, a])   
\tag{$2.2.1$}
\end{equation}
Since $x\cdot t(t>0)$ drops into $q$, there exists a $T_x>0$ such that $x\cdot T_x$ is in the interior of 
$B_r(q)$. Therefore there is a closed ball $D_\d(x)$ such that for any $y\in D_\d(x)\cap \pat B_r(p), y\cdot
 T_x$ is in the interior of $B_r(q)$. We can choose $\d$ small enough such that the set $D_\d(x)\cdot [0,
  T_x]\cap \pat B_s(p)=\emptyset$. Since $y\cdot T_x\in g^{-1}((-\infty,a]),\;y\cdot [T_x,+\infty)\in 
  g^{-1}((-\infty,a])$. Therefore $D_\d(x)\cdot [0,+\infty)\cap\pat B_s(p)=\emptyset$
  in view of (2.2.1). This shows that $\pat
   B_r(p)-B^+_{r,s}(p)$ is open in $\pat B_r(p)$.\\
\indent
To prove (2), let $x\in B^+_{r,s}(p)$ and let $t_x<+\infty$ be the arrival time. We need only to prove that for 
any small $\ape>0$, there exists a neighborhood $D_\d(x)$ of $x$ such that for any $y\in D_\d(x)\cap B^+_{r,s}
(p)$,
\begin{equation}  
t_y\ge t_x-\ape
\tag{$2.2.2$}
\end{equation}
Since $t_x$ is the arrival time of the trajectory $x\cdot t(t>0)$ to $\pat B_{\frac{r}{2}}(p)$, $x\cdot [0, t_x-\ape]$
 has a positive distance from $\pat B_{\frac{r}{2}}(p)$. By the continuity of the flow $v$, there is a neighborhood 
 $D_x$ such that $\bar{D}_x\cdot [0, t_x-\ape]$ has a positive distance from $\pat B_{\frac{r}{2}}(p)$ as well. Hence for 
 $y\in D_x\cap B^+_{r,s}(p)$, (2.2.2) is true.\qed\\
\par 
Since the Lyapunov function $g(x)$ is strictly decreasing along any nonconstant trajectory, there exists 
a $\d_x>0$ such that $g(x)-g(x\cdot t_x)=2\d_x$. By lemma 2.2.1 (2) and the property of $g(x)$, 
the function $g(y)-g(y\cdot t_y)$ is lower semicontinuous with 
respect to $y$ on $B^+_{r,s}(p)$. Hence for 
each $x\in B^+_{r,s}(p)$, there exists a neighborhood $U_x$ in $  B^+_{r,s}(p)$ such that 
$\forall y\in U_x,\,g(y)-g(y\cdot t_y)>\d_x$. 
This implies that $\d^+_{r,s}(p):=\min_{x\in B^+_{r,s}(p)}{\d_x}$ is positive and 
$\forall x\in B^+_{r,s}(p),\,g(x)-g(x\cdot t_x)>\d^+_{r,s}(p)$. In the same way, we can obtain
 a $\d^-_{r,s}(p)>0$ such that $\forall x\in B^-_{r,s}(p),$ there is 
 $g(x\cdot(-t_x))-g(x)>\d^-_{r,s}(p)$.\\
\indent
Let $0<s_1<s_2<\frac{r}{2}$, then we have the following conclusions:
$$
\aligned
&(1)\;\;B^+_{r,s_1}(p)\subset B^+_{r,s_2}(p);B^-_{r,s_1}(p)\subset B^-_{r,s_2}(p)\\
&(2)\;\;\d^\pm_{r,s_2}(p)\le \d^\pm_{r,s_1}(p)\\
&(3)\;\;\mbox{There are $s_0>0$, such that for any $0<s<s_0$, there is } B^+_{r,s}(p)\cap B^-_{r,s}(p)=\emptyset
\endaligned
$$
The first two conclusions are obvious. For (3), we can choose $s_0>0$ satisfying $\osc_{x\in B_{s_0}(p)}g(x)<\min{\d^\pm_{r,s_0}(p)}$. If there is a $0<s<s_0$ such that $ B^+_{r,s}(p)\cap B^-_{r,s}(p)\neq\emptyset$, then this means that there is a trajectory in $\bar{B}_r(p)$ which has non empty intersection with $\partial B_r(p)$ and $\partial B_s(p)$. But this is absurd, because $g(x)$ is strictly decreacing along any trajectory. After traveling from $\partial B_s(p)$ to $\partial B_r(p)$, then back to $\partial B_s(p)$, the value of $g(x)$ will decrease at least $\d^+_{r,s}(p)+\d^-_{r,s}(p)$ which contradicts the fact that $\osc_{x\in B_s(p)}g(x)\le\osc_{x\in B_{s_0}(p)}g(x)<\min{\d^\pm_{r,s}(p)}$. Hence for $0<s<s_0$, (3) holds.\\
  \\
\indent
We introduce some notations that we will use later in this paper.\\
\indent
Let $v$ be a flow on a compact metric space $(X,d)$, whose nonwandering set has finitely many components. Denote them by $A=\{A_n,\cdots, A_1\}$. Define $B_r(A_i)(i=1,\cdots,n)$ to be the open neighborhood of $A_i$ having distance less than $r$ to $A_i$. Now consider the flow $v$ restricted to the punctured space $X\backslash\cup^n_{i=1}B_r(A_i)$. Then all the trajectories can be classified into three types:\\
\indent
(1) trajectories with domain $[a,b],-\infty<a<b<\infty$;\\
\indent
(2) trajectories with domain $[a, +\infty)$ or $[-\infty,b],-\infty<a,b<\infty$;\\
\indent
(3) $(-\infty,+\infty)$.\\
\indent
We denote the sets of three type trajectories by $\G_1^A(r),\G_2^A(r)$ and $\G_3^A(r)$ respectively. If $A_i$ is a point $p_i$ in $X$ for $i=1,\cdots,n$, then we denote $\G_i^A(r)$ simply by $\G_i(r)$.\\
\indent
Since the flow $v$ can be viewed as the action of the Lie group ${\Bbb R}$ on the compact metric space $X$, we assume in this paper that this action is Lipschitz with a finite Lipschitz constant. Therefore we can talk about the tangent vector of a trajectory of $v$ and integrate the tangent vector to get the length of a segment of a trajectory.\\
 \\
\noindent
{\bf\large 3 Flows carrying a cocycle $\a$}\\
  \\
\noindent
{\bf 3.1 Bounded Alexander-Spanier cohomology}\\
\par
In his paper [8], Churchill introduced a bounded Alexander-Spanier cohomology theory. If we consider the category of compact pairs, this bounded Alexander-Spanier cohomology theory is identical with the usual Alexander-Spanier theory. One can see the proof of equivalence in [8]. In this section we will formulate such a cohomology theory.\\
\indent
Let $\mathcal{R}$ be an ordered ring and $X$ be any set. A function $\a:\,X\rightarrow \mathcal{R}$ is called bounded if there is an element $M\in \mathcal{R}^+$, called a bound for $\a$, such that $|\a(x)|\le M$ for any $x\in X$. Let $X^{q+1}=X\times\cdots\times X(q+1\;\mbox{times},\;q\ge 0)$, and $\bar{C}^q(X,\mathcal{R})$ be the module of all bounded function from $X^{q+1}$ to $\mathcal{R}$ with addition and scalar multiplication defined pointwise. A coboundary homomorphism $\d:\;\bar{C}^q(X,\mathcal{R})\rightarrow \bar{C}^{q+1}(X,\mathcal{R})$ is defined as
$$
\d\a(x_0,\cdots,x_{q+1})=\sum^{q+1}_{i=0}(-1)^i\a(x_0,\cdots,\hat{x}_i,\cdots, x_{q+1})
$$
where $(x_0,\cdots,x_{q+1})\in X^{q+2}$ and $(x_0,\cdots,\hat{x}_i,\cdots,x_{q+1})$ denotes the $(q+1)$-tuple $(x_0,\cdots,x_{i-1}, x_{i+1},\cdots,x_{q+1})$. It is easy to prove that $\d\d=0$ and so $(\bar{C}^*(X,\mathcal{R}),\d)$ is a complex.\\
\indent
An element $\a\in \bar{C}^q(X,\mathcal{R})$ is said to be locally zero, if $\a$ as a function from $X^{q+1}$ to $\mathcal{R}$ vanishes in a neighborhood of the diagonal of $X^{q+1}$. All the locally zero elements of $\bar{C}^q(X,\mathcal{R})$ constitute a submodule of $\bar{C}^q(X,\mathcal{R})$, denoted by $C^q_0(X,\mathcal{R})$. $(C_0^*(X,\mathcal{R}),\d)$ is a subcomplex, since if $\a$ is locally zero in $X^{q+1}$, then $\d\a$ is locally zero in $X^{q+2}$. Define $C^*(X,\mathcal{R})$ to be the quotient cochain complex of $(\bar{C}^*(X,\mathcal{R}),\d)$ by $(C_0^*(X,\mathcal{R}),\d)$. Now the bounded Alexander-Spanier cohomology group $\mathcal{H}^*(X,\mathcal{R})$ is the cohomology group of the cochain complex $(C^*(X,\mathcal{R}),\d)$.\\
\indent
For any topological space $Y$ and any continuous map $f:X\rightarrow Y$, there is the induced homomorphism $f^\#:\bar{C}^*(Y,\mathcal{R})\rightarrow\bar{C}^*(X,\mathcal{R})$ defined by
$$
f^\#\a(x_0,\cdots,x_q)=\a(f(x_0),\cdots,f(x_q)),
$$
where $\a\in\bar{C}^q(Y,\mathcal{R})$ and $(x_0,\cdots,x_q)\in Y^{q+1}$. It is easily seen if $\a\in C^*_0(Y,\mathcal{R})$, then $f^\#\a\in C^*_0(X,\mathcal{R})$. Therefore $f^\#$ induces a homomorphism from $C^*(Y,\mathcal{R})$ to $C^*(X,\mathcal{R})$, and we still denote it by $f^\#$.\\
\indent
Let $i$ be the inclusion map from a subspace $A$ to $X$, then we have the cochain map $i^\#:C^*(X,\mathcal{R})\rightarrow C^*(A,\mathcal{R})$. Let $C^*(A,X,\mathcal{R}):=\ker{i^\#}$. Now the bounded Alexander-Spanier cohomology module $\mathcal{H}^*(A,X,\mathcal{R})$ is defined to be the cohomology module of the cochain complex $C^*(A,X,\mathcal{R})$.\\
\par
Let $(B,Y)$ be a topological pair and let $f:(A,X)\rightarrow (B,Y)$ be continuous, then the cochain map $f^\#:C^*(B,Y)\rightarrow C^*(A,X)$ induces the homomorphism $f^*:\mathcal{H}^*(B,Y,\mathcal{R})\rightarrow \mathcal{H}^*(A,X,\mathcal{R})$.\\
\indent
As mentioned in the beginning of this section, the cohomology theory constructed above is identical to the usual Alexander-Spanier cohomology theory when restricted to the category of compact pairs. Therefore it satisfies the Eilenberg-Steenrod axioms and the continuity which is a special property of Alexander-Spanier cohomology theory.\\
\indent
Let $[\a]\in \mathcal{H}^1(X,\mathcal{R})$, then a cocycle $\a\in [\a]$ is a bounded function from $X^2$ to $\mathcal{R}$ with $\d\a$ being locally zero. If $[\a]=0$, then the representation $\a$ is called coboundary and there is a bounded function $\b:\,X\rightarrow \mathcal{R}$ such that $\d\b-\a$ is locally zero. If this is the case, we write $\d\b\simeq \a$.\\
\indent
Let $[\a]\in \mathcal{H}^1(X,\mathcal{R})$ and $A$ be a path connected component of $X$. Suppose that $[\a]|_A$ is zero and there are two bounded function $\b$\ and $\b'$ on $A$ such that $\d\b\simeq\d\b'\simeq\a$, then $\b-\b'$ is constant on $A$. Namely, because $\d(\b-\b)$ is locally zero on $A$, for any point $x\in A$ there is a small closed neighborhood such that $\b-\b'$ is constant on it. Since $A$ is path connected, any two points on $A$ can be joined by a chain of such closed neighborhoods. Hence $\b-\b'$ is a constant on $A$. Define $I_\a(x,y):=\b(y)-\b(x)$, where $x,y\in A$. It is obvious that $I_\a(x,y)$ depends only on $[\a]$.\\
\indent
Now let $X$ be a topological space. Let $C(X)$ be the group of all the continuous maps from $X$ to $S^1$.(C(X) can also be seen as the set of continuous functions of absolute value $1$ defined on $X$. then we can talk about the angle of $f(x)$). Let $R(X)$ be the subgroup of $C(X)$ consisting of all the continuous maps having the form $\exp{2\pi iH(x)}$ with $H(x)$ being a real-valued continuous function on $X$.\\
  \\
{\sl{\bf Proposition 3.1.1[8][25]}$\;\;$ If $X$ is a compact polyhedron, then the map $J:\,C(X)/R(X)\rightarrow \mathcal{H}^1(X,{\Bbb R})$ is a group embedding, and the range of $J$ spans $\mathcal{H}^1(X,{\Bbb R})$.}\\
  \\
\indent
Let $\g$ be a curve in $X$ and let $\tr_\g\arg{f}$ be the change of angular variable of $f$ along $\g$. It is well-known that for any closed curve $\g$ in $X$, $\frac{1}{2\pi}\tr_\g\arg{f}$ is an integer depending only on the homology class of $\g$ and the coset of $f$ in the group $C(X)/R(X)$. Hence $C(X)/R(X)$ is in fact the integral 1-dimensional cohomology group.\\
\indent
Take a group of closed  curves $\{l_1,\cdots,l_m\}(m=\dim{X})$ such that the homology classes $\{[l_1],\cdots,[l_m]\}$ constitute a basis of $\mathcal{H}_1(X,{\Bbb R})$. With respect to this basis, we have the dual basis $\{[f_1],\cdots,[f_m]\}$. These circle-valued maps $f_1,\cdots,f_k$ are called rank 1 integral cocycles.\\
  \\
{\bf 3.2 Integration of cocycles along chains}\\
  \\
\indent
Let $(X,d)$ be a compact metric space with metric $d$. For any point $p\in X$, we can choose a sequence of open balls, $\bar{B}_1\supset B_1\supset \bar{B}_2\supset\cdots$ such that $\cap^\infty_{i=1}B_i=\{p\}$. Since the bounded Alexander-Spanier cohomology has continuity, it follows that $\lim_{i\rightarrow\infty}H^1(\bar{B_i})=H^1(\{p\})=0$. Therefore for a fixed cocycle $\a$ , there is a small open ball $B_p$ such that $\a|_{\bar{B}_p\times \bar{B}_p}$ is a coboundary. Then there is a bounded function $\b_p$ defined on $\bar{B}_p$ satisfying $\d\b_p\simeq\a|_{\bar{B}_p\times\bar{B}_p}$. Since $X$ is compact, for any arbitrary small $r>0$, there are finitely many open balls $\{B_j\}^l_{j=1}$ with radius $r$ covering $X$. Furthermore, on every closed ball $\bar{B}_j$, there is $\b_j$ such that $\a|_{\bar{B}_j\times\bar{B}_j}\simeq\d\b_j$ and $|\b_j|$ has an upper bound $M_j$ on $\bar{B}_j$. We call this covering an $r$-covering and denote it by $\mathcal{U}(r)$. Define $M_\a(r)=\max_j{M_j}$.\\
\indent
Take a curve $\g:[a,b]\rightarrow X$, and choose a partition $a=t_0<t_1<\cdots<t_k=b$ such that for $0\le j\le k-1,\,\g([t_j, t_{j+1}])\subset B_{i_j}$ for some $1\le i_j\le l$. Now we can define the integration of $\a$ along $\g$ to be
$$
\int_\g\a=\sum^{k-1}_{j=0}\b_{i_j}(\g(t_{j+1}))-\b_{i_j}(\g(t_j))
$$
Notice when $\g([a,b])\subset B_{i_0}$, then $\int_\g\a=I_\a(\g(a),\g(b))$.\\
  \\
{\sl {\bf Proposition 3.2.1[8]}$\;\;$ The integral $\int_\g\a$ is independent of the partition of $\g$, the covering $\{B_i\}^l_{i=1}$ and the functions $\{\b_i\}$ selected on $\{\bar{B}_i\}$. It depends only on the cocycle $\a$ and the relative homotopy class of $\g$.
}\\
 \\
\indent
Let $v$ be a flow on $X$ and choose an $r$-covering $\mathcal{U}$ for small $r$. Let $\g=\{x_0=x,x_1,\cdots,x_k=y|t_0,\cdots,t_{k-1}\}$ be a $(\mathcal{U}, T)$-chain from $x$ to $y$. We define the integral of $\a$ along the chain $\g$ as 
$$
\int_\g\a=\sum^{k-1}_{j=0}\int_{[x_j,x_j\cdot t_j]}\a+\sum^{k-1}_{j=0}\b_{i_j}(x_{j+1})-\b_{i_j}(x_j\cdot t_j),
$$
where $(x_j\cdot t_j,x_{j+1})$ lies in $B_{i_j}$ for some $B_{i_j}\in\mathcal{U}$.\\
\indent
If we connect point pairs $(x_j\cdot t_j, x_{j+1})(j=1,\cdots,k-1)$ by any curve lying in $B_{i_j}$, we get an induced curve $\tilde{\g}$ from $x$ to $y$. It is obvious that 
$$
\int_\g\a=\int_{\tilde{\g}}\a
$$
It is easily seen that the integration $\int_\g\a$ has the following properties:\\
\indent
(1) If $\la\in\mathcal{R}$, then $\int_\g\la\a=\la\int_\g\a$.\\ 
\indent
(2) If $\g$ is a trivial curve, i.e., represented by one point, then $\int_\g\a=0$\\
\indent
(3) Let $\g_1$ and $\g_2$ be two curves for which the end point of $\g_1$ is the starting point of $\g_2$. Let $\g_1\ast\g_2$ denote the curve $\g_1$ followed by $\g_2$, then 
$$
\int_{\g_1\ast\g_2}\a=\int_{\g_1}\a+\int_{\g_2}\a\
$$
\\ \\
{\sl {\bf Proposition 3.2.2}$\;\;$ Let $\a$ be a cocycle on the compact metric space $(X,d)$. Let $\g$ be a curve on $M$ with length $l(\g)\le L$, then there is a constant $C_L$ depending only on $L$ and $\{M_j\}$ such that 
$$
\Big|{\int_\g\a}\Big|\le C_L
$$
}
  \\
Proof$\;\;$ With respect to $\a$ there is an $r$-covering $\mathcal{U}=\{B_j\}$ such that $\a|_{\bar{B}_j\times\bar{B}_j}\simeq\d\b_j$ and $|\b_j|$ has upper bound $M_j$ on $\bar{B}_j$. Since $\g$ is a curve with length $L$, we can use at most $[\frac{L}{2r}]+1$ $r$-balls to cover $r$. Therefore
$$
\Big|\int_\g\a\Big|\le 2([\frac{L}{2r}]+1)\max_{B_j\in\mathcal{U}}M_j
$$
\qed
  \\
\\{\sl{\bf Corollary 3.2.3}$\;\;$ Let $v$ be a flow on a compact metric space $(X,d)$ and let $\a$ be a cocycle on it. If $\g(T)$ is a trajectory with time interval $T\le T_0$, then there is a constant $C_{T_0}$ depending on $T_0$ and $\{M_j\}$ such that 
$$
\Big|\int_{\g(T)}\a\Big|\le C_{T_0}
$$
}
  \\
Proof$\;\;$ For any $\g(T)$ with time interval $T\le T_0$, we have 
$$
l(\g(t))=\int_{\g(T)}|\dot{\g}(t)|\;dt \le T_0\max{|\dot{\g}(t)|},
$$
hence by applying proposition 3.2.2, we get the conclusion. \qed\\
  \\
{\bf  3.3 Integration of cocycles in a gradient-like flow}\\
  \\
\indent
In this section we assume that $v$ is a gradient-like flow with finitely many fixed points $\{p_j;j=1,2,\cdots,n\}$. Let $\a$ be a cocycle on $X$. Then there is an associated $r$-covering $\mathcal{U}(r)$ such that when restricted to the closure of each open ball in $\mathcal{U}(r)$, $\a$ is a coboundary. Let $r$ be small enough such that the Lyapunov function $g(x)$ corresponding to the flow $v$ has mutually non-intersecting images $g(B_r(p_j))$ for $j=1,2,\cdots,n$.\\
\indent
Define $\d_0(r)=\min_{1\le j\le n}\d^\pm_{r,\frac{r}{4}}(p_j)$. Since $g(x)$ is uniformly continuous, there exists $s_0\le \frac{r}{4}$ satisfying $\osc_{B_{s_0}}{g(x)}<\frac{1}{3}\d_0(r)$. By our choice and in view of the analysis of the local structure of gradient-like flows (see section 2.2), it is easy to see that the following inequality holds:
\begin{equation}
3\max_x\osc_{y\in B_{s_0}(x)}g(y)<\d_0(r)\le\d^\pm_{r,\frac{r}{4}}(p_j)\le \d^\pm_{r,s_0}(p_j)\tag{$3.3.1$}
\end{equation}
for $ j=1,2,\cdots,n$.\\
  \\
{\sl{\bf Proposition 3.3.1}$\;\;$ Let $v$ be a gradient-like flow with finitely many fixed points. $\a$ is a cocycle with an associated $r$-covering. Let $\mathcal{U}(s)$ be an $s$-covering with $s<s_0<\frac{r}{4}$ for some $s_0$ satisfying (3.3.1). Then there exist $M>0$ and $T(s_0)>0$ depending on $r$ but not on $s$ such that for any $(\mathcal{U}(s), T(s_0))$-chain $\tilde{\g}$, we have
\begin{equation}
\Big|\int_{\tilde{\g}}\a\Big|\le M \tag{$3.3.2$}
\end{equation}
}
 \\
Proof$\;\;$ Firstly we will prove that given any $ s\le r$, for any trajectory $\g\in\G_1(s)$ with domain $[a,b]$, there exists a constant $T(s)$ satisfying $b-a\le T(s).$\\
\indent
Let $x\in\partial B_s(p_j)$, then the trajectory $x\cdot t$ or $x\cdot(-t)$ for $t\ge 0$ will flow into 
some different fixed point of $v$. (There does not exist any trajectory joining one point to itself.) 
Without loss of generality, we assume that $(x\cdot t)(t\ge 0)$ flows into the point $p_k$ for $k\neq j$.
 (It is possible that the trajectory may pass through some ball $B_s(p_l)$, but this does not influence 
 the result below.) Hence there exists a $T_x>0$ such that $x\cdot T_x\cap B_{\frac{s}{2}}(p_k)
 \neq\emptyset$. By the continuity of the flow, there exists a small closed ball $D_{s_x}(x)$ which 
 satisfies that for any point $y\in D_{s_x}(x)\cap\pat B_s(p_j)$, the trajectory 
 $[y, y\cdot T_x]\cap\pat B_{\frac{3s}{4}}(p_k)\neq\emptyset$. This shows that for any trajectories starting from $D_{s_x}(x)\cap \pat B_s(p_j)$ their time intervals are not greater than $T_x$, and covering the sphere $\pat B_s(p_j)$ by finitely such small closed balls, then it is easy to  see that any trajectory in $\G_1(s)$ starting from $\pat B_s(p_j)$ has time interval not greater than a constant $ T_j(s)$. Let $T(s)=\max_{1\le j\le n}{T_j(s)}$, then all the trajectories in $\G_1(s)$ have time interval not greater than $T(s)$.\\
\indent
Let $\g=\{x_0,y_0,x_1,y_1,\cdots,x_{k-1},y_{k-1},x_k|t_0,\tau_0,\cdots,t_{k-1},\tau_{k-1}\}$ be a chain. It is called an $(r,s_0)$-chain if it satisfies the following conditions:\\
\indent
(1) $x_j\in\pat B_r(p_{l_j})$ for $j=0,1,\cdots,k-1$ and $x_k\in B_r(p_{l_k})$;\\
\indent
(2) $y_j\cdot\tau_j\in\pat B_r(p_{l_{j+1}})$ for $j=0,1,\cdots,k-1$;\\
\indent
(3) the point pair $(x_j\cdot t_j,y_j)\in B_s(m_j)$, for $j=0,1,\cdots,k-1$, where those $B_s(m_j)'s$ are elements in the $s$-covering $\mathcal{U}(s)$ of $X$.\\
\indent
(4) $B_s(m_j),j=0,1,\cdots,k-1$ and $B_s(p_{l_j}),j=0,1,\cdots,k$, are $2k-1$ open balls and any two of them do not intersect.\\
\indent
(5) $t_j,\tau_j\le T(s_0)$, where $T(s_0)$ is an upper bound for the time interval of all the trajectories in $\G_1(s_0)$.\\
\indent
Now compute the integration of $\a$ along the $(r,s_0)$-chain $\g$.\\
\begin{equation}
\aligned
\Big|\int_{\g}\a\Big|&\le \sum^{k-1}_{j=0}\Big|\int_{[x_j,x_j\cdot t_j]}\a\Big|+\sum^{k-1}_{j=0}\Big|\int_{[y_j,y_j\cdot\tau_j]}\a\Big|\notag\\
&+\sum^{k-1}_{j=0}\Big|\int_{[x_j\cdot t_j,y_j]}\a\Big|+\sum^{k-1}_{j=0}\osc_{x\in B_r(p_{l_{j+1}})}{\b_{p_{l_{j+1}}}}\notag\\
&=I+I\!\!I+I\!\!I\!\!I+I\!\!V
\endaligned\tag{$3.3.3$}
\end{equation}
Since $0<t_j,\tau_j\le T(s_0)$, applying corollary 3.2.3, we obtain 
$$
I+I\!\!I\le 2kC_1
$$
where $C_1$ is a constant depending only on $\max{|\dot{v}(t)|},\a$ and $r$.\\
\indent
Connect the point pair $(x_j\cdot t_j,y_j)$ with a line segment, then its length is at most $2s$. Using proposition 3.2.2, we have 
$$
I\!\!I\!\!I\le kC_2
$$
Here $C_2$ depends only on $s_0,\a$ and $r$.\\
\indent
For $I\!\!V$,  
$$
I\!\!V\le 2kM_\a(r)
$$
Combining the above estimates, we have
\begin{equation}
|\int_{\g}\a|\le kC\tag{$3.3.4$}
\end{equation}
where $C$ depends on $r,s_0,\a$ and $\max{|\dot{v}(t)|}$.\\
\indent
Now the proof of the proposition is changed to the problem to reduce each $(\mathcal{U}(s),T(s_0))$-chain to an $(r,s_0)$-chain, while keeping the integration of $\a$ invariant.\\
\indent
Let $\tilde{\g}=\{x_0,x_1,\cdots,x_k|t_0,\cdots,t_{k-1}\}$ be a $(\mathcal{U}(s),T(s_0))$-chain. Consider the following cases.\\
\indent
(a) If the chain $\tilde{\g}$ is contained in a ball $\bar{B}_r(p_{i_0})$ for $0\le i_0\le n$, then
\begin{align*}
&|\int_{\tilde{\g}}\a|=|\b_{p_{i_0}}(x_k)-\b_{p_{i_0}}(x_0)|\\
& \le\osc_{x\in B_r(p_{i_0})}{\b_{p_{i_o}}(x)}\le 2M_\a(r)
\end{align*}
\indent
(b) There is $1\le k_0\le k$ such that the subchain $\tilde{\g}_{0k_0}=\{x_0,x_1,\cdots,x_{k_0}|t_0,\cdots,\\t_{k_0-1}\}$ is contained in $B_r(p_{i_0})$ but the point $x_{k_0}\cdot t_{k_0}$ is not in $\cup^n_{j=0}\bar{B}_r(p_j)$. Hence there is a first intersection point $\bar{x}_0$ of $[x_{k_0},x_{k_0}\cdot t_{k_0}]$ and $\pat B_r(p_{i_0})$. Let $\bar{t}_0$ satisfy $\bar{x}_0\cdot\bar{t}_0=x_{k_0}\cdot t_{k_0}$ and let $\bar{y}_0=x_{k_0+1}$. Denote the ball containing the point pair $(\bar{x}_0\cdot\bar{t}_0,\bar{y}_0)$ by $B_s(m_0)$. Since $t_{k_0}\ge T(s_0)$ and $\tilde{\g}_{0k_0}\sbs B_r(p_{i_0})$, this implies that $x_{k_0}\in B_{s_0}(p_{i_0})$, hence the trajectory $[\bar{x}_0,\bar{x}_0\cdot\bar{t}_0]$ does not intersect with $B_{s_0}(p_{i_0})$. Otherwise,$[x_{k_0},x_{k_0}\cdot t_{k_0}]$ starts from $B_{s_0}(p_{i_0})$, intersects $\pat B_r(p_{i_0})$, then goes back to  $B_{s_0}(p_{i_0})$. This conclusion contradicts (3.3.1) when we check the change of the Liapounov function $g(x)$ along $[x_{k_0},x_{k_0}\cdot t_{k_0}]$. Therefore $[\bar{x}_0,\bar{x}_0\cdot\bar{t}_0]$ is part of a trajectory in $\G_1(s_0)$, and we have $\bar{t}_0\le T(s_0)$. Since $\bar{x}_0\cdot\bar{t_0}=x_{k_0}\cdot t_{k_0}$ is not in $\bar{B}_r(p_{i_0})$, by the choice of $s$, we know that $\bar{y}_0:=x_{k_0+1}\notin B_{s_0}(p_{i_0})$. Hence $\bar{y}_0\cdot t_{k_0+1}$ meets firstly $B_{s_0}(p_{l_1})$ for some fixed point $p_{l_1}$. Let $\bar{y}_0\cdot \bar{\tau}_0$ be the last intersection point of $[\bar{y}_0,\bar{y}_0\cdot t_{k_0+1}]$ with $\pat B_r(p_{i_1})$. In the same way, we can prove that $\bar{\tau}_0\le T(s_0)$. Joining $x_0$ and $\bar{x}_0$ by a curve $l_{x_0\bar{x}_0}$, then the subchain $\tilde{\g}_{0k_0}$ is reduced to a curve $l_{x_0\bar{x}_0}$ combined with a chain $\{\bar{x}_0,\bar{y}_0|\bar{t}_0,\bar{\tau}_0\}$ and it is clear that the reduction keeps the integration of $\a$ invariant.\\
\indent
(c) There is a $2\le k_0\le k$ such that the subchain 
$\tilde{\g}_{0,k_0-1}=\{x_0,x_1,\cdots,x_{k_0-1}|t_0,\\t_1,\cdots,t_{k_0-2}\}$ and the trajectory 
$[x_{k_0-1},x_{k_0-1}\cdot t_{k_0-1}]$ is contained in $B_r(p_{i_0})$ but $x_{k_0}\notin B_r(p_{i_0})$. 
Then we take a point $\bar{x}_0\in\pat B_r(p_{i_0})\cap B_s(m_0)$, where $B_s(m_0)$ is the ball 
containing $(x_{k_0-1}\cdot t_{k_0-1},x_{k_0})$. Since $t_{k_0}\ge T(s_0)$, the trajectory 
$[x_{k_0},x_{k_0}\cdot t_{k_0}]$ will firstly meet $B_{s_0}(p_{i_1})$. We denote $x_{k_0}$ by 
$\bar{y}_0$ and represent the last intersection point of $[x_{k_0},x_{k_0}\cdot t_{k_0}]$ with 
$\pat B_r(p_{i_1})$ by $\bar{y}_0\bar{\tau}_0$ for some $\bar{\tau}_0>0$. Due to the same argument as in 
(b), we know that $\bar{\tau}_0\le T(s_0)$. So $\tilde{\g}_{0k_0}$ is reduced to the combination of a curve $l_{x_0\bar{x}_0}$ and a chain $\{\bar{x}_0,\bar{y}_0|0,\bar{\tau}_0\}$. It is clear that such inductions also keep the integration of $\a$ invariant. Note that $x_{k_0-1}$ must be in $B_{s_0}(p_{i_0})$.\\
\indent
(d) There is a subchain $\tilde{\g}_{0k_0}=\{x_0,x_1,\cdots,x_{k_0}|t_0,\cdots,t_{k_0-1}\}$ contained in 
$B_r(p_{i_0})$ and the trajectory $ [x_{k_0},x_{k_0}\cdot t_{k_0}]$ goes through $B_r(p_{i_1})$. 
Let $\bar{x}_0=\bar{y}_0$ be the first intersection point of  $[x_{k_0},x_{k_0}\cdot t_{k_0}]$ with 
$\pat B_r(p_{i_1})$. Also we can prove that $\bar{\tau}_0\le T(s_0)$ as in (b) and (c). Then the chain $\tilde{\g}_{0k_0}$ followed by $[x_{k_0},\bar{y}_0\cdot\bar{\tau}_0]$ is reduced to the combination of $l_{x_0\bar{x_0}}$ with the chain $\{\bar{x}_0,\bar{y}_0|0,\bar{\tau}_0\}$, while the integration of $\a$ is invariant under change.\\
\indent
(e) If $x_0\notin\cap^n_{j=1}B_r(p_j)$, then $x_0\cdot t_0$ will firstly meet some $B_s(p_{i_0})$. Let $x'_0=x_0\cdot t'_0$ be the last intersection point of $[x_0,x_0\cdot t_0]$ with $\pat B_r(p_{i_0})$. In the same way, we can prove $t'_0\le T(s_0)$.\\
\indent
Now using the above steps (b)-(e) continuously, any $(\mathcal{U}(s),T(s_0))$-chain 
$\tilde{\g}=\{x_0,x_1,\cdots,x_k|t_0,\cdots,t_{k-1}\}$ can be reduced to the combination of a curve 
$l_{x_0\bar{x}_0}$ ( or $[x_0,x_0\cdot t_0']$ in case (e)) with an $(r,s_0)$-chain 
$\bar{\g}=\{\bar{x}_0,\bar{y}_0,\cdots,\bar{x}_l|\bar{t}_0,\bar{\tau}_0,\cdots,\bar{\tau}_{l-1}\}$ 
if we can prove that (4) in the definition of a $(r,s_0)$-chain holds.\\
\indent
If (4) is not true, then there exists an $(r,s_0)$-cycle $\bar{\g}_c=\{\bar{x}_{i_0},\bar{y}_{i_0},\cdots,\bar{y}_{i_l},\bar{x}_{i_0}|\bar{t}_{i_0},\\ \bar{\tau}_{i_0},\cdots,\bar{\tau}_{i_l}\}$ where the point pair $(\bar{y}_{i_l}\cdot\bar{\tau}_{i_l},\bar{x}_{i_0})\in\pat B_r(p_{i_0})$. Now extend the trajectory $\bar{y}_{i_j}\cdot t(t\ge 0)$ for $j=0,\cdots,l.$ Since $\bar{\tau}_c$ is obtained by reducing the $(\mathcal{U}(s),T(s_0))$-chain, hence $\bar{y}_{i_j}\cdot t$ will go through $B_{s_0}(p_{i_{j+1}})$. We let $\bar{y}_{i_j}\cdot\tau_{i_j}(\tau_{i_j}>\bar{\tau}_{i_j})$ be a point in $B_{s_0}(p_{i_{j+1}})$. Take the inverse process corresponding to (2)-(4), then there is a trajectory starting from some point $x_{i_j-1}\in B_{s_0}(p_{l_j-1})$ to some point $x_{i_j-1}\cdot t_{i_j-1}$ in the $s$-ball containing $\bar{x}_{i_{j-1}}\cdot t_{i_{j-1}}$ and $\bar{y}_{i_j}$. Therefore frome the cycle $\bar{\g}_c$ we get a $(\mathcal{U}(s),T(s_0))$-chain $\tilde{\g}_c=\{x_{i_0},\bar{y}_{i_0},\cdots,\bar{y}_{i_l},x_{i_0}|t_{i_0},\tau_{i_0},\cdots,t_{i_l},\tau_{i_l}\}$.\\
\indent 
Consider the change of the function $g(x)$ along $\tilde{\g}_c$. On one hand, since the two ends of $\tilde{\g}_c$ are $x_{i_0}$ and $\bar{y}_{i_l}\cdot\tau_{i_l}$, by (3.3.1)
\begin{align*}
\Big|\int_{\tilde{\g}_c}\d g\Big|\le \osc_{x\in B_{s_0}(p_{i_0})}g(x)<\frac{1}{3}\d_0(r)
\end{align*}
On the other hand, each $\bar{y}_{i_j}\cdot t (0<t\le\tau_{i_j})(j=0,\cdots,l)$ goes across $B_r(p_{i_j})\backslash B_{\frac{r}{4}}(p_{i_j})$, we have 
\begin{align*}
\int_{\tilde{\g}_c}\d g\ge& (l+1)\d^+_{r,\frac{r}{4}}-2(l+1)\max_{x}{\osc_{y\in B_s(x)}g(y)}
>&\d_0(r)-\frac{2}{3}\d_0(r)=\frac{1}{3}\d_0(r)
\end{align*}
This is absurd. The contradiction shows that the chain $\bar{\g}$ we get from a $(\mathcal{U}(s),T(s_0))$-chain $\tilde{\g}$ is indeed an $(r,s_0)$-chain. Let $\bar{\g}=\{\bar{x}_0,\bar{y}_0,\cdots,\bar{x}_l|\bar{t}_0,\bar{\tau}_0,\cdots,\bar{\tau}_{l-1}\}$. Since $v$ is gradient-like, the index $l$ is less than $n$. Therefore applying the estimate (3.3.4) and the fact that the reduction from $\tilde{\g}$ to $\bar{\g}$ keeps the integration of $\a$, we have 
$$
\Big|\int_{\tilde{\g}}\a\Big|=\Big|\int_{\bar{\g}}\a\Big|\le nC.
$$
where $C$ depends on $r, s_0, \a $ and $\max{|\dot{v}(t)|}$. Proposition 3.3.1 now is proved. \qed\\
 \\
{\bf 3.4 Flows carrying a cocycle $\a$}\\
  \\
\indent 
Consider a fixed cocycle $\a$ on $X$. With respect to $\a$ there is an $r$-covering for small $r$ such that the restriction of $\a$ to the closure of every $r$-ball is a coboundary. Let $v$ be a flow on $X$ having $n$ fixed points $\{p_i\}^n_{i=1}$. Consider the flow $v$ restricted to the punctured space $X\backslash\cup^n_{i=1}B_r(p_i)$. Recall that in the end of section 2.2, we have defined three sets $\G_1(r),\G_2(r)$ and $\G_3(r)$. Now we have the following definitions.\\
  \\
{\sl {\bf Definition 3.4.1}$\;\;$ The flow $v$ defined on the compact metric space $(X,d)$ is said to be an $\a$-flow, if there exist a cocycle $\a$ , a small $r>0$ and a $T_0>0$ such that for some $\r>0$ and $0\le\la<1$, the following conditions are satisfied:\\
\indent
(1)$\;\; v$ is gradient-like in $B_r(p_i)$\\
\indent
(2)$\;\;\max_{(x,y)\in \bar{B}_r(p_i)\times\bar{B}_r(p_i)}{|I_\a(x,y)|}\le \la\r$, for $1\le i\le n$.\\
\indent
(3)$\;\;$ for any trajectory $\g\in\G_1(r)$, 
$$
\int_\g\a\ge\r
$$
\indent
(4)$\;\;$ if $\g(T_0)$ denotes any sub-trajectory of $\g\in\G_2(r)\cup\G_3(r)$ with time interval $T_0$, then
$$
\int_{\g(T_0)}\a\ge \r
$$
} 
  \\
Similarly, if we replace the fixed points by connected nonwandering sets of $v$ and replace the set $\G_i(r)$ by $\G_i^A(r)$ which is also defined at the end of section 2.2, then we can get the definition of generalized $\a$-flows with respect to the nonwandering set $A$.\\
  \\
{\sl {\bf Definition 3.4.2}$\;\;$ The flow $v$ defined on the compact metric space $(X,d)$ is said to be a generalized $\a$-flow with respect to the nonwandering set $A=\{A_n,\cdots,A_1\}$, if there exist a cocycle $\a$ , a small $r>0$ and a $T_0>0$ such that for some $\r>0$ and $0\le\la<1$, the following conditions are satisfied:\\
\indent
(1)There is a Lyapunov function $g_i(x)$ defined in $B_r(A_i)$ for $i=1,\cdots,n$, such that $g_i(x)$ is constant on $A_i$.\\
\indent
(2)$\;\;\max_{(x,y)\in \bar{B}_r(A_i)\times\bar{B}_r(A_i)}{|I_\a(x,y)|}\le \la\r$, for $1\le i\le n$.\\
\indent
(3)$\;\;$ for any trajectory $\g\in\G_1^A(r)$, 
$$
\int_\g\a\ge\r
$$
\indent
(4)$\;\;$ if $\g(T_0)$ is any sub-trajectory of $\g\in\G_2^A(r)\cup\G_3^A(r)$ having time interval $T_0$, then 
$$
\int_{\g(T_0)}\a\ge \r
$$
} 
  \\
{\sl {\bf Definition 3.4.3}$\;\;$ A flow $v$ is called carrying a cocyle $\a$ if $v$ is an $\a$-flow or a generalized $\a$-flow.}\\
 \\
\indent 
We will give examples and some propositions to show that the (generalized) $\a$-flows include many interesting and important flows. Let us give some examples first.\\
  \\
\noindent 
{\bf Example 3.4.1}$\;\;$
Let $v$ be the gradient flow generated by a Morse function $f$ on a smooth manifold $X$. We say that the flow $v$ is a $df$-flow. In fact, in this case every function $\b_j$ defined on the closed ball $B_j$ in the $r$-covering can be taken as the restriction of $f$ to $B_r(p)$. Now let $S(f)$ be the critical point set and $p_i$ denote a critical point. Let $\r=\frac{1}{2}\min{|f(p_i)-f(p_j)|}$, where the minimum is taken over all pairs $(p_i,p_j)\in S_i(f)\times S_j(f), i\neq j$,
 and choose $r$ small enough such that $\max_{p_i\in S(f)}{\osc_{B_r(p_i)}{|f|}}\le \frac{\r}{2}$, then 
 the flow $v$ satisfies conditions (1) and (2). To prove (3), let $\g\in\G_1(r)$ be a segment of a trajectory joining $B_r(p_i)$ and $B_r(p_j)$. Since $v$ is 
 a gradient-flow, the $p_i$ and $p_j$ must be two critical points with different index $i$ and $j$. Then
\begin{align*}
&\int_\g df=f(\g(b))-f(\g(a)) \\
&=f(p_i)-f(p_j)+f(\g(b))-f(p_i)+f(p_j)-f(\g(a))\\
&\ge f(p_i)-f(p_j)-2\max_{p_i\in S(f)}{\osc_{B_r(p_i)}{|f|}}\\
&\ge 2\r-\r=\r,
\end{align*}
hence (3) is proved. Since $v$ is a gradient flow, $\G_2(r)\cup \G_3(r)$ is empty. Hence the condition (4) is trivial.\\
  \\
\noindent
{\bf Example 3.4.2}$\;\;$
Let $v$ be the flow generated by the vector field $V$ dual to the closed Morse 1-form $\o$. Take an $r$-covering of $X$ such that $\o$ is exact when restricted to the closed ball $\bar{B}_r(p)$, for instance, let $\o=df_p|_{\bar{B}_r(p)}$. By the preceding discussion, we know that $I_\o(x,y)$ and $\int_\g\o$ are independent of the choice of the $r$-covering and the related functions $\{f_p\}$. Therefore 
\begin{align*}
&\max_{(x,y)\in \bar{B}_r(p_i)\times\bar{B}_r(p_i)}{|I_\o(x,y)|}
=|\hbox{osc}_{x\in\bar{B}_r(p_i)}f_{p_i}(x)|\\
&\le 2r\max_{x\in\bar{B}_r(p_i)}|df_{p_i}|
=2r\max_{x\in\bar{B}_r(p_i)}{|\o(x)|}
\tag{$3.4.1$}
\end{align*}
On the other hand, there exists a $r_0>0$ such that a good coordinate system can be chosen in 
$B_{r_0}(p_i)$ for any $i=1,\cdots,n$, to make the dual vector field have the form 
$(-x_1,\cdots,-x_s,x_{s+1},\cdots \,x_m)$. Furthermore, $r_0$ can be chosen small enough such that 
$$
r_0<\min_{i\neq j}d(B_{r_0}(p_i), B_{r_0}(p_j)),
$$
where $p_i\in S_i(\o)$.\\
\indent
Let $r<\frac{r_0}{2}$, then by analyzing the explicit expression of the trajectory in 
$B_{r_0}(p_i), i=1,\cdots,n$, it is easy to see that for any $\g\in\G_1(r)$,
$$
\int^b_a|\dot{\g}(t)|\;dt\ge\frac{r_0}{2}.
$$
\indent
Using the above inequality, for any $\g\in \G_1(r)$, we have 
\begin{equation}
\aligned
\int_\g\o=&\sum_j f_{x_{i_j}}(\g(t_{j+1}))-f_{x_{i_j}}(\g(t_j))\notag\\
=&\sum_j\int^{t_{j+1}}_{t_j}|df_{x_{i_j}}|^2=\int^b_a|\o(\g(t))|^2\;dt\notag\\
\ge &\min_{x\notin\cup_i B_r(p_i)}|\o(x)|\cdot\int^b_a|\dot{\g}(t)|\;dt\notag\\
\ge &\min_{x\notin\cup_i B_r(p_i)}|\o(x)|\cdot C_{r_0}
\endaligned\tag{$3.4.2$}
\end{equation}
where $C_{r_0}=\frac{r_0}{2}$\\
Therefore if we can choose $r,\r$ satisfying 
\begin{equation}
4r\max_{x\in\cup_i\bar{B}_r(p_i)}|\o(x)|\le \r\le C_{r_0}\min_{x\notin\cup_i B_r(p_i)}|\o(x)|\tag{$3.4.3$}
\end{equation}
then conditions (2) and (3) are satisfied ((1) is obvious). Let $|\cdot|_0$ be the Euclidean metric, then there exist constants $C_0$ and $C_1$ such that for any vector field $V$ near any point $x\in X$,  
$$
C_0|V|_0\le |V|\le C_1|V|_0
$$
Then 
$$
\min_{x\notin\cup_i B_r(p_i)}|\o(x)|\ge C_0 r
$$
and 
$$
\max_{x\in\cup_i\bar{B}_r(p_i)}|\o(x)|\le C_1 r
$$
Hence (3.4.3) holds if we can choose $\r$ satisfying
\begin{equation}
4C_1 r^2\le \r\le C_{r_0}C_0 r\tag{$3.4.4$}
\end{equation}
This can be done by choosing $r$ small enough. Now we can choose $r$ and $\r$ such that the conditions (2) and (3) are true.\\
\indent
For the condition (4), we can choose $T_0$ large enough such that for fixed $r$ and $\r$,
$$
\r\le  \min_{x\notin\cup_i B_r(p_i)}|\o(x)|^2 T_0\le \int_{\g(T_0)}|\o(\g(t))|^2\;dt=\int_{\g(T_0)}\o
$$
The above argument shows that the flow $v$ satisfies our definition, i.e., $v$ is a $\o$-flow, which fits with the fact that $v$ is generated by the vectorfield dual to $\o$.\\
  \\
\noindent
{\bf Remark}$\;\;$ Example 3.4.1 can also be treated as in the example 3.4.2, but in this case every trajectory is in the set $\G_1(r)$. In the case of Example 3.4.2, it is possible for some closed 1-forms that the flows they generate have trajectories in $\G_2(r)\cup\G_3(r)$.\\
  \\
\noindent
{\bf Example 3.4.3}$\;\;$
We let $\o$ be a closed 1-form as in the example 3.4.2. Now we consider a family of flows depending on the parameter $\ape,0\le \ape<1$:
$$
\dot{v}_\ape(t)=V(v(t))+\ape U(v(t)),
$$
where $V(x)$ is the dual vector field of $\o$ and $U(x)$ is a vector field vanishing on a small neighborhood 
of each zero point of $\o$ and satisfying the norm inequality:
$$
|U(x)|\le |V(x)|=|\o(x)|.
$$
We claim that the flow $v_\ape$ is also an $\o$-flow. This can be seen as follows. Firstly (3.4.1) holds. Secondly we study the integration of $\o$ along  $\g_\ape\in\G^\ape_1(r)$, we have 
$$
\aligned
\int_{\g_\ape}\o=&\sum_j f_{x_{i_j}}(\g_\ape(t_{j+1}))-f_{x_{i_j}}(\g_\ape(t_j))\\
=&\sum_j \int^{t_{j+1}}_{t_j}df_{x_{i_j}}(\dot{\g}_\ape(t_j))\\
=&\sum_j\int^{t_{j+1}}_{t_j}\langle\o(\g_\ape(t)), V(\g_\ape(t))+\ape U(\g_\ape(t))\rangle dt\\
\ge& (1-\ape)\int^b_a|\o(\g_\ape(t))^2 dt
\endaligned
$$
Since $U(x)$ vanishes in a small neighborhood of each zero point, the argument behind (3.4.2) also holds 
in this case. Hence (3.4.3) should be modified as
$$
4r\max_{x\in\cup_i\bar{B}_r(p_i)}|\o(x)|\le \r\le (1-\ape)C_{r_0}\min_{x\notin\cup_i B_r(p_i)}|\o(x)|
$$
Thus we can choose $r,\r$ small enough such that the conditions (2) and (3) in the definition are true.\\
\indent
For (4), we should choose $T_0$ large enough such that 
$$
\r\le (1-\ape) \min_{x\notin\cup_i B_r(p_i)}|\o(x)|^2 T_0\le\int_{\g_\ape(T_0)}\o
$$
Therefore, if we choose $r,\r,T$ suitably, the conditions (1)-(4) are also true for the flow $v_\ape$.\\
 \\
{\bf Example 3.4.4}$\;\;$ This example will provide a method for changing a gradient flow of a Morse function to an $\a$-flow with $\a$ being a nontrivial cocycle.\\
\indent
It is well-known that in any closed smooth manifold $M_0$ there exists a Smale function with exactly one local maximum point and one local minimum point. Hence the gradient flow of such a Smale function has only one minimal attractor and one minimal repeller. Now take small $m$-disc $D_0^m$ and $D_1^m$ whose center points are the minimum point and the maximum point respectively. Now cutting the two discs and connecting the two holes by a ``bridge'' $S^{m-1}\times D^1$, we get a new manifold $M_1$. Take the flow on $S^{m-1}\times D^1$ to be the one generated by the vector field $(0,\pm\frac{d}{dt})$ and let the flow go from the part of the original minimum point to the part of the original maximum point. After connecting the vector fields smoothly in the connecting collar, we get an $\a$-flow on $M_1$ (This fact can be shown in the same way as in Example 3.4.9.). It is obvious that this flow on $M_1$ has no attractor or repeller and the set $\G_2(r)\cup\G_3(r)\neq\emptyset$. \\
\indent
If a gradient flow of a Morse function has more than two minimum or maximum points, we can proceed as above 
and connect some minima and maxima pairwise
to construct the $\a$-flow with nontrivial cocycle $\a$ and having attractors and repellers.\\
  \\
{\bf Example 3.4.5 (Flows carrying a cohomology class)}$\;\;$ We consider an extreme case in the 
definition of an $\a$-flow $v$, i.e., that $v$ has no fixed point in $X$. In this case, the set 
$\G_1(r)\cup\G_2(r)=\emptyset$ and the only condition that makes $v$ an $\a$-flow is that there 
exist constants $\r>0$ and $T_0>0$ such that for any trajectory $\g(T_0)$ with time interval $T_0$, we have 
\begin{align*}
\int_{\g(T_0)}\a\ge \r
\end{align*}
\indent
Now the $\a$-flow $v$ becomes a so called ``flow carrying a cohomology class'' as introduced by 
R.C.Churchill in his paper [8]. The reason that the flow is called ``carrying a cohomology class'' is 
that the above condition is independent of the choice of the representative in the cohomology class 
$[\a]$. In fact, if $\a_1\in[\a]$ is another cocycle, then there exists a coboundary 
$\b\in\mathcal{H}^1(X,{\Bbb R})$, such that 
\begin{align*}
\a_1-\a\simeq \b
\end{align*}
and so 
\begin{align*}
\int_{\g(kT_0)}\a_1&=\int_{\g(kT_0)}\a+\b(e(\g(kT_0)))-\b(s(\g(kT_0)))\\
&\ge k\r-2M_\b
\end{align*}
where $s(\g)$ and $e(\g)$ are the start point and the end point of the trajectory $\g$, and $M_\b$ is the bound of $\b$. Hence if we choose $k>[\frac{2M_\b+\r}{\r}]+1$, then we have 
\begin{align*}
\int_{\g(kT_0)}\a_1\ge \r.
\end{align*}
\indent
The existence of a flow carrying a cohomology class in a manifold will induce the vanishing theorem for the Novikov numbers. This result will be given in section 4.5.\\
\indent
The flow in the following example is a simple example of such an $\a$-flow which also carries a 
cohomology class. It was first explained in this way in [8].\\
  \\
{\bf Example 3.4.6}$\;\;$ Let $S^1$ be the unit circle. We will define a nontrivial cocycle $\a$ on $S^1$ as follows. Define two open sets in $S^1$:
\begin{align*}
U_1&=\Big\{e^{i\ta}\;|\;-\frac{5\pi}{8}<\ta<\frac{5\pi}{8}\Big\}\\
U_2&=\Big\{e^{i\ta}\;|\;\frac{3\pi}{8}<\ta<\frac{13\pi}{8}\Big\}
\end{align*}
Then $U_1\cup U_2$ is an open covering of $S^1$. Define two bounded functions in $U_1$ and $U_2$. Let 
$\b_2=0$ in $U_2$ and let 
\begin{align*}
\b_1(e^{i\ta})=\Big\{
\begin{array}{ll}
1 & \;\ta\in [0,\frac{5\pi}{8}]\\
0 & \;\ta\in [-\frac{5\pi}{8},0)
\end{array}
\end{align*}
Now we can define a nontrivial cocycle $\a:S^1\times S^1\to {\Bbb R}$,
\begin{align*}
\a(e^{i\ta_1},e^{i\ta_2})=\Big\{
\begin{array}{ll}
\b_i(e^{i\ta_2})-\b_i(e^{i\ta_1}) &\;\hbox{if}\;\; |\ta_1-\ta_2|<\frac{\pi}{4},\;i=1,2\\
0 &\;\hbox{if}\;\; |\ta_1-\ta_2|\ge \frac{\pi}{4}
\end{array}
\end{align*}
\indent
Let $\pi_i:T^m=S^1\times\cdots\times S^1\to S^1$ be the projection map which projects from the $m$-dimensional torus to the $i$-th factor space. Define $\a_i(x,y):T^m\times T^m\to {\Bbb R}$ to be the pull-back map $(\pi_i\times \pi_i)^*\a$. It is easy to see that $[\a_i](i=1,\cdots,m)$ is a linearly independent basis of the bounded Alexander-Spanier cohomology group $\mathcal{H}^1(X,{\Bbb R})$\\ 
\indent
Let $\o=(\o_1,\cdots,\o_m)$ be a constant irrational vector field in ${\Bbb R}^m$. Here ``irrational'' means that $\o_i(i=1,\cdots,m)$ is an irrational number. Then this vector field can be projected to the torus $T^m$ such that $\o$ is a vector field on $T^m$. Consider the irrational flow generated by $\o$,
$$
\dot{v}=\o(v)
$$
Taking any point $x$ on $T^m$ and integrating $\a_i$ along $[x,x\cdot t]$, we can get
$$
\int_{[x,x\cdot t]}\a_i=[\frac{\o_i t}{2\pi}].
$$
This shows that the irrational flow is an $\a_i(-\a_i)$-flow for $i=1,\cdots,m$ , if $\o_i>0(\o_i<0)$. In particular, it is a flow carrying cohomology class $[\a_i](-[\a_i])$ if $\o_i>0(\o_i<0)$.\\
  \\
{\bf Example 3.4.7}$\;\;$ In Example 3.4.1 and Example 3.4.2, we have considered the flow generated by a dual vector field of a Morse function or a closed Morse 1-form. In those cases, the fixed points of the flows are non-degenerate fixed points. In some cases, those flows generated by a function or a closed 1-form may have Bott type nondegenerate critical manifolds. For instance, let $Z$ be a non-degenerate critical manifold of a function $f$. This means the following condition holds:\\
\indent
(1) grad$f(x)=0,\,\forall x\in Z$.\\
\indent
(2) At any point $x\in Z$, there is a decomposition $T_x(M)=T_x^\perp(Z)\otimes T_x(Z)$ and $d^2f(x)$ as bilinear form is nondegenrate along the vertical tangent subspace $T^\perp_x(Z)$.\\
\indent
Since locally a closed 1-form $\o$ is the differential of a smooth function, hence the Bott type nondegenerate zero set of the closed 1-form $\o$ can be defined as above.\\
\indent
Therefore using the same argument as in Example 3.4.2, it can be proved that such a flow generated by $\o$ is a generalized $\o$-flow with respect to all the Bott type nondegenerate zero sets.\\
 \\
{\bf Example 3.4.8}$\;\;$ In the paper [27], S.Smale introduced the Morse-Smale flow. We give the definition below.\\
\indent
A flow $v$ on a manifold $M$ is called a Morse-Smale flow if it satisfies:\\
\indent
(1) The chain recurrent set of $v$ consists of a finite number of hyperbolic closed orbits and hyperbolic fixed points. ( The reader can see the meaning of ``hyperbolic'' in section 4.3)\\
\indent
(2) The unstable manifold of any closed orbit or fixed point has transversal intersection with the stable manifold of any closed orbit or fixed point.\\
\indent
It was proved that such flows have ``global'' gradient-like structures and have a Morse decomposition which induces the Morse inequalities.\\
\indent
However, in some cases, although the nonwandering set of the flow contains only the hyperbolic periodic 
orbits and the hyperbolic fixed points, the flow is not a Morse-Smale flow because of the existence of 
``cycles'' which consist of some orbits ``connecting'' different invariant sets and form a closed curve. 
We can give a definition of such flows when restricted to the category of flows carrying a cocycle.\\
 \\
{\bf $\a$-Morse-Smale flow}$\;\;$ Let $v$ be a generalized $\a$-flow with respect to the nonwandering set $A=\{A_n,\cdots,A_1\}$. If $A$ contains only the hyperbolic orbits or hyperbolic fixed points, then $v$ is called an $\a$-Morse-Smale flow.\\
The following example provides a concrete $\a$-Morse-Smale flow.\\
  \\
{\bf Example 3.4.9}$\;\;$ Let $W$ be a compact oriented connected manifold with boundary $\pat W=\pat_+ W\cup \pat_- W$. $W$ is called a flow manifold if it satisfies the following two conditions:\\ 
\indent
(1) $\pat_+ W\cap \pat_- W=\emptyset$;\\
\indent
(2) The Euler characteristic number $\chi(W)=\chi(\pat_- W)$.\\
\indent
In his paper [1] and [2], D.Asimov has proved the facts that if $W$ is a flow manifold such that $\dim(W)\neq 3$ and $W$ is not a M\"obius band, then $W$ has a round handle decomposition and a manifold with a round handle decomposition has a non singular Morse-Smale vector field, i.e., the flow generated by this vector field  has only hyperbolic periodic orbits and at the boundary $\pat_- W(\pat_+W)$, the vector field points outward (inward). Actually, he obtained the following result (see [15]):\\
\indent
Let $f:W\to [a,b]$ be a Morse function with two critical points $p$ and $q$ of index $k$ and $k+1$ respectively, such that the unstable manifold of $q$ does not intersect the stable manifold $p$, then there is a non-singular Morse-Smale vector field $X$ on $W$ satisfying,\\
\indent
(1) $X=-\na f$ on a neighborhood of $\pat W$, and\\
\indent
(2) the flow of $X$ has exactly one closed orbit and this orbit has index $k$ and is untwisted.\\
\indent
Now we can use this fact to construct many $\a$-Morse-Smale flows. A simple example is that we can 
reverse the flow and then double the manifold to get a closed manifold $2W$ and a flow on it such that 
the nonwandering set of the flow are the untwisted hyperbolic periodic orbits with index $k$ and 
index $n-1-k$. This is actually an $\a$-flow. The cocycle $\a$ can be defined as follows. Take a metric $d$ on $2W$, and define two open sets on $2W$ to be 
\begin{align*}
&U_1=\{x\in 2W\;|\;d(x,\pat_-W)\le 2\ape\}\\
&U_2=\{x\in 2W\;|\;d(x, \pat_-W)\ge \ape\}.
\end{align*}
\indent
Then $U_1\cup U_2$ is a covering of $2W$. Denote the two half parts of $2W$ by $W_1$ and $W_2$. Define $\b_2(x)=0$ on $U_2$ and define $\b_1(x)$ in $U_1$ to be 
\begin{align*}
\b_1(x)=\Big\{
\begin{array}{ll}
0 & \;\hbox{if}\; d(x,W_1)=0\\
1 & \;\hbox{if}\; d(x,W_1)>0
\end{array}
\end{align*}
Now the cocycle $\a:2W\times 2W\longrightarrow {\Bbb R}$ is defined as 
\begin{align*}
\a(x,y)=\Big\{
\begin{array}{ll}
\b_i(y)-\b_i(x) &\;\hbox{if}\;d(x,y)<\ape,\;\hbox{for}\;i=1,2\\
0 &\;\hbox{if}\;d(x,y)\ge\ape
\end{array}
\end{align*}\\
\indent
In the same way we can piece together many flow manifolds with boundary to get $\a$-Morse-Smale flows that possesses hyperbolic periodic orbits and hyperbolic fixed points with different indexes.\\
   \\
{\sl {\bf Theorem 3.4.1}$\;\;$ Let $v$ be an $\a$-flow on the compact metric space $(X,d)$. If $\a$ is 
a trivial cocycle, then $v$ is a gradient-like flow. Conversely, if $v$ is a gradient-like flow, then 
$v$ is a $\d g$-flow with $g$ being the associated Lyapunov function of the flow $v$.\\}
 \\
Proof$\;\;$ Firstly we will show that if $v$ is an $\a$-flow with $\a$ being a trivial cocycle, then 
$\G_2(r)\cup\G_3(r)$ is empty. In fact, if we assume that there exists a trajectory 
$\g\in\G_2(r)\cup\G_3(r)$, then for large $T$ the integral $\int_{\g(T)}\a$ can be larger than 
any pre-given number. Now using a shortest curve $K$ joining the two ends of $\g(T)$, then 
$\int_{\g(T)\ast K}\a=0$, since $\a$ is a trivial cocycle. This gives $\int_{\g(T)}\a=|\int_K\a|$, 
however, this is not possible, since $|\int_K\a|$ is bounded for any finitely long curve $K$.\\
\indent
Secondly, we will show that there is no trajectory in $\G_1(r)$ forming a "cycle". If this case happened, 
then there would exist a $(\mathcal{U}(r), T)$-cycle $\tilde{\g}=\{p'_{i_0},\cdots,p'_{i_k}=p'_{i_0}|t_0,\cdots,t_{k-1}\}$ 
where the point pair $(p'_{i_j}\cdot t_j, p'_{i_{j+1}})\sbs B_r(p_{i_j})$ for 
$j=0,\cdots,k-1$. Since $\a$ is a trivial cocycle, we get $\int_{\tilde{\g}}\a=0$. On the other hand, we have 
$$
\aligned
\int_{\tilde{\g}}\a\ge &\sum^{k-1}_{j=0}\int_{[p'_{i_j},p'_{i_j}\cdot t_j]}\a-\sum^{k-1}_{j=0}\osc_{B_r(p_{i_j})}\b_{i_j}\\
\ge &k\r-\la k\r>0
\endaligned
$$
This is a contradiction. Therefore this result with the first one implies that for any $x\in X$, its 
limit sets $\o(x)$ and $\o^*(x)$ are in $\cup^n_{i=1}B_r(p_i)$, furthermore in $\{p_i,\;i=1,\cdots,n\}$
 in view of (1) of the definition 3.4.1. Hence $v$ is a gradient-like flow.\\
\indent
For the second conclusion, since each gradient-like flow $v$ has an associated Lypunov function $g$, if we define $\a:X^2\to{\Bbb R}$ to be 
$$
\a(x,y)=g(y)-g(x),\;\forall (x,y)\in X^2
$$
then it is easy to see that $v$ is a $\d g$-flow.\qed\\
  \\
{\sl {\bf Theorem 3.4.2}$\;\;$ Let $v$ be an $\a $-flow on the compact metric space $(X,d)$. Let $\{s_i\}\to 0$ as $i\to\infty$. If $\a$ is a nontrivial cocycle and for any $M>0$ and $T>0$ there is a $(\mathcal{U}(s_i),T)$-chain $\tilde{\g}_i$ for each $s_i$ such that 
\begin{equation}
\Big|\int_{\tilde{\g}}\a\Big|\ge M\tag{$3.4.5$}
\end{equation}
then $v$ is not a gradient-like flow. Furthermore, if it is known that for any $M>0$ there is a trajectory $\tilde{\g}$ satisfying (3.4.5) or there is an oriented cycle $\tilde{\g}_0$ consisting of some orbits joining fixed points where the direction is determined by the forward direction of the flow $v$ such that 
\begin{equation}
\Big|\int_{\tilde{\g}}\a\Big|\ge 1\tag{$3.4.6$}
\end{equation}
then $v$ is not a gradient-like flow.
}\\
  \\
\indent
Proof$\;\;$ The first conclusion is a direct corollary of proposition 3.3.1. and the third one is obvious. We only consider the second case. From proposition 3.2.3, if the time intervals of the trajectories have an upper bound, then for any cocycles on $X$ the absolute value of the integration, $|\int_{\tilde{\g}}\a|$ has a uniform bound. Therefore the trajectory $\tilde{\g}$ with respect to the arbitrary large $M$ has arbitrary large time interval and it is the chain needed for the hypothesis in the first conclusion.\qed\\
  \\
{\sl {\bf Theorem 3.4.3}$\;\;$ Let $v$ be a $\a$-flow with a non-trivial cocycle $\a$ on a compact metric space
 $(X,d)$. If $v$ is not a gradient-like flow, then for any small sequence $\{s_i\}$ which tends to zero 
 as $i\to\infty$, and for any $M>0, T>0$ there is a $(\mathcal{U}(s_i),T)$-chain $\tilde{\g}_i$ for each 
 $s_i$ such that 
\begin{equation}
\Big|\int_{\tilde{\g}}\a\Big|\ge M\tag{$3.4.7$}
\end{equation}
}\\
 \\
\noindent
Proof$\;\;$ If $v$ is not a gradient-like flow, then there is a non fixed point $x_0$ in the chain 
recurrent set of $v$. Therefore for any $s$ and $T>0$ there is a $(\mathcal{U}(s),T)$-chain 
$\tilde{\g}=\{x_0,\cdots,x_k=x_0|t_0,\cdots,t_{k-1}\}$. \\
Now we consider two possibilities:\\
\indent
(1) If $x_0\cdot t\notin \cup^n_{i=1}B_r(p_i)$\\
\indent
In this case, we have the integration:
$$
\int_{[x,x\cdot t]}\a\ge [\frac{t}{T_0}]\r
$$
Since $x_0$ is a point in the chain recurrent set of $v$, there is a $(\mathcal{U}(s),T)$-chain 
$\tilde{\g}=\{x_0,\cdots,x_k=x_0|t_0,\cdots,t_{k-1}\}$, and we can modify the chain $\tilde{\g}$ to get a 
$(\mathcal{U}(s),kT)$-chain $\tilde{\g}_1=\{x_0,x_0|kT\}$. We can let $\tilde{\g}$ to be a 
$(\mathcal{U}(s_0),T)$-chain for small $s_0<s$ such that there exists a $x'_0$ in the ball $B_s(m_{i_1})$ 
containing the point pair $(x_0\cdot T,x_1)$ satisfying $x'_0\cdot(-T)\in B_s(p_{i_0})$ and 
$x'_0\cdot T\in B_s(m_{i_2})$. The point is that a $(\mathcal{U}(s_0),T)$-chain for small $s_0$ can 
provide a chain $(\mathcal{U}(s),T)$ for large $s$, while decreasing the number of trajectories it 
contains. Hence we can obtain a $(\mathcal{U}(s),T)$-chain $\tilde{\g}_1=\{x_0,x_0|T\}$ for any $s$ and 
$T$. The integration of $\a$ along $n\tilde{\g}_1$ has the following estimates:
$$
\int_{n\tilde{\g}_1}\a\ge n([\frac{T}{T_0}]\r-\osc_{B_s}\b_{x_0})\ge n([\frac{T}{T_0}]\r-2M_\a)
$$
Therefore if we choose $T\ge \frac{2M_\a\cdot T_0}{\r}+1$ and $n$ large enough, then (3.4.7) is true.\\
\indent
(2) If $x_0\cdot t_0\in\cup^n_{i=1}B_r(p_i)$ for some $t_0$\\
\indent
In this case we can replace $x_0$ by $x_0\cdot t_0$, since $x_0\cdot t_0$ is also a point in the chain recurrent set of the flow $v$. So we assume that $x_0\in B_r(p_{i_0})$. Using the same argument as in (1), there exists an $s_0>0$ depending on $x_0$ such that for any $s<s_0$, there is a $(\mathcal{U}(s),T)$-chain $\tilde{\g}=\{x_0,\cdots,x_k=x_0|t_0,\cdots,t_{k-1}\}$ satisfying $(x_j\cdot t_j,x_{j+1})\in B_s(m_{j+1})\sbs B_r(p_{i_{j+1}})$ for $j=0,1,\cdots,k-1$. Thus we have 
$$
\int_{n\tilde{\g}}\a\ge n(k\r-\la k\r)
$$
If $n$ is large enough, the integration will be larger than any given number. So (3.4.7) holds. \qed\\
\indent
In view of Theorem 3.4.1-3.4.3, the following corollary is obvious.\\
  \\
{\sl {\bf Corollary 3.4.4}$\;\;$ Let $v$ be an $\a$-flow on the compact metric space $(X,d)$. If for any chain $\tilde{\g}$,  
$$
|\int_{\tilde{\g}}\a|\le M
$$
for some $M>0$, then $v$ is a gradient-like flow.
}\\
  \\
{\bf Example 3.4.10}$\;\;$ Theorem 3.4.2 and 3.4.3 give a sufficient and necessary condition (3.4.5) for an $\a$-flow to be a non gradient-like flow on $X$. The following simple example shows that an $\a$-flow with $\a$ being a nontrivial cocycle can be a gradient-like flow if it does not satisfy the condition (3.4.5).\\
\indent
Let $S^1$ be the unit circle with the standard metric. Choose $r>0$ small enough and consider the following four open sets
\begin{align*}
U_1&=\{\ta\;|\;|\ta-\frac{\pi}{2}|<2r\}\\
U_2&=\{\ta\;|\;|\ta-\frac{3\pi}{2}|<2r\}\\
U_3&=\{\ta\;|\;\frac{\pi}{2}+r<\ta<\frac{3\pi}{2}-r\}\\
U_4&=\{\ta\;|\;-\frac{\pi}{2}+r<\ta<\frac{\pi}{2}-r\}
\end{align*}
Then $\cup^4_{i=1}U_i$ is a covering of $S^1$. Define bounded functions $\b_i(i=1,2,3,4)$ on $U_i$ as follows 
$$
\b_i(x)=0,\;\hbox{if}\;i=1,2
$$
and
\begin{align*}
\b_3(e^{i\ta})=\Big\{
\begin{array}{ll}
0 & \;\hbox{if}\;\frac{\pi}{2}+r<\ta<\pi\\
1 & \;\hbox{if}\;\pi\le \ta<\frac{3\pi}{2}-r
\end{array}
\end{align*}
\begin{align*}
\b_4(e^{i\ta})=\Big\{
\begin{array}{ll}
2 & \;\hbox{if}\;-\frac{\pi}{2}+r<\ta\le 0\\
0  & \;\hbox{if}\;0<\ta<\frac{\pi}{2}-r
\end{array}
\end{align*}
Consequently, we can define a cocycle $\a(x,y):S^1\times S^1\to {\Bbb R}$ to be\begin{align*}
\a(e^{i\ta_1},e^{i\ta_2})=\Big\{
\begin{array}{ll}
\b_i(e^{i\ta_2})-\b_i(e^{i\ta_1})&\;\hbox{if}\; |\ta_1-\ta_2|<r,\;\hbox{for}\; i=1,2,3,4 \\
0 & \;\hbox{if}\;|\ta_1-\ta_2|\ge r
\end{array}
\end{align*}
This cocycle is a nontrivial cocycle. Since if we let $\g$ be the oriented curve starting from the point 
$\ta=\frac{\pi}{2}$ and then going around the circle in the clockwise direction, then for any integer $l$,
 we have
$$
\int_{l\g}\a=l.
$$
\indent
Now we consider the gradient-like flow $v$ on $S^1$, that has two fixed points at $\ta=\pm\frac{\pi}{2}$, and flows from the point $\ta=\frac{\pi}{2}$ to the point $\ta=-\frac{\pi}{2}$. Then it is easy to check that $v$ carries the cocycle $\a$ with the parameter $\la=0,\r=1$ and $\G_2(r)\cup\G_3(r)=\emptyset$ in the definition 3.4.1.\\
 \\
{\bf Example 3.4.11}$\;\;$ Let $v$ be a flow on $S^1$ which has three fixed points at $\ta=0,\frac{2\pi}{3},\frac{4\pi}{3}$, and the forward direction of the flow is the anticlockwise direction. It is easy to see that $v$ carries a nontrivial cocycle and is a non gradient-like flow.\\
  \\
{\bf Example 3.4.12}$\;\;$ Let $v$ be the flow in Example 3.4.2, i.e., generated by a closed Morse 1-form $\o$. If $[\o]$ is a nontrivial cohomology class, then $v$ is a non-gradient-like flow.\\
\indent
This result can be seen as follows. We can lift $\o$ to a covering space $\bar{X}$ such that the pull-back form by the projection map is an exact form, say, $d\bar{f}$. Then the lifting flow $\bar{v}$ is the gradient flow of $\bar{f}$ on any cobordism $\bar{f}^{-1}[a,b]$ for $-\infty<a<b<\infty$. Since $\bar{f}$ is strictly decreasing along the nonconstant trajectory of $-\bar{v}$, there does not exist a global attractor in the interior of $\bar{f}^{-1}[a,b]$. Therefore, there is always a trajectory or a chain consisting of orbits joining fixed points such that the start point is at the level set $\bar{f}^{-1}(b)$ but the end point is at $\bar{f}^{-1}(a)$. Since $|a|, |b|$ can be arbitrarily large, hence by Theorem 3.4.2, $v$ on the underlying manifold $X$ is not a gradient-like flow.\\
 \\
{\bf Remark}$\;\;$ However, for general $\a$-flows, things would be complicated, since in the definition of the $\a$-flow, we use a ``global'' integration. Hence this allows some ``fluctuation'' to occur in the flow\\
  \\
{\bf 3.5 $\pi$-Morse decomposition of $\a$-flows}\\
  \\
\indent
From proposition 3.1.1, we know that there is a group embedding $\;J:\;\\C(X)/R(X)\to \mathcal{H}^1(X; {\Bbb R})$ whose image spans $\mathcal{H}^1(X; {\Bbb R})$. Let $\{f_i\}^m_{i=1}$ be functions in $C(X)$ such that the image $\a_i=J(f_i)(i=1,2,\cdots,m)$ is the cocycle representing the rank 1, indivisible integral cohomology class $[\a_i]$, while making $\{[\a_i]\}^m_{i=1}$ linearly independent. Let $[l_i]$ be the integral homology class dual to $[\a_i]$ for $i=1,\cdots,m$, i.e.,\\
\begin{align*}
\int_{[l_j]}[\a_i]=\int_{l_j}\a_i=\Big\{
\begin{array}{ll}
1 & i=j\\
0 & i\neq j 
\end{array}
\end{align*}
Let $\pi$ be a subgroup of the fundamental group of $X$ spanned by $\{l_1,\cdots,l_s\}$. Since $X$ is a compact polyhedron, there is a covering space $\bar{X}$ with deck transformation group $\pi$. Let $P_\pi:\;\bar{X}\to X$ be the projection map. Take a function $f_i:\;X\to S^1$ for $1\le i\le s$. There is a function $\hat{f}_i:=f_i\cdot P_\pi:\;\bar{X}\to S^1$. Since the pair
$$
\aligned
&\langle(f_i)_*\cdot(P_\pi)_*\cdot \pi_1(\bar{X},\bar{x}_0), d\ta\rangle=
\langle(P_\pi)_*\pi_1(\bar{X},\bar{x}_0), (f_i)^*(d\ta)\rangle\\
&=2\pi\langle(P_\pi)_*\pi_1(\bar{X},\bar{x}_0), \a_i\rangle=0,
\endaligned
$$
for $1\le i\le s$, each $\hat{f}_i$ has a lifting map $\bar{f}_i:\;\bar{X}\to {\Bbb R}$ such that $\hat{f}_i(\bar{x})=\exp{2\pi i\bar{f}_i(\bar{x})}$.\\
\indent
Assume that 1 is a regular value of $f_i$ for $1\le i\le s$. Define $\bar{N}_i(k)=\bar{f}_i^{-1}(k),\;
\bar{W}_i(k)=\bar{f}_i^{-1}([k,k+1])$ and $N_i=f_i^{-1}(1)$, then $P_\pi^{-1}(N_i)=\cup^{\infty}_
{k=-\infty}\bar{N}_i(k)$ and $l_i\cdot\bar{N}_i(0)=\bar{N}_i(1)$. Define $\bar{X}_0=
\{\bar{x}\in\bar{X};0\le\bar{f}_i(\bar{x})\le 1,1\le i\le s\}$, then $\bar{X}_0$ is a fundamental domain 
in $\bar{X}$. Let $\bar{N}_{i0}(1):=\bar{X}_0\cap\bar{N}_i(1)$ and let $\bar{N}_{i0}(0):=
\bar{X}_0\cap\bar{N}_i(0)$. Then the boundary of $\bar{X_0}$ is $\cup^s_{i=1}(\bar{N}_{i0}(0)
\cup\bar{N}_{i0}(1))$. The functions $\{\bar{f}_i;1\le i\le s\}$ divide the boundary of $\bar{X_0}$ into
 two parts:
$$
\bar{N}^+_0=\cup^s_{i=1}\bar{N}_{i0}(1);\\
\bar{N}^-_0=\cup^s_{i=1}\bar{N}_{i0}(0)
$$
  \\
{\sl{\bf Definition 3.5.1 ($\pi$-Morse decomposition)}$\;\;$ Let $v$ be a flow in a compact polyhedron with metric $d$. Let $\bar{v}$ be the lifting flow of $v$ in the covering space $\bar{X}$ which has a  deck transformation group $\pi$ as above. If $v$ satisfies the following conditions:\\
\indent
(1) The nonwandering set $\bar{A}$ of $\bar{v}$ has finitely many connected components $\{\bar{A}_n,\cdots,\bar{A}_1\}$ in the fundamental domain $\bar{X}_0$.
\\
\indent
(2) $\{\bar{N}^+_0,\bar{A}_n,\cdots,\bar{A}_1\}$ is a relative Morse decomposition of 
$\bar{v}$ in $\bar{X}_0$, where ``relative'' means that 
$\forall \bar{x}\in\bar{X}_0,\;\o(\bar{x})\cup\o^*(\bar{x})\not\subset\cup^n_{i=1}\bar{A}_i$, the trajectory $\bar{x}\cdot t$ goes through $\bar{N}^+_0(\bar{N}^-_0)$ at forward (backward) time and $\o(\bar{x})\cup\o^*(\bar{x})\notin\bar{X}_0$.\\
\indent
Then we say that $v$ has a $\pi$-Morse decomposition.}\\
 \\
{\sl {\bf Lemma 3.5.1}$\;\;$ Let $v$ be a flow in a compact space $X$ and let $P_\pi:\,\bar{X}\longrightarrow X$ be the covering projection mentioned above. If $\bar{A}$ is a nonwandering set in $\bar{X}$, then $P_\pi(A)$ is also a nonwandering set in $X$.}\\
  \\
Proof$\;\;$ Let $\bar{A}$ be a nonwandering set in $\bar{X}$ and let $\bar{x}\in\bar{A}$. Take an open neighborhood $\bar{U}$ of $\bar{x}$ such that $P_\pi:U\longrightarrow P_\pi(U)$ is a local homeomorphism. Since $\bar{x}\in \o(\bar{U})$, there is a sequence $\bar{x}_i\in\bar{U}$ and $t_i\in{\Bbb R}, t_i\to\infty$, such that for any smaller neighborhood $\bar{V}$ of $\bar{x}$, when $i$ becomes large enough, then $\bar{x}_i\cdot t_i\in V$. Fix $\bar{x}_i$ and project the integral curve $\bar{x}_i\cdot t$ to $X$, then we get an integral curve $P_\pi(\bar{x}_i\cdot t)=P_\pi(\bar{x}_i)\cdot t$ starting at $P_\pi(\bar{x}_i)\in P_\pi(\bar{U})$ ending at $P_\pi(\bar{x}_i\cdot t)\in P_\pi(V)$. Therefore $P(\bar{x})$ is a nonwandering point in $X$, which shows that $P_\pi(\bar{A})$ is a nonwandering set in $X$.\qed\\
\indent
By lemma 3.5.1, the set $P_\pi(\bar{A}_i)(i=1,\cdots,n)$ is a nonwandering set of $v$ in $X$, and is 
called a {\bf {$\pi$-stable nonwandering set}} of $v$. Denote the set 
$\{P_\pi(\bar{A}_n),\cdots,\\P_\pi(\bar{A}_1)\}$ by $\mathcal{N}_\pi(v)$.\\
\indent
Now we let $v$ be an $\a$-flow on $X$ and the cocycle $\a$ is nontrivial with rank $s$. Without loss of generality we can assume that the cohomology class $[\a]$ has the following representation:\\
\begin{equation}
[\a]=\sum^s_{i=1}\la_i[\a_i],\;\;\la>0\tag{$3.5.1$}
\end{equation}
for $i=1,\cdots,s$. (3.5.1) implies that there is a bounded function $\b:\,X\to{\Bbb R}$ such that 
\begin{equation}
\a-\sum^s_{i=1}\la_i\a_i\simeq \d\b\tag{$3.5.2$}
\end{equation}
We still use $\pi$ to represent the subgroup of the fundamental group of $X$ spanned by $\{l_1,\cdots,l_s\}$.\\ 
 \\
{\sl {\bf Theorem 3.5.2}$\;\;$ Let $v$ be an $\a$-flow with $\a$ being a nontrivial cocycle having rank s on a compact polyhedron $X$ with metric $d$. Then $v$ has a $\pi$-Morse decomposition with rank $\pi=s$ and the $\pi$-stable nonwandering set is $\mathcal{N}_\pi(v)=\{p_i\}^n_{i=1}$.
}\\
 \\
Proof$\;\;$ As before, we denote the lifting flow of $v$ in $\bar{X}$ by $\bar{v}$ and denote the lifting trajectory of any $\g$ by $\bar{\g}$. To show that $v$ has a $\pi$-Morse decomposition, we should prove the following facts:\\
\indent
(1) There is no oriented cycle in $\bar{X}_0$.\\
\indent
(2) For any $\bar{x}\in\bar{X}_0$, either $\bar{x}\cdot t$ goes through $\bar{N}^+_0(\bar{N}^-_0)$ and then out of $\bar{X}_0$ forever at forward(backward) time, or $\o(x)\cup \o^*(x)$ is contained in the union of the fixed points in $\bar{X}_0$.\\
\indent
Since an oriented cycle in $\bar{X}_0$ when projected down to $X$ will be an oriented cycle in $X$, therefore for the proof of (1), we will show that any oriented cycle in $X$ can be ``untied'' when lifting to the covering space $\bar{X}$, i.e., the lifting trajectories go from the infinitely far place to the infinitely far place.\\
\indent
Let $\g_c$ be an oriented cycle consisting of fixed points $p_{i_0},\cdots,p_{i_k}=p_{i_0}$ and trajectories $[p_{i_0},p_{i_1}],\cdots,[p_{i_{k-1}}, p_{i_k}]$. Consider the integration of $\a$ along $\g_c$. Applying the definition that $v$ is an $\a$-flow, we have 
$$
\int_{\g_c}\a\ge \sum^{k-1}_{j=0}\int_{[p_{i_j}, p_{i_{j+1}}]_r}\a\\
-\sum^{k-1}_{j=0}\osc_{x\in B_r(p_{i_j})}{\b_{p_{i_j}}(x)}\ge k\r-\la k\r>0
$$
where $[p_{i_j}, p_{i_{j+1}}]_r\in \G_1(r)$ is a sub-trajectory of $[p_{i_j}, p_{i_{j+1}}]$.\\
Notice that
$$
\sum^s_{i=1}\la_i\int_{\g_c}\a_i=\int_{\g_c}\a>0
$$
Hence there is an index $1\le i\le s$ such that
$$
\int_{\g_c}\a_i\ge 1
$$
Here we use the fact that $\a_i$ is an integral cocycle. Let $\bar{\g}_c$ be the lifting curve of $\g_c$ on $\bar{X}$. It is easy to see that the following holds:
$$
\int_{\bar{\g}\cap\bar{W}_i(0)}\d \bar{f}_i=\frac{1}{2\pi}\tr_{\g_c}f_i=\int_{\g_c}\a_i\ge 1
$$
and furthermore, for any integer $m>0$,
$$
\int_{\bar{\g}_c\cap(\cup^{m-1}_{k=0}\bar{W}_i(k))}\d \bar{f}_i\ge m
$$
This shows $\bar{\g}_c$ goes from the infinitely far place to infinitely far place. Therefore (1) is proved.\\
\indent
To prove (2), we only need to prove for any $x\in X$ with its limit set $\o(x)(\o^*(x))$ not in the union of all the fixed points of $v$ , that the lifting flow $\bar{x}\cdot t$ will pass through $\bar{N}^+_0(\bar{N}^-_0)$ then go to infinity at forward (backward) time.\\
\indent
Assume that $\o(x)\notin\cup^n_{i=1}\{p_i\}$. There are two cases which may occur for the trajectory $x\cdot t(t>0)$.\\
\indent
Case 1. $x\cdot t$ meets infinitely many balls $\{B_r(p_{i_0}),B_r(p_{i_1}),\cdots, B_r(p_{i_l}),\cdots\}$. Then the integration of $\a$ is 
$$
\int_{[x,x\cdot t]}\a\ge \sum^{l(t)-1}_{k=0}\int_{[x,x\cdot t]_k}\a-\sum^{l(t)-1}_{k=0}\osc_{x\in B_r(p_{l_k})}{\b_{p_{l_k}}(x)}
$$
where $l(t)$ is an integer representing the number of balls that $[x,x\cdot t]$ has met and $[x,x\cdot t]_k\in \G_1(r)$ is a trajectory of $[x,x\cdot t]$ between $B_r(p_{i_k})$ and $B_r(p_{i_{k+1}})$.\\
\indent
We have the estimate
\begin{equation}
\int_{[x,x\cdot t]}\a\ge l(t)(1-\la)\r\tag{$3.5.3$}
\end{equation}
Since $\a\simeq\sum^s_{i=1}\la_i\a_i+\d\b$ and $\b:\;X\to{\Bbb R}$ is a bounded function(let its bound be $M_\b$), then
\begin{gather}
\sum^s_{i=1}\la_i\int_{[x,x\cdot t]}\a_i\ge l(t)(1-\la)\r-|\b(x\cdot t)-\b(x)| \nonumber \\
\ge l(t)(1-\la)\r-2M_\b\tag{$3.5.4$}
\end{gather}
If $t$ is large enough such that $l(t)>\frac{2M_\b+C}{(1-\la)\r}$, then we have $$
\sum^s_{i=1}\la_i\int_{[x,x\cdot t]}\a_i>C
$$
for arbitrary large $C$. Now we can use the same argument as in (1) to show that the trajectory $[x,x\cdot t]$ must pass through $\bar{N}_0^+$ and then travel to the infinitely far place of $\bar{X}$.\\
\indent
Case 2. After meeting finitely many balls, $x\cdot t$ travels in the region $X\backslash \cup^n_{i=1} B_r(p_i)$.\\
\indent
In this case, we have (assume that $[x,x\cdot t]$ meets no balls).
$$
\int_{[x,x\cdot t]}\a\ge [\frac{t}{T_0}]\r
$$
In the same way as in case 1, we can prove that $[x,x\cdot t]$ will pass through $\bar{N}^+_0$ and then go to infinity.\\
\indent
We can also discuss the case that $\o^*(x)\notin\cup^n_{j=1}\{p_i\}$. In that case the trajectory 
$x\cdot t$ will pass through $\bar{N}^-_0$ and then go to infinity as $t\to -\infty$.\\
\indent
Now we have proved (2), and the theorem is proved.\qed\\
\indent
Using the same argument as in the proof of Theorem 3.5.2, we get \\
  \\
{\sl {\bf Theorem 3.5.3}$\;\;$ Let $v$ be a generalized $\a$-flow with respect to the nonwandering set 
$A=\{A_n,\cdots,A_1\}$ on a compact polyhedron $X$ with metric $d$. If $\a$ is a nontrivial cocycle 
having rank $s$, then $v$ has a $\pi$-Morse decomposition with respect to $A$ with rank $\pi=s$ and the $\pi$-stable nonwandering set is $\mathcal{N}_\pi(v)=A$.
}\\
 \\
{\bf\large 4 Novikov-Morse Theory}\\
  \\
\noindent
{\bf 4.1 Deformation of complexes relative to $\a$-flows}\\
  \\
\indent
In this section, we let $X$ be an $m$-dimensional compact polyhedron. Let $v$ be an $\a$-flow on $X$, where $\a$ is a nontrivial cocycle.\\
\indent
By Theorem 3.5.2, the flow $v$ has a $\pi$-Morse decomposition, where $\pi$ is a subgroup of the fundamental group of $X$ determined by the cocycle $\a$. Assume that $\pi$ is spanned by $G_+:=\{l_1,\cdots,l_s\}$, where $s$ is the rank of $\a$. Let $\bar{X}$ be the covering space of $X$ with deck transformation group $\pi$. As in section 3.5, the boundary of the fundamental domain $\bar{X}_0,\,\pat \bar{X}_0$ consists of two parts, $\bar{N}^+_0=\cup^s_{i=1}\bar{N}_{i0}(1)$, and $\bar{N}^-_0=\cup^s_{i=1}\bar{N}_{i0}(0)$ and satisfying for $i=1,2,\cdots,s$,
$$
l_i\cdot\bar{N}_{i0}(0)=\bar{N}_{i0}(1)
$$
\indent
Define $\pi_+$ to be the monoid constructed by $(G_+,e)$ with the group action from $\pi$.\\ 
\indent
Now the triangulation of $X$ provides a triangulation of $\bar{X}_0$ such that $\pat\bar{X}_0$ is an 
$m-1$-dimensional subcomplex. Define $R$ to be the subcomplex formed by those $m-1$-cells in 
$\bar{N}^-_0$. For any $p$-dimensional cell $e$ in $R$, we can define a map $\k$ from $R$ to $G_+$ such 
that $\k(e)$ is a replacement map sending $e$ to $\k(e)\cdot e\in \bar{N}^+_0$. \\
\indent
Now we fix the triangulation of $\pat\bar{X}_0$ and modify the inner simplices in $\bar{X}_0$. In detail, we will replace those cells supported on the boundary by cylinderical cells.\\
\indent
Note that $\pat\bar{X}_0$ is a simplicial complex, we can retract $\pat\bar{X}_0$ a small distance $\ape$ into the interior (if we choose a metric). $\ape$ is so small such that each cell in $\pat\bar{X}_0$ does not degenerate during the movement. In this way, we actually get a CW-decomposition of $\bar{X}_0$ with the cells around the boundary being cylindrical cells. Denote the $m-1$-dimensional subcomplex obtained by retracting $\bar{N}^+_0$ by distance $\ape$ by $R_I$.\\
\indent
Lifting  $\bar{X}_0$ to the universal covering $\tilde{X}$, we can get a fundamental domain $\tilde{X}_0$ in $\tilde{X}$. At this point, we forget all the information about the flow.\\
\indent
Now the lifted part of the boundary $\pat\bar{X}_0$ forms some part of the boundary $\pat\tilde{X}_0$. Hence the elements of $G_+$ are still replacement maps. We can view $\k(e)$ for $e\in R$ as a replacement map acting on $\tilde{R}$, the lifting of $R$.\\
\indent
Define a set $\bar{N}_0$ in $\bar{X}_0$ to be the space getting rid of all the cylindrical cells supported on $R_I$ having the form: $e\times\ape I, I=(0,1)$. Lifting $\bar{N}_0$ to $\tilde{X}_0$, we get $\tilde{N}_0$. Let $i_-:R\to\tilde{R}\cap\pat\tilde{X}_0$ and $i_+:R\to\tilde{R}_I\cap\pat\tilde{N}_0$ be two injections. Then $i_\pm$ are injections from $R$ to $\tilde{N}_0$. Define $f_q$ to be a map from $C_q(R)$ to $C_q(\tilde{N}_0)$:
$$
f_q=\k\cdot(i_-)_*(\cdot)-(i_+)_*(\cdot)
$$
where $\k\cdot(i_-)_*(\cdot)$ is a ${\Bbb Z}$-linear extension of the map defined on the basis as follows:\\
$\forall e\in C_q(R)$,
$$
\k\cdot(i_-)_*(e)=\k(e)\cdot(i_-)_*(e)
$$
Here $\k(e)$ is the element in $G_+$.\\
\indent
$\k\cdot(i_-)_*(\cdot)$ also commutes with the differential $d$. Take a simplex $e\in C_q(R)$ and let $d_R e=\sum_{e'}e'$, then we have 
$$
\aligned
&\k\cdot(i_-)_*(d_R e)=\sum_{e^\prime}\k\cdot(i_-)_*(e^\prime)=\sum_{e^\prime}\k(e^\prime)\cdot(i_-)_*(e^\prime)\\
&=\k(e)\cdot(\sum_{e^\prime}(i_-)_*(e^\prime))=\k(e)\cdot d_{\tilde{N}_0}(i_-)_*(e)\\
&=d_{\tilde{N}_0}(\k\cdot(i_-)_*(e))
\endaligned
$$
Therefore $f=\{f_q=\k\cdot(i_-)_*-i_+\}$ is a chain map from $C_*(R)$ to $C_*(\tilde{N}_0).$\\
\indent
Consequently, we can construct the algebraic mapping cone $(\hbox{con}(f)_*,d_c)$, whose $q$-th term and differential are
$$
\hbox{con}(f)_q=C_{q-1}(R)\oplus C_q(\tilde{N}_0)
$$
$$
d_c=\left(\begin{array}{cc}
-d_R & 0\\
f     & d_{\tilde{N}_0}
\end{array}
\right)
$$
\indent
Now our aim is to construct the chain equivalence between $(\text{con}(f)_*,d_c)$ and $(C_*(\tilde{X}_0),d_{\tilde{X}_0}))$.\\
\indent
Since the basis of the ${\Bbb Z}$-module chain complex $C_*(\tilde{X}_0)$ contains two kinds of cells, the cylindrical cells supported on $i_+(R)$ and the other cells in $\tilde{N}_0$. Each element $\b\in C_q(\tilde{X}_0)$ can be represented linearly and uniquely in the following form:
$$
\b=(-1)^{q-1}(i_+)_*(e)\times\ape I+s
$$
where $e\in C_{q-1}(R), I=(0,1)$ and $s\in C_q(\tilde{N}_0)$.\\
\indent
Define $\ta_q:\hbox{con}(f)_q\longrightarrow C_q(\tilde{X}_0)$ as 
$$ 
\ta_q(e\oplus s)=(-1)^{q-1}(i_+)_*(e)\times\ape I+s=\b
$$
and define $\e_q:C_q(\tilde{X}_0)\longrightarrow\hbox{con}(f)_q$ as 
$$
\e_q(\b)=e\oplus s
$$
Clearly, $\ta=\{\ta_q\}$ and $\e=\{\e_q\}$ are chain maps satisfying $\ta\cdot\e=\hbox{Id}_{\tilde{X}_0}$ and $\e\cdot\ta=\hbox{Id}_c$\\
  \\
{\sl {\bf Proposition 4.1.1}$\;\;$ Let $(\hbox{con}(f)_*,d_c),\ta$ and $\e$ be defined as above, then the following diagram commutes and $\ta$ defines a chain equivalence between $(\hbox{con}(f)_*,d_c)$ and $(C_*(\tilde{X}_0),d_{\tilde{X}_0})$.\\}
$$
\diagram
\hbox{con}(f)_q \dto^{\ta_q} \rto^{d_c} &\hbox{con}(f)_{q-1}\dto_{\ta_{q-1}}\\
C_q(\tilde{X}_0) \rto^{d_{\tilde{X}_0}} & C_{q-1}(\tilde{X}_0) 
\enddiagram
$$
  \\
Proof$\;\;$ Let $e\oplus s\in C_{q-1}(R)\oplus C_q(\tilde{N}_0)$, then 
$$
\aligned
&\ta_{q-1}\cdot d_c(e\oplus s)\\
&=\ta_{q-1}((-d_R e)\oplus(\k(e)(i_-)_*(e)-(i_+)_*(e)+d_{\tilde{N}_0}s))\\
&=-(-1)^{q-2}(i_+)_*(d_R e)\times\ape I+\k(e)(i_-)_*(e)-(i_+)_*(e)+d_{\tilde{N}_0}s
\endaligned
$$
On the other hand,
$$
\aligned
&d_{\tilde{X}_0}\cdot\ta_q(e\oplus s)\\
&=d_{\tilde{X}_0}((-1)^{q-1}(i_+)_*(e)\times\ape I+s)\\
&=(-1)^{q-1}(i_+)*(d_R e)\times\ape I+(i_+)_*(e)\times\ape-(i_+)_*(e)\times 0+d_{\tilde{N}_0}s\\
&=(-1)^{q-1}(i_+)_*(d_R e)\times\ape I+\k(e)\cdot(i_-)_*(e)-(i_+)_*(e)+d_{\tilde{N}_0}s
\endaligned
$$
Hence the diagram is commutative, and the claim is true.\qed\\
\indent
Let $\r_e:\pi_1(X)\to \pi$ be the extension of $\pi$ by the normal group $\pi_1(\bar{X})$. Since $\pi_+$ is a monoid in $\pi$, the set $\r^{-1}_e(\pi_+)$ is also a monoid of $\pi_1(X)$. We denote $\r^{-1}_e(\pi_+)$ by $\pi_1(X)_+$. ${\Bbb Z}\pi_1(X)_+$ is a subring of ${\Bbb Z}\pi_1(X)$.\\
\indent
Now we can tensor $C_*(R),C_*(\tilde{X}_0)$ and $\hbox{con}(f)_*$ with the ring ${\Bbb Z}\pi_1(X)_+$, then we have ${\Bbb Z}\pi_1(X)_+$-module chain complexes $({\Bbb Z}\pi_1(X)_+\otimes C_*(R),I\otimes d_R), ({\Bbb Z}\pi_1(X)_+\otimes C_*(\tilde{X}_0),I\otimes d_{\tilde{X}_0})$, and ${\Bbb Z}\pi_1(X)_+\otimes\hbox{con}(f)_*,I\otimes d_c)$, where
$$
I\otimes d_c=\left(\begin{array}{cc}
-I\otimes d_R & 0\\
I\otimes f & I\otimes d_{\tilde{N}_0}
\end{array}
\right)
$$
and $I\otimes f=I\otimes(\k\cdot(i_-)_*-(i_+)_*)$ is a ${\Bbb Z}\pi_1(X)_+$-chain map from 
${\Bbb Z}\pi_1(X)_+\otimes C_*(R)$ to ${\Bbb Z}\pi_1(X)_+\otimes C_*(\tilde{N}_0)$.\\
\indent
Let $\tilde{R}_+=\cup_{g\in\pi_1(X)_+}g\cdot i_-(R)$, then $C_*(\tilde{R}_+)$ becomes a
 ${\Bbb Z}\pi_1(X)_+$-module chain complex. Let $\tilde{X}_+=\cup_{g\in\pi_1(X)_+}g\cdot\tilde{X}_0$, 
 then $C_*(\tilde{X}_+)$ is also a ${\Bbb Z}\pi_1(X)_+$-module chain complex. Define a chain equivalence 
 $s_*:{\Bbb Z}\pi_1(X)_+\otimes C_*(\tilde{X}_0)\longrightarrow C_*(\tilde{X}_+)$ as follows,\\
for $g\otimes e\in {\Bbb Z}\pi_1(X)_+\otimes C_*(\tilde{X}_0)$,
$$
s_*(g\otimes e)=g\cdot e
$$
Extending the maps $(i_\pm)_*, f$ equivariantly to $(\tilde{i}_\pm)_*,\tilde{f}$, acting on $C_*(\tilde{R}_+)$, e.g., let $g\in \pi_1(X)_+, e\in C_*(R)$, then $g\cdot e\in C_*(\tilde{R}_+)$ and $\tilde{f}(g\cdot e)=g\cdot\tilde{f}(e)=g\cdot f(e)=g\cdot[\k\cdot(i_-)_*(e)-(i_+)_*(e)]$.\\
\indent
The following diagram is commutative
\begin{align*}
\diagram
{\Bbb Z}\pi_1(X)_+\otimes C_*(R) \rto^{I\otimes f}\dto^{s_*} &{\Bbb Z}\pi_1(X)_+\otimes C_*(\tilde{N}_0)\dto_{s_*}\\
C_*(\tilde{R})\rto^{\tilde{f}} &C_*(\tilde{N}_+)\\
\enddiagram
\tag{$4.1.1$}
\end{align*}
where $\tilde{N}_+=\cup_{g\in \pi_1(X)_+}g\cdot\tilde{N}_0$.\\
\indent
Consequently, we have the following commutative diagram with exact rows and with the vertical chain maps $s_*$ being ${\Bbb Z}\pi_1(X)_+$-isomorphic.
$$
\diagram
0 \rto &{\Bbb Z}\pi_1(X)_+\otimes C_{*-1}(R) \dto^{s_*}\rto &{\Bbb Z}\pi_1(X)_+\otimes\hbox{con}(f)_* \dto^{ss_*}\rto & {\Bbb Z}\pi_1(X)_+\otimes C_*(\tilde{N}_0) \dto^{s_*}\rto & 0\\
0 \rto &C_{*-1}(\tilde{R}_+)\rto &\hbox{con}(\tilde{f})_* \rto & C_*(\tilde{N}_+) \rto & 0
\enddiagram   
$$
where $\hbox{con}(\tilde{f})_*=C_{*-1}(\tilde{R}_+)\oplus C_*(\tilde{N}_+)$ is a ${\Bbb Z}\pi_1(X)_+$-module chain complex with differential,
$$
\tilde{d}_c=\left(\begin{array}{cc}
-d_{\tilde{R}_+} & 0\\
\tilde{f} & d_{\tilde{N}_+}
\end{array}
\right)
$$
Here the fact that $ss_*:=s_*\oplus s_*$ is a chain map is justified by the commutative diagram (4.1.1).\\
\indent
From the Five-Lemma for complexes, $ss_*$ is a chain equivalence. Thus we get the following corollary from proposition 4.1.1.\\
 \\
{\sl {\bf Corollary 4.1.2}$\;\;$ The two ${\Bbb Z}\pi_1(X)_+$-module chain complexes $(C_*(\tilde{X}_+),d_{\tilde{X}_+})$ and $(\hbox{con}(\tilde{f})_*,\tilde{d}_c)$ are chain equivalent.}\\
  \\
\indent
Similarly, if we extend the maps $\k\cdot(i_-)_*,(i_+)_*$ and $f$ to ${\Bbb Z}\pi_1(X)$-equivariant maps from $C_*(\tilde{R})$ to $C_*(\tilde{N})$, we have\\
  \\
{\sl {\bf Corollary 4.1.3}$\;\;$ The two ${\Bbb Z}\pi_1(X)$-module chain complexes $(C_*(\tilde{X}),\tilde{d}_{\tilde{X}})$ and $(\hbox{con}(\tilde{f})_*,\tilde{d}_c)$ are chain equivalent. Here the $q$-th term and differential of the mapping cone $(\hbox{con}(\tilde{f})_*,\tilde{d}_X)$ is 
$$
\hbox{con}(\tilde{f})_q=C_{q-1}(\tilde{R})\oplus C_q(\tilde{N})
$$
$$
\tilde{d}_c=\left(\begin{array}{cc}
-d_{\tilde{R}} & 0\\
\tilde{f} & d_{\tilde{N}}
\end{array}
\right)
$$
where $\tilde{f}=\widetilde{\k\cdot(i_-)_*}-\widetilde{(i_+)_*}$.
}\\
 \\
{\bf 4.2 Monodromy representations and Novikov numbers}\\
   \\
\indent
In this section, we will construct the Novikov complex associated with a cocycle $\a$ and establish some properties of this Novikov complex.\\
 \\
{\bf Abelization of fundamental groups}\\
  \\
\indent
Recall the group $\pi$ is relative to the $\a$-flow. Take the generator set $G_+$ of $\pi$, 
$G_+=\{l_1,\cdots,l_s\}$. There is a natural homomorphism
$$
\r_A:\,\pi\longrightarrow\mathbb{Z}^s,\,s=\hbox{rank}\pi
$$
which sends $l_i$ to $(0,\cdots,0,e_i,0\cdots,0)$, where $e_i$ is the unit 1.\\
\indent
Combining $\r_A$ with the extension homomorphism, we get a group homomorphism:
$$
\r_{\pi_1}=\r_A\cdot\r_e:\,\pi_1(X)\longrightarrow\mathbb{Z}^s
$$
This group homomorphism can be extended to a ring homomorphism:
$$
\r_{\pi_1}:\,\mathbb{Z}[\pi(X)]\longrightarrow\mathbb{Z}[\mathbb{Z}^s]
$$
However, there is a natural ring homomorphism $\r_q$ from $\mathbb{Z}[\mathbb{Z}^s]$ to the Laurent polynomial ring $Q_s=\mathbb{Z}[t_i,t_i^{-1};\,i=1,2,\cdots,s]$. $\r_q$ is defined as follows. Let $g=\sum z_j E_j$, where $E_j=(a_1,\cdots,a_s)$ is an $s$-dimensional vector with integral entries, then
$$
\r_q(g)=\sum z_j t^{E_j}=\sum z_j t_1^{a_1}\cdots t_s^{a_s}
$$
Hence combined with $\r_{\pi_1}$, we get the representation
$$
\r_{\pi_1 q}=\r_q\cdot \r_{\pi_1}:\,\mathbb{Z}[\pi_1(X)]\longrightarrow Q_s
$$
In fact, $\r_{\pi_1 q}$ is fully determined by the group $\pi$ and its representation. Restricting $\r_{\pi_1 q}$ to the subring $\mathbb{Z}[\pi_1(X)_+]$, we can get a ring homomorphism
$$
\r_P=\r_{\pi_1 q}|_{\mathbb{Z}[\pi_1(X)_+]}:\,\mathbb{Z}[\pi_1(X)_+]\longrightarrow P_s=\mathbb{Z}[t_1,\cdots,t_s]
$$
  \\
{\bf Monodromy representations}\\
  \\
\indent
Let $\tilde{E}$ be a local system of free abelian groups on the compact polyhedron $X$, then $\tilde{E}$ is determined by its monodromy representation $\r_{\tilde{E}}$:
$$ 
\r_{\tilde{E}}:\,\pi_1(X,x_0)\longrightarrow \text{Aut}(\tilde{E}_0)=GL(k,\mathbb{Z})
$$
where $\tilde{E}_0$ is the fibre of the free abelian group at $x_0$ and $k=\text{rank}(\tilde{E}_0)$. Let $E=\tilde{E}\otimes\mathbb{C}$, then $E$ is a complex flat vector bundle with the holonomy $\r_E$
$$
\r_E:\,\pi_1(X,x_0)\longrightarrow GL(k;\mathbb{Z})\otimes\mathbb{C}
$$
\indent
Now the tensor product of the representations $\r_q\otimes\r_{\tilde{E}}$ gives a representation of a $\mathbb{Z}[\pi_1(X)]$-ring to the linear space $(Q_s)^k$, where $Q_s$ is the polynomial space with $s$ variables over $\mathbb{Z}$.\\
\indent
Since $\r_{\tilde{E}}$ is an anti-homomorphism, i.e., $\forall g,g'\in \mathbb{Z}[\pi_1(X)], \r_E(g\cdot g')=\r_E(g')\cdot\r_E(g)$, hence $\r_P\otimes\r_{\tilde{E}}$ gives a right $\mathbb{Z}[\pi_1(X)]_+$-module structure on $P^k_s$. With the $P_s$-module structure of itself, $P_s^k$ becomes a $(P_s,\mathbb{Z}[\pi_1(X)_+])$-bimodule.\\
\indent
Define $D_*=P_s^k\otimes_{\mathbb{Z}[\pi_1(X)]_+}C_*(\tilde{X}_+)$, then $D_*$ is a $P_s$-module chain 
complex.\\
 \\
{\bf Evaluation representations}\\
  \\
\indent
Take any complex $s$-vector $a=(a_1,\cdots, a_s)\in \mathbb{C}^s$. The complex number field $\mathbb{C}$ can be given a $P_s$-module structure, whose module structure is provided by the action: for a polynomial $P(t_1,\cdots,t_s)$, $P(t_1,\cdots,t_s)\cdot x=P(a_1,\cdots,a_s)\cdot x=P(a)\cdot x$ for $x\in \mathbb{C}$. We denote the $P_s$-module of $\mathbb{C}$ evaluated at $t=a$ by $\mathbb{C}_a$. Similarly for any $a\in (\mathbb{C}^*)^s,\,\mathbb{C}$ can be viewed as a $Q_s$-module. If $p$ is a prime number, then the field $\mathbb{Z}_p$ also has a $P_s$-module structure which is given by the evaluation at $t=0$. We consider the complexes $\mathbb{C}_a\otimes_{P_s}D_*$ and $\mathbb{Z}_p\otimes_{P_s}D_*$. The following theorem for $s=1$ is given in the paper [13]. Here we consider the case $s\ge 1$  and use a different argument, which is based on section 4.1.\\
 \\
{\sl{\bf Theorem 4.2.1}$\;\;$ Let $D_*=P^k_s\otimes_{\mathbb{Z}[\pi_1(X)_+]}C_*(\tilde{X}_+)$. We have \\
\indent
(1) For any nonzero complex vector $a\in(\mathbb{C}^*)^s$, the homology $H_*(\mathbb{C}_a\otimes_{P_s}D_*)$ is isomorphic to $H_*(X,\,a^\a\otimes E)$, which is viewed as the homology of the presheaf $a^\a\otimes E$ on $X$.\\
\indent
(2) Let $p$ be a prime number and let $\mathbb{Z}_p$ have the $P_s$-module structure which is provided by the evaluation at $t=(t_1,\cdots,t_s)=0$. Then the homology $H_*(\mathbb{Z}_p\otimes_{P_s}D_*)$ is isomorphic to $H_*(\bar{X}_0,|R_I|\times \ape\bar{I};\,\mathbb{Z}_p\otimes P_\pi^*\tilde{E})$, where $\tilde{E}$ is a local system on $X$ and $P_\pi^* E$ is the pull-back local system on $\bar{X}_0$ by the projection $P_\pi:\,\bar{X}\longrightarrow X$.\\
\indent
(3) $H_*(C_0\otimes_{P_s}D_*)$ is isomorphic to $H_*(\bar{X}_0,|R_I|\times\ape\bar{I};\,P_\pi^*E)$.
}\\
  \\
Proof $\;\;$ Since the proof of (3) is the same as that of (2), we only give the proofs of (1) and (2).\\
\noindent
Proof of (1)\\
\indent
Since the $\mathbb{Z}[\pi_1(X)_+]$-basis of $C_*(\tilde{X}_+)$ is finite, all the complexes related to $C_*(\tilde{X}_+)$ are finitely generated, and hence all the homology groups are finitely generated.\\
\indent
For $a\in (C^*)^s$, we have the isomorphism
$$
\aligned
\mathbb{C}_a\otimes_{P_s}D^*\cong& \mathbb{C}_a\otimes_{P_s}((P_s)^k\otimes_{\mathbb{Z}[\pi_1(X)_+]}C_*(\tilde{X}_+))\\
\cong&(\mathbb{C}_a\otimes_{Q_s}(Q_s)^k)\otimes_{\mathbb{Z}[\pi_1(X)]}(\mathbb{Z}[\pi_1(X)]\otimes_{\mathbb{Z}[\pi_1(X)_+]}C_*(\tilde{X}_+))\\
\cong&(\mathbb{C}_a\otimes_\mathbb{Z}\mathbb{Z}^k)\otimes_{\mathbb{Z}[\pi_1(X)]}C_*(\tilde{X})\\
\cong&\mathbb{C}^k\otimes_{\mathbb{Z}[\pi_1(X)]}C_*(\tilde{X})
\endaligned
$$
Here the representation of $\mathbb{Z}[\pi_1(X)]$ is given by
$$ 
g\longrightarrow \r_{P(a)}(g)\otimes\r_E(g)
$$
Hence the homology of $\mathbb{C}^k\otimes_{\mathbb{Z}[\pi_1(X)]}C_*(\tilde{X})$ is the same as the homology of the presheaf $a^\a\otimes E$ on $M$ that corresponds to the flat vector bundle produced by the above holonomy representation. (1) is proved.\\
\indent
Proof of (2):\\
\indent
$$
\mathbb{Z}_p\otimes_{P_s}D^*\cong\mathbb{Z}_p\otimes_{P_s}(P^k_s\otimes_{\mathbb{Z}[\pi_1(X)_+]}C_*(\tilde{X}_+))\cong \mathbb{Z}^k_p\otimes_{\mathbb{Z}[\pi_1(X)_+]}C_*(\tilde{X}_+)
$$
Here the representation of $\mathbb{Z}[\pi_1(X)_+]$ on $\mathbb{Z}^k_p$ is 
\begin{equation}
g\longrightarrow \r_{P(0)}(g)\otimes\r_E(g)\tag{$4.2.1$}
\end{equation}
Since $\r_{P(0)}(g)=\r_{P(0)}\cdot\r_A\r_e(g)$, except in the case that $g\in\mathbb{Z}[\pi_1(X)_+]$ satisfies $\r_e(g)=0\in\mathbb{Z}^k$, the evaluation representation will make the final representation vanish. Hence (4.2.1) becomes
$$
\aligned
g\longrightarrow \r_E(g),\;\;&\hbox{if}\;\r_e(g)=0\\
g\longrightarrow 0,\;\;&\hbox{if}\;\r_e(g)\neq 0
\endaligned
$$
By corollary 4.1.2, in order to prove (2), we need to prove that $\mathbb{Z}^k_p\otimes_{\mathbb{Z}[\pi_1(X)_+]}con(\tilde{f})_*$ is equivalent to $(\mathbb{Z}_p\otimes P_\pi^*\tilde{E})\otimes C_*(\bar{X}_0,|R_I|\times\ape\bar{I})$. Now we have
$$
\mathbb{Z}^k_p\otimes_{\mathbb{Z}[\pi_1(X)_+]}con(\tilde{f})_q=(\mathbb{Z}^k_p\otimes_{\mathbb{Z}[\pi_1(X)_+]}C_{q-1}(\tilde{R}_+))\oplus(\mathbb{Z}^k_p\otimes_{\mathbb{Z}[\pi_1(X)_+]}C_q(\tilde{N}_+))
$$
\begin{equation}
I\otimes_{\mathbb{Z}[\pi_1(X)_+]}d_c=\left(
\begin{array}{cc}
-I\otimes_{\mathbb{Z}[\pi_1(X)_+]}\otimes\tilde{d}_{\tilde{R}_+} & 0\notag\\
I\otimes_{\mathbb{Z}[\pi_1(X)_+]}\tilde{f} & I\otimes_{\mathbb{Z}[\pi_1(X)_+]}\tilde{d}_{\tilde{N}_+}\notag
\end{array}
\right)
\end{equation}
If we take $\tilde{E}\otimes C_*(R)$ as the complexes of $R$ with twisted coefficients $\tilde{E}$, 
then $\mathbb{Z}_p^k\otimes_{\mathbb{Z}[\pi_1(X)_+]}C_*(\tilde{R}_+)\cong \mathbb{Z}_p
\otimes(\tilde{E}\otimes C_*(R))$. Similarly $\mathbb{Z}_p^k\otimes_{\mathbb{Z}[\pi_1(X)_+]}C_*(\tilde{N}_+)\cong \mathbb{Z}_p\otimes(P_\pi^*\tilde{E}\otimes C_*(\bar{X}_0,|R_I|\times\ape\bar{I}))$. Therefore 
$$
\mathbb{Z}^k_p\otimes_{\mathbb{Z}[\pi_1(X)_+]}con(\tilde{f})_*\cong (\mathbb{Z}_p\otimes(\tilde{E}
\otimes C_{*-1}(R)))\oplus(\mathbb{Z}_p\otimes(P_\pi^*\tilde{E}\otimes C_*(\bar{X}_0,|R_I|
\times\ape\bar{I})))
$$
and its differential is the same as
\begin{equation}
d_c=\left(
\begin{array}{cc}
-d_R & 0\notag\\
-(i_+)_* & d_{\bar{X}_0}\notag
\end{array}
\right)
\end{equation}
The conclusion of (2) now comes from the following lemma.\\
  \\
\indent
{\sl {\bf Lemma 4.2.2} $\;\;$ Let $N$ be a compact manifold with a piece of boundary denoted by $N_+$. Let $R$ be a topological space having a continuous map $i_+:\; R\longrightarrow N_+$ which is a homotopy equivalence, then the homology of the mapping cone $\mathcal{C}=con((i_+)_*:\;C_*(R)\to C_*(N))$ is isomorphic with that of $C_*(N,N_+)$.
}\\
  \\
Proof $\;\;$ We have the short exact sequence of complexes
$$
\diagram
0 \rto &C_*(N_+) \rto^{j} & C_*(N) \rto &C_*(N,N_+) \rto & 0
\enddiagram
$$  
Then we have the long exact sequence of homology groups:
$$
\diagram
\cdots \rto^{j_*} &H_q(N) \rto & H_q(N,N_+) \rto & H_{q-1}(N_+)\rto^{j_*} & H_{q-1}(N) \rto &\cdots
\enddiagram
$$    
Similarly, from the short exact sequence of the mapping cone
$$
\diagram
0 \rto & C_*(N) \rto &con(i_+)_* \rto &C_{*-1}(R) \rto & 0  
\enddiagram
$$
we get the long exact sequence
$$
\diagram
\cdots \rto & H_q(R) \rto^{(i_+)_*} & H_q(N)\rto &H_q(con(i_+)) \rto &H_{q-1}(R)\rto & \\
 & &\qquad \rto&H_{q-1}(N) \rto &\cdots& 
\enddiagram
$$ 
The two long exact sequences yield the following two short exact sequences
$$
\diagram
0 \rto & \text{Coker}(j_*:\,H_q(N_+)\to H_q(N)) \rto &H_q(N,N_+)\rto &\qquad \\
 \qquad\rto & \text{Ker}(j_*:\,H_{q-1}(N_+)\to H_{q-1}(N)) \rto & 0& 
\enddiagram
$$ 
and
$$
\diagram
0 \rto &\text{Coker}((i_+)_*:\,H_q(R)\to H_q(N)) \rto &  H_q(con(i_+)_*)\rto & \qquad\\
\qquad\rto &\text{Ker}((i_+)_*:\,H_{q-1}(R)\to H_{q-1}(N)) \rto & 0& 
\enddiagram
$$  
Note that $(i_+)_*:\,H_q(R)\to H_q(N)$ is the combination of the isomorphisms $H_q(R)\to H_q(N_+)$ and $j_*:\,H_q(N_+)\to H_q(N)$, and therefore using the Five-Lemma, we get for any $q\ge 0$
$$
H_q(N,N_+)\cong H_q(con(i_+)_*)
$$
\qed
  \\
{\bf Novikov numbers}\\
 \\
\indent
In Theorem 4.2.1, we have considered the complex $D_*=P^k_s\otimes_{{\Bbb Z}[\pi_1(X)_+]}C_*(\tilde{X}_+)$. In this part, we always let the vector bundle $E$ that appears in the previous section be a trivial line bundle. Therefore the complex $D_*$ we consider here has the form $D_*=P_s\otimes_{{\Bbb Z}[\pi_1(X)_+]}C_*(\tilde{X}_+)$. Since the representation of ${\Bbb Z}[\pi_1(X)_+]$ in $P_s$ is completely determined by the cohomology class $[\a]$, we denote the homology group $H_*({\Bbb C}_a\otimes_{P_s}D_*)$ as $H_*(X,a^\a)$, or in other words, view the homology group $H_*({\Bbb C}_a\otimes_{P_s}D_*)$ as the homology group of the presheaf $a^\a$ on $X$ which is given by the monodromy representation $\r_P:{\Bbb Z}[\pi_1(X)_+]\to P_s={\Bbb Z}[t_1,\cdots,t_s]$.\\
 \\
{\bf Definition}$\;\; b_i([\a]):=\hbox{rank}(H_i(D_*))$ for $i=0,1,\cdots,m$ are called the Novikov numbers.\\
  \\
\indent
The following is essentially given in [20] and the proof here is essentially the same as in the paper [12] except that we consider here $s\ge 1$ and $X$ a compact polyhedron. For the sake of completeness, we give the proof here.\\
 \\
{\sl {\bf Theorem 4.2.3}$\;\;$ Let $X$ be a compact polyhedron. Define a function for fixed $i$ to be 
$$
a\in ({\Bbb C}^*)^s\longrightarrow \dim_{{\Bbb C}}H_i(X,a^\a),
$$ 
then it has the following properties:\\
\indent
(1) It is generically constant, more precisely, except on a proper algebraic subvariety $\mathcal{L}$ in $({\Bbb C}^*)^s$, the dimension $\dim_{{\Bbb C}}H_i(X,a^\a)$ is constant and this constant is just the Novikov number $b_i([\a])$ we defined above.\\
\indent
(2) For any point $\hat{a}\in\mathcal{L}$, 
$$
\dim_{{\Bbb C}}H_i(X,\hat{a}^\a)>b_i([\a]) 
$$
}\\
 \\
Proof$\;\;$ we denote by $B_i(a)$ the rank of the following linear map
$$
I\otimes d:{\Bbb C}_a\otimes_{P_s}D_i\longrightarrow {\Bbb C}_a\otimes_{P_s}D_{i-1}
$$
It is easy to see that this map $I\otimes d$ can be represented by a matrix with entries being polynomials with $s$ variables on $({\Bbb C}^*)^s$. Therefore, we know that except on a proper algebraic variety in $({\Bbb C}^*)^s$, $B_i(a)$ is constant and on that algebraic variety $B_i(a)$ is smaller than the generic constant.\\
\indent
Consider the following truncated complex,
$$
\cdots\longrightarrow {\Bbb C}_a\otimes D_i\longrightarrow {\Bbb C}_a\otimes D_{i-1}\longrightarrow 0,
$$
and using the Euler-Poincar\'e formula, we obtain 
\begin{align*}
\sum^\infty_{j=0}(-1)^j\hbox{rank}({\Bbb C}_a\otimes D_{i-1+j})=
\sum^\infty_{j=0}(-1)^j\dim_{{\Bbb C}}(H_j({\Bbb C}_a\otimes D_{i-1+*}))
\end{align*}
Therefore we obtain
\begin{align*}
\dim_{{\Bbb C}}H_i(X;a^\a)&=\sum^\infty_{j=0}(-1)^j\hbox{rank}({\Bbb C}_a\otimes D_{i+j})\\
-&\sum^\infty_{j=1}(-1)^j\hbox{rank}(H_{i+j}({\Bbb C}_a\otimes D_*))-B_i(a)
\end{align*}
Now it is easy to see that (1) and (2) are true.\\
 \\
{\bf Remark}$\;\;$ If $s=1$ and $X$ is a closed manifold, then it was proved in the paper [11] that the Novikov number we defined here is the same as the Novikov number which is defined as the dimenion of the homology group of the Novikov complex and which was originally introduced by S.S.Novikov in his papers [18],[19] and [20] (see the introduction).\\
   \\
{\bf 4.3 Morse type inequalities for a filtration}\\
  \\
\indent
In this part, we always assume $\mathcal{R}$ to be a commutative Noetherian ring, and $C$ be a f.g. finite $\mathcal{R}$-module chain complex. Suppose that $C$ has a filtration of finite length:$\{0\}\sbs C_{n+1}\sbs\cdots\sbs C_1\sbs C_0=C$. We will give a Morse type inequality for $C$ under the above filtration.\\
\indent
First, we consider a triple $(\tau',\tau'',\tau)$ of complexes, satisfying the inclusion relations $\tau'\sbs\tau''\sbs\tau$. Then we have a short exact sequence of complexes
\begin{equation}
0\longrightarrow \tau''/\!\tau'\longrightarrow\tau/\!\tau'\longrightarrow\tau/\!\tau''\longrightarrow 0\tag{$4.3.1$}
\end{equation}
Define $C(\mu,\tau)=\tau/\mu$ to be the quotient complex of complexes $\tau$ and $\mu$. Let $H_*(\mu,\tau)=H_*(C(\mu,\tau))$ be the homology groups of the quotient complex $C(\mu,\tau)$. We get from (4.3.1) the long exact sequence of the homology groups of the triple $(\tau',{\tau}'',\tau)$,\\
\begin{equation}
\cdots\longrightarrow H_1({\tau}'',\tau)\stackrel{\d_1}{\longrightarrow} H_0({\tau}',{\tau}'')\longrightarrow H_0({\tau}',\tau)\longrightarrow H_0({\tau}'',\tau)\stackrel{\d_0}{\longrightarrow} 0\tag{$4.3.2$}
\end{equation}
where $\d_j:\;H_j({\tau}'',\tau)\longrightarrow H_{j-1}({\tau}',{\tau}''), j=0,1,\cdots$ is the connecting homomorphism.\\
\indent
Introduce some notations as follows\\
$$
\aligned
b&_j(\mu,\tau)=\hbox{rank}(H_j(\mu,\tau))\\
d&_j({\tau}',{\tau}'',\tau)=\hbox{rank}(\hbox{im} \d_j)\\
p&(\mu,\tau,;t)=\sum_{j\ge 0}b_j(\mu,\tau)t^j\\
q&(\tau',\tau'',\tau;t)=\sum_{j\ge 0}d_j(\tau',\tau'',\tau)t^j
\endaligned 
$$
Then we have\\ 
  \\
{\sl {\bf Lemma 4.3.1}$\;\;$ If $({\tau}',{\tau}'',\tau)$ is the triple of f.g. $\mathcal{R}$-module chain complexes in (4.3.1), then
\begin{equation}
p({\tau}'',\tau;t)-p({\tau}',\tau;t)+p({\tau}',{\tau}'';t)=(1+t)q({\tau}',{\tau}'',\tau;t)\tag{$4.3.3$}
\end{equation}
}\\
  \\
Proof$\;\;$Truncating the long exact sequence of (4.3.2), we obtain the exact sequence for $j=0,1,\cdots$ (we define $\d_{-1}=0$),
$$
0\rightarrow\mbox{im}(\d_{j})\rightarrow H_j({\tau}',{\tau}'')\rightarrow H_j({\tau}',\tau)\rightarrow H_j({\tau}'',\tau)\rightarrow\mbox{im}\d_{j-1}\to 0
$$
Now it is easy to get the result by using the standard method.\qed\\
\indent
For the given filtration of $C$, we let $({\tau}',{\tau}'',\tau)=(C_{n+1},C_j,C_{j-1}),\;j=1,\cdots,n$. By lemma 4.3.1, we have for $j=1,\cdots,n$
$$\aligned
&p(C_j,C_{j-1};t)-p(C_{n+1},C_{j-1};t)+p(C_{n+1},C_j;t)\\
&=(1+t)q(C_{n+1},C_j,C_{j-1};t)
\endaligned
$$
Hence we get\\
  \\
{\sl{\bf Proposition 4.3.2}$\;\;$ If $C$ is a f.g.finite $\mathcal{R}$-module chain complex with filtration $ \{0\}\sbs C_{n+1}\sbs C_n\sbs\cdots \sbs C_0=C$, then we have
$$
\sum^{n+1}_{j=1}p(C_j,C_{j-1};t)=p(C_{n+1},C;t)+(1+t)\sum^n_{j=1}q(C_{n+1},C_j,C_{j-1};t)
$$
}
  \\
\noindent
{\bf Morse type inequalities for a filtration of a $G$-CW complex}\\
  \\
\indent
Let $G$ be a topological group acting on a CW-complex $X$. Suppose that $X$ is a finite $G$-CW complex(i.e., $X/\! G$ is a finite CW-complex) with a $G$-invariant filtration
\begin{equation}
\emptyset\sbs X_{n+1}\sbs X_n\sbs\cdots X_0=X\tag{$4.3.4$}
\end{equation}
\noindent
{\bf Remark}$\;\;$ For a $G$-CW complex $X$, there is a natural filtration filtered by its own skeleton, 
e.g.,
$$
\emptyset=X^{-1}\sbs X_0\sbs\cdots X^p\sbs \cdots \cup_{p\ge 0}X^p=X
$$
The $p$-th chain group is $C_p(X)=H_p(X^p,X^{p-1})$, and its basis consists of the $p$-dimensional 
equivariant cells. The differential $d_p:\;C_p(X)\to C_{p-1}(X)$ is the connecting homomorphism of the 
triple $(X^p,X^{p-1},X^{p-2})$.\\
  \\
\indent
Let $\mathcal{R}$ be a commutative Noetherian ring and $\r:\;{\Bbb Z}G\to\hbox{End}(\mathcal{R})$ be a representation of the group $G$ in $\mathcal{R}$. Then (4.3.4) induces a filtration of the $\mathcal{R}$-module chain complex $\mathcal{R}\otimes_{{\Bbb Z}G}C_*(X)$,
\begin{equation}
0\sbs \mathcal{R}\otimes_{{\Bbb Z}G}C_*(X_{n+1})\sbs\cdots\sbs\mathcal{R}\otimes_{{\Bbb Z}G}C_*(X_0)=\mathcal{R}\otimes_{{\Bbb Z}G}C_*(X)\tag{$4.3.5$}
\end{equation}
Denote $\mathcal{R}\otimes_{{\Bbb Z}G}C_*(X_j)$ by $D_j$ for $j=1,\cdots,n+1$, then (4.3.5) becomes
\begin{equation}
0\sbs D_{n+1}\sbs\cdots\sbs D_0=D\tag{$4.3.6$}
\end{equation}
Since $X$ is a finite $G$-CW complex, the $\mathcal{R}$-module chain complex $D$ is finitely generated. Hence we have\\
  \\
{\sl {\bf Proposition 4.3.3}$\;\;$ Let $\mathcal{R}$ be a commutative Noetherian ring with a $G$-representation. Let $X$ be a finite $G$-CW complex with a $G$-invariant filtration (4.3.4), then
\begin{equation}
\sum^{n+1}_{j=1}p(D_j,D_{j-1};t)=p(D_{n+1},D;t)+(1+t)\sum^n_{j=1}q(D_{n+1},D_j,D_{j-1};t)\tag{$4.3.7$}
\end{equation}
where $D_j=\mathcal{R}\otimes_{{\Bbb Z}G}C_*(X_j)$, and $p(\cdot,\cdot;t),q(\cdot,\cdot,\cdot;t)$ are defined as above, e.g.,
$$
\aligned
&p(D_j,D_{j-1};t)=\sum_{l\ge 0}b_l(D_j,D_{j-1})t^l\\
&b_l(D_j,D_{j-1})=\hbox{rank}(H_l(\mathcal{R}\otimes_{{\Bbb Z}G}C_*(X_j,X_{j-1})))
\endaligned
$$
}  
  \\
{\bf Comparison of Poincar\'e polynomials for prime ideals}\\
  \\
\indent
Let $\mathcal{R}$ be the ring mentioned above and $C$ a chain complex over the ring $\mathcal{R}$ . Take a prime ideal $\mathcal{P}$ in $\mathcal{R}$, and we can define a chain complex $C(\mathcal{P})$ over the quotient field $Q(\mathcal{R}/\mathcal{P})$ of the entire ring $\mathcal{R}/\mathcal{P}$,\\
$$
C(\mathcal{P}):=Q(\mathcal{R}/\mathcal{P})\otimes_\mathcal{R} C
$$
Define the $j$-th Betti number of $C(\mathcal{P})$:
$$
b_j(C;\mathcal{P})=b_j(\emptyset,C;\mathcal{P})=\dim_{Q(\mathcal{R}/\mathcal{P})}(H_j(Q(\mathcal{R}/\mathcal{P})\otimes_\mathcal{R} C))
$$
Let $\mathcal{Q}$ be another prime ideal such that $\mathcal{P}\sbs\mathcal{Q}$. We can define the torsion number $T_j(\mathcal{P},\mathcal{Q})$ of $C(\mathcal{P})$ with respect to $\mathcal{Q}$ as:
$$
T_j(\mathcal{P},\mathcal{Q})=\hbox{rank}_\mathcal{P}(\hbox{im}d_{j+1}(\mathcal{P}))-\hbox{rank}_{\mathcal{Q}}(\hbox{im}d_{j+1}(\mathcal{Q}))
$$
where $d_j(\mathcal{P})$ represents the boundary operator in $C(\mathcal{P})$.\\
\indent
In the same way, we define the Poincar\'e polynomial of $C(\mathcal{P})$ as:
$$
P(C;t,\mathcal{P}):=\sum_{j\ge 0}b_j(C,\mathcal{P})t^j
$$
\indent
The following result is proved in [13]. However, since we need more information than provided in [13], we will give a different proof below.\\
  \\
{\sl {\bf Proposition 4.3.4}$\;\;$ Let $\mathcal{R}$ be a commutative Noetherian ring and $C$ a chain complex on $\mathcal{R}$. If $\mathcal{P},\mathcal{Q}$ are two prime ideals in  $\mathcal{R}$ such that $\mathcal{P}\sbs\mathcal{Q}$, then there is a polynomial $Q(\mathcal{P},\mathcal{Q};t)$ with nonnegative integer coefficients such that 
\begin{equation}
P(C;t,\mathcal{Q})=P(C;t,\mathcal{P})+(1+t)Q(\mathcal{P},\mathcal{Q};t)\tag{$4.3.8$}
\end{equation}
where $(1+t)Q(\mathcal{P},\mathcal{Q};t)=\sum_{j\ge 0}(T_j(\mathcal{P},\mathcal{Q})+T_{j-1}(\mathcal{P},\mathcal{Q})t^j$.
}\\
 \\
Proof$\;\;$ For any $j=0,1\cdots,$we have two short exact sequences.
$$
0\longrightarrow\ker{d_j(\mathcal{P})}\longrightarrow C_j(\mathcal{P})\longrightarrow\text{im}d_j(\mathcal{P})\longrightarrow 0
$$
$$
0\longrightarrow\text{im}d_{j+1}(\mathcal{P})\longrightarrow\ker{d_j(\mathcal{P})}\longrightarrow H_j(C(\mathcal{P}))\longrightarrow 0
$$
Therefore
$$\aligned
\text{rank}_\mathcal{P}C_j(\mathcal{P})=&\text{rank}_\mathcal{P}(\ker{d_j(\mathcal{P})})+\text{rank}_\mathcal{P}(\text{im}d_j(\mathcal{P}))\\
=&\text{rank}_\mathcal{P}(H_j(C(\mathcal{P})))+\text{rank}_\mathcal{P}(\text{im}d_j(\mathcal{P}))+\text{rank}_\mathcal{P}(\text{im}d_{j+1}(\mathcal{P}))
\endaligned
$$
In the same way, for the prime ideal $\mathcal{Q}$, we have 
$$
\text{rank}_\mathcal{Q}C_j(\mathcal{Q})=\text{rank}_\mathcal{Q}(H_j(C(\mathcal{Q})))+\text{rank}_q(\text{im}d_j(\mathcal{Q}))+\text{rank}_\mathcal{Q}(\text{im}d_{j+1}(\mathcal{Q})).
$$
Since $C_j$ is free module on $\mathcal{R}$, we have 
$$
\text{rank}_\mathcal{P}C_j(\mathcal{P})=\text{rank}_\mathcal{Q}C_j(\mathcal{Q})
$$
Thus, by taking the difference we get
$$
b_j(C;\mathcal{Q})-b_j(C;\mathcal{P})=T_j(\mathcal{P},\mathcal{Q})+T_{j-1}(\mathcal{P},\mathcal{Q})
$$
and then
$$
P(C;t,\mathcal{Q})-P(C;t,\mathcal{P})=\sum_{j\ge 0}(T_j(\mathcal{P},\mathcal{Q})+T_{j-1}(\mathcal{P},\mathcal{Q}))t^j
$$
Note that the Euler-Poincar\'e equality holds:
$$
P(C;-1,\mathcal{P})=P(C;-1,\mathcal{Q})
$$
which implies that the right hand side can be written as $(1+t)Q(\mathcal{P},\mathcal{Q};t)$ with $Q(\mathcal{P},\mathcal{Q};t)$ a polynomial with nonnegative integer coefficients.\\
  \\
{\bf 4.4 Isolated invariant sets and Conley index}\\
  \\
\indent
In this section, we will explain the concepts of isolated invariant sets and Conley index which were introduced by C.Conley (see [9],[10]) as a generalization of critical points and Morse index of a Morse function. Using the Conley index pair, we can generalize the Novikov inequalities for counting the zeros of a closed Morse 1-form to the Novikov-Morse type inequalities for the isolated $\a$-flows which will be defined in section 4.5.\\
\noindent
{\bf Isolated invariant sets}\\
\indent
Let a flow $v$ be defined on the space $X$ and let $N\sbs X$ be a compact set. $N$ is called an isolated invariant set of the flow if the invariant set $I(N)\sbs \text{int}N$.\\
\indent
It is easy to prove that if $D=\{M_n,\cdots,M_1\}$ is a Morse decomposition of an isolated invariant set $S$, then each Morse set $M_i(i=1,\cdots,n)$ is an isolated set.\\
\noindent
{\bf Index pair}\\
\indent
Let $S$ be an isolated invariant set. A compact pair $(N_0,N_1)$ in $X$ is called an index pair for $S$, if \\
\indent
(1) $\overline{N_1\backslash N_0}$ is an isolated neighborhood for $S$.\\
\indent
(2) $N_0$ is positively invariant relative to $N_1$.\\
\indent
(3) if $v\in N_1$ and $v\cdot {\Bbb R}^+\not\sbs N_1$, then there is a $t\ge 0$ such that $v\cdot[0,t]\sbs N_1$ and $v\cdot t\in N_0$\\
\indent
Here condition (2) means that if $v\in N_0$ and $v\cdot[0,t]\sbs N_1$, then $v\cdot [0,t]\sbs N_0$. (3) implies that for any point $v\in N_1$, either $v$ flows into the invariant set or flows out of $N_1$ at finite time through the ``exit set'' $N_0$.\\
\indent
The index pair has the following two important properties:\\
\noindent
{\bf Homotopy invariance}\\
\indent
Let $(N_0,N_1)$ and $(\bar{N}_0,\bar{N}_1)$ be two index pairs of an isolated invariant set $S$. Then the two index pairs as space pairs are homotopically equivalent.\\
\noindent
{\bf Continuation principle}\\
\indent
Since we do not use this property in this paper, we refer the interested reader to the fourth chapter of [9].
 \\
\indent
Because of the homotopy invariance, the index pair defined by C.Conley is called the Conley index, and we define $h(S)$ to be the homotopy type of the index pair $(\ast,N_1/\! N_0)$ of the isolated invariant set $S$.\\
 \\
{\bf Computation of the Conley index}\\
  \\
\noindent
{\bf Sum formula}\\ 
\indent
If $S_1$ and $S_2$ are isolated invariant sets with $S_1\cap S_2=\emptyset$, then $S_1\cup S_2$ is an isolated invariant set and $h(S_1\cup S_2)=h(S_1)\vee h(S_2)$, where $\vee$ denotes the join of pointed topological spaces.\\
\indent
This formula can be achieved as follows. Let $(N_0,N_1)$ and $(\bar{N}_0,\bar{N}_1)$ be the index pairs of $S_1$ and $S_2$ respectively, then $(N_0\cup\bar{N}_0, N_1\cup\bar{N}_1)$ is the index pair for the isolated invariant set $S_1\cup S_2$. Hence $h(S_1\cup S_2)=(\ast, N_1\cup \bar{N}_1\backslash \! N_0\cup\bar{N}_0)=h(S_1)\vee h(S_2)$\\
\noindent
{\bf Product formula}\\
\indent
Let $S_1$ and $S_2$ be isolated invariant sets with index pairs $(N_0,N_1)$ and $(\bar{N}_0,\bar{N}_1)$ respectively, then $S_1\times S_2$ is an isolated invariant set with index pair $(N_0,N_1)\times(\bar{N}_0,\bar{N}_1)=(\pat N_0\times\bar{N}_1\cup\pat\bar{N}_0\times N_1, N_1\times \bar{N}_1)$. Hence $h(S_1\times S_2)=(\ast, N_1\times\bar{N}_1\backslash\pat N_0\times\bar{N}_1\cup\pat\bar{N}_0\times N_1)=(\ast, N_1\backslash N_0)\wedge (\ast,\bar{N}_1\backslash\bar{N}_0)=h(S_1)\wedge h(S_2)$. Here ``$\wedge$'' is the smash product between two pointed spaces.\\
\noindent
{\bf Index of hyperbolic fixed points of the flow}\\
\indent
Let $x_0$ be a fixed point of the flow $\dot{v}=V(v)$. Consider its linearized equation at $x_0$,
\begin{equation}
\dot{v}=DV(x_0)\cdot v\tag{$4.4.1$}
\end{equation}
$x_0$ is called a hyperbolic fixed point if the real part of each eigenvalue of the coefficient metrix $DV(x_0)$ is not zero. Hence the tangent space $T_{x_0}M=E^+\oplus E^-$, where $E^+(E^-)$ represents the eigenspace of the eigenvalues having positive(negative) real part. Let $k=\text{rank}(E^+)$ and $D^k$ be the $ k$-dimensional unit disc in $E^+$. Since the local unstable manifold and the local stable manifold are transversal at the hyperbolic fixed point $x_0$, the Conley index $h(v_0)=(\ast,D^k/\!\pat D^k)=(\ast, S^k)$.\\
\noindent
{\bf Index of an orientable hyperbolic periodic orbit}\\
\indent
Let $x_0$ be any point in the hyperbolic periodic orbit $v_0$ and let $T^\bot_{x_0}M$ be the tangent subspace vertical to $\dot{v}_0$. Then the derivative of the Poincar\'e map defines the stable eigenspace $E^-$ and the unstable eigenspace $E^+$. Let $k=\dim{E^+}$, then this means that the eigenspace of the unstable manifold is orientable, the Conley index $h(v_0)=(\ast, S^1\amalg \{\ast\})\wedge(\pat D^k,D^k)$, where $(\ast,S^1\amalg \{\ast\})$ is the Conley index of the flow $\dot{\ta}=1$ in $S^1$. Therefore,
$$ 
\aligned
h(v_0)&\cong (\ast,S^1\amalg\{\ast\})\wedge (\ast, S^k)\\ 
&\cong (\ast,S^1\amalg\{\ast\})\wedge (\ast, S^1)\wedge (\ast,S^{k-1})\\
&\cong (\ast,S^2\vee S^1)\wedge (\ast, S^{k-1})\\
&\cong (\ast,S^k)\vee (\ast,S^{k+1})
\endaligned
$$
i.e., the index of $v_0$ is the sum of a pointed $k$-sphere and a pointed $k+1$-sphere.\\
\noindent
{\bf Index of an unorientable hyperbolic periodic orbit}\\
\indent
Let $P_{x_0}$ be the Poincar\'e map at a point of the hyperbolic orbit. Without loss of generality, we can assume that $dP_{x_0}\in O(m-1)$. Now the unstable manifold is the integral submanifold of $E^+$. If the unstable manifold is unorientable, $E^+\times_P S^1$ is a twisted vector bundle on $S^1$ twisted by $dP_{x_0}$. We can connect $dP_{x_0}$ with $-I\in O(m-1)$ by a path, so the homotopy type of $E^+\times_{-I}S^1$ is the same as $E^+\times_P S^1$. Hence we can assume that $dP_{x_0}$ reverses only $1$ eigensubspace and the Conley index of $v_0$ is the product of a pointed $k-1$-sphere with a M$\ddot{o}$bius band collapsing its boundary, hence 
$$
h(v_0)=(\ast,{\Bbb RP}^2\amalg \{\ast\})\wedge(\ast, S^{k-1})
$$
{\bf Index of a critical manifold}\\
\indent
Let $Z$ be a non-degenerate critical manifold of a smooth function $f$. This means that the following two conditions:\\
\indent
(1) $\hbox{grad}f(x)=0,\forall x\in Z$.\\
\indent
(2) at any point $x\in Z$, there is a decomposition $T_x M=T^\perp_x(Z)\oplus T_x(Z)$ and $d^2f(x)$ as bilinear form is nondegenerate along the vertical tangent subspace $T^\perp_x(Z)$\\
\indent
The subspace $T^\perp_x(Z)$ can be decomposed again into the direct sum of $E^+_x$ and $E^-_x$ which correspond to the eigenspaces of positive eigenvalues and negative eigenvalues of $d^2f|_{T^\perp_x(Z)}$.\\
\indent
If $Z$ is connected, $\text{rank}E^-_x$ does not depend on the point $x$. Hence $\text{ind}(Z)=\dim{(\text{rank}E^-_x)}$ is an invariant of $Z$ and is called an index of $Z$ with respect to the Morse function $f$.\\
\indent
Let $k=m-\dim{Z}-\text{ind}Z$ and let $D^k_x,S^{k-1}_x$ be the unit disc and unit sphere in $E^+_x$. We let $D^k(Z), S^{k-1}(Z)$ represent the fibre bundles on $Z$ with fibres $D^k$ and $S^{k-1}$ respectively.\\
\indent
Now giving a metric to the manifold $M$, there is an $\ape$-small closed tubular neighborhood $U(Z)$ such that $U(Z)$ is homeomorphic to the normal bundle $T^\perp(Z)$. Let $\pat U^+(Z)$ be the intersection of the boundary $\pat U(Z)$ with the unstable manifold $U^+(Z)$ of the gradient flow of $f$. Then $(\ast,U(Z)\backslash\pat U^+(Z))$ is the index pair of $Z$ which is homotopic to $(\ast, U^+(Z)\backslash\pat U^+(Z))$. The latter pointed space is homeomorphic to $(\ast, D^k(Z)\backslash S^{k-1}(Z))$. Let $\widetilde{o(Z)}$ be the local system on $Z$ such that $\mathbb{C}\otimes\widetilde{o(Z)}$ is the orientation bundle $o(Z)$ on $Z$. Let $\tilde{E}$ be a local system of abelian groups on $M$ and $E=\tilde{E}\otimes\mathbb{C}$. Applying the Leray-Hirsch theorem (see[28]) to the fibre-bundle pair $(S^{k-1}(Z), D^k(Z))$, we have
$$
\aligned
&H_*(\ast, D^k(Z)\backslash S^{k-1}(Z);\widetilde{o(Z)}\otimes E|_Z)\cong H_*(S^{k-1}(Z),D^k(Z);\widetilde{o(Z)}\otimes E|_Z)\\
&\cong H_*(Z,\widetilde{o(Z)}\otimes E|_Z)\otimes H_*(\ast, S^k;\mathbb{C})\cong H_{\ast-k}(Z,\widetilde{o(Z)}\otimes E|_Z)
\endaligned
$$
Hence 
$$
H_*(h(Z),\widetilde{o(Z)}\otimes E|_Z)\cong H_{*-k}(Z,\widetilde{o(Z)}\otimes E|_Z)
$$
  \\
{\bf 4.5 Novikov-Morse type inequalities for flows carrying a cocycle}\\
  \\
\indent
Let $v$ be a flow on a compact polyhedron $X$ having a $\pi$-Morse decomposition. Assume that the $\pi$-stable nonwandering set of $v$, $\mathcal{N}_\pi(v)=\{A_n,\cdots,A_1\}$. We have the following definitions.\\ 
 \\
{\sl {\bf Definition 4.5.1 (Isolated $\pi$-stable nonwandering sets)}$\;\;$ Let the flow $v$ be as above. 
A $\pi$-stable nonwandering set $A_i$ in $A$ is called an isolated $\pi$-stable nonwandering set if 
$A_i$ is also an isolated invariant set of $v$.}\\
 \\
\indent
We denote the set containing all the isolated $\pi$-stable nonwandering sets by $\mathcal{IN}_\pi(v)$.\\
 \\
{\sl {\bf Definition 4.5.2}$\;\;$ Let $v$ be a generalized $\a$-flow with respect to 
$A=\{A_n,\cdots,A_1\}$. Then by Theorem 3.5.3, $v$ has an associated $\pi$-Morse decomposition and 
$\mathcal{N}(v)=A$. If each $A_i$ is an isolated $\pi$-stable nonwandering set, then $v$ is called 
an isolated generalized $\a$-flow or an isolated flow that carries a cocycle $\a$.}\\
  \\
\indent
It is easy to see that all the previous examples are isolated (generalized) $\a$-flows.\\
\indent
Before formulating our main result in this section, we introduce some algebraic notations.\\
 \\
{\bf Ideals}\\
 \\
\indent
Let $P_s={\Bbb Z}[t_1,\cdots,t_s]$ be the polynomial ring with $s$ variables over the integers ${\Bbb Z}$. Define a prime ideal $I$ in $P_s$ to be $I=1+\langle t_1,\cdots,t_s\rangle$, where $\langle t_1,\cdots,t_s\rangle$ is the ideal generated by $t_1,\cdots,t_s$. Then the zero set of $I$ is a codimension 1 arithmetic variety in ${\Bbb C}^s$. Let $a\in ({\Bbb C}^*)^s$ be not in the zero set of $I$. Define $I_a$ to be the prime ideal in the polynomial ring $P_s$ consisting of the polynomials vanishing at $a$. By the choice of $a$, the free terms of all the polynomials $f(t)\in I_a$ are divisible by some prime number $p$. Therefore we obtain
$$
I_a\sbs I_p=\langle p\rangle+\langle t_1,\cdots,t_s\rangle
$$
  \\
{\sl {\bf Theorem 4.5.1}$\;\;$ Let $X$ be a compact polyhedron with a metric $d$. Let $v$ be a flow on 
$X$ having a $\pi$-Morse decomposition with rank $\pi=s$ and satisfying 
$\mathcal{IN}_\pi(v)=\mathcal{N}_\pi(v)$. Let $\tilde{E}$ be a local system of free 
abelian groups and let $E={\Bbb C}\otimes\tilde{E}$. If $a\in({\Bbb C}^*)^s$ is not in the arithmetic 
variety associated to the prime ideal $I=1+\langle t_1,\cdots,t_s\rangle\sbs P_s$, then there is a prime 
number $p$ relative to $a$ such that 
\begin{align*}
&\sum_{A\in \mathcal{IN}_\pi(v)}p(h(A);t,{\Bbb Z}_p\otimes J_A^* E)=p(X;t,a^\phi\otimes E)\notag\\
&+(1+t)Q_1(I_a,I_p;t)+(1+t)Q_2(I_p;t)\tag{$4.5.1$}
\end{align*}
where $J_A$ is the inclusion map from the isolated neighborhood of an isolated invariant set to $X$, 
$Q_1(I_a,I_p;t)$ and $Q_2(I_p;t)$ are all polynomials with nonnegative integer coefficients, and 
$a^\phi$ is a complex line bundle determined by the representation of the group $\phi$.}
\\
  \\
Proof $\;\;$ Since $v$ has a $\pi$-Morse decomposition, the lifting flow $\bar{v}$ of $v$ has a relative Morse decomposition 
$$
D=\{\bar{N}^+_0,\bar{A}_n,\cdots,\bar{A}_1\}
$$
and it is easy to see that $\bar{A}_i(i=1,\cdots,n)$ are isolated invariant sets in $\bar{X}_0$ because of that $A_i(i=1,\cdots,n)$ are isolated invariant sets. Now we use the following lemma about the existence of a filtration of a Morse decomposition which was given in [10], with a slight improvement here.\\
 \\
{\sl {\bf Lemma 4.5.2}$\;\;$ Let $S$ be an isolated invariant set in a compact space $X$ and let $\{M_n,\cdots,M_1\}$ be a Morse decomposition of $S$. Then there exists an increasing sequence of compact sets $N_n\sbs N_{n-1},\cdots,\sbs N_0$ such that, for any $j\ge i$, the pair $(N_{j-1},N_i)$ is an index pair for $M_{ji}$. In particular, $(N_n,N_0)$ is an index pair for $S$, and $(N_j,N_{j-1})$ is an index pair for $M_j$. Those compact sets $N_j$ can be chosen to be CW-complexes if the space $X$ is a CW-complex.\\}
 \\
\indent
Applying this lemma, we can get a filtration $\bar{R}_I\times\ape\bar{I}\sbs\bar{N}_n\sbs\bar{N}_{s-1}\sbs\cdots\sbs\bar{N}_1\sbs\bar{X}_0$. Let $\bar{U}$ be a set in $\bar{X}_0$. Let the lifting of $\bar{U}$ in the fundamental domain of the universal covering be $\tilde{U}_0$. We can define a semi-equivariant set in the universal covering $\tilde{X}$ to be $\tilde{U}_+:=\cup_{g\in\pi_1(X)_+}g\cdot\tilde{U}_0$. With that notations we get a $\pi_1(X)_+$-equivariant filtration in $\tilde{X}$:
\begin{equation}
\widetilde{\bar{R}_I\times\ape\bar{I}}\sbs (\tilde{N}_n)_+\sbs (\tilde{N}_{n-1})_+\sbs\cdots\sbs\tilde{X}_+\tag{$4.5.2$}
\end{equation}
This filtration of spaces induces the following filtration of complexes:
$$
0\sbs P^k_s\otimes_{{\Bbb Z}[\pi_1(X)_+]}C_*(\widetilde{\bar{R}_I\times\ape\bar{I}})\sbs P^k_s\otimes_{{\Bbb Z}[\pi_1(X)_+]}C_*((\tilde{N}_n)_+)\cdots\sbs P^k_s\otimes_{{\Bbb Z}[\pi_1(X)_+]}C_*(\tilde{X}_+)
$$
Here the representation of the subring ${\Bbb Z}[\pi_1(X)_+]$ in the space $P^k_s$ is described in section 4.2.\\
\indent
Continue to tensor the above filtration by ${\Bbb Z}_p\otimes_{P_s}$ and use proposition 4.3.3, we get the Poincar\'e polynomial
\begin{align*}
&\sum^{n+1}_{j=1}p({\Bbb Z}_p\otimes J^* E\otimes C_*(\bar{N}_j),{\Bbb Z}_p\otimes J^* E\otimes C_*(\bar{N}_{j-1})\,;\,t)\notag\\ 
&=p({\Bbb Z}_p\otimes J^* E\otimes C_*(\bar{R}_I\times\ape\bar{I})\,,\,{\Bbb Z}_p\otimes J^* E\otimes C_*(\bar{X}_0)\,;\,t)\notag\\
&+(1+t)\sum^n_{j=1}q({\Bbb Z}_p\otimes J^* E\otimes C_*(\bar{R}_I\times\ape\bar{I})\,,\,\tag{$4.5.3$}\\
&{\Bbb Z}_p\otimes J^* E\otimes C_*(\bar{N}_j)\,,\,{\Bbb Z}_p\otimes J^* E\otimes C_*(\bar{N}_{j-1})\,;\,t)\notag
\end{align*}
where we define $\bar{N}_0={\Bbb Z}_p\otimes J^* E\otimes C_*(\bar{X}_0)$ and $\bar{N}_{n+1}={\Bbb Z}_p\otimes J^* E\otimes C_*(\bar{R}_I\times\ape\bar{I})$. Denote the sum term on the right hand side of (4.5.3) by $Q_2({\Bbb Z}_p;t)$, then (4.5.3) can be written as 
\begin{align*} 
 &\sum^{n+1}_{j=1}p({\Bbb Z}_p\otimes J^* E\otimes C_*(\bar{N}_j\,,\,\bar{N}_{j-1})\,;\,t)\\
 &=p({\Bbb Z}_p\otimes J^* E\otimes C_*(\bar{R}_I\times\ape\bar{I})\,,\,{\Bbb Z}_p\otimes J^* E\otimes C_*(\bar{X}_0)\,;\,t)\tag{$4.5.4$}\\ 
 &+(1+t)Q_2({\Bbb Z}_p;t)
\end{align*}
Since $a$ is not in the arithmetic variety associated to the prime ideal $I$, we have $I_a\sbs I_p$, 
where $I_p$ is defined before this theorem. By Theorem 4.2.1, (2), we have 
\begin{align*}
 &p({\Bbb Z}_p\otimes J^* E\otimes C_*(\bar{R}_I\times\ape\bar{I}),{\Bbb Z}_p\otimes J^* E\otimes C_*(\bar{X}_0)\,;t)\\
 &=p({\Bbb Z}_p\otimes_{P_s}(P^k_s\otimes_{{\Bbb Z}[\pi_1(X)_+]}C_*(\tilde{X}_+));t)\\
 &=p({\Bbb Z}_p\otimes_{P_s}D_*;t)\\
 &=\sum_{j\ge 0}b_j({\Bbb Z}_p\otimes_{P_s}D_*)t^j\\
 &=\sum_{j\ge 0}\dim_{Q(P_s/\!I_p)}{H_j(Q(P_s/\!I_p)\otimes_{P_s}D_*)}t^j\\
 &=p(D_*;t,I_p)
\end{align*}
Now by proposition 4.3.4,
\begin{equation}
p(D_*;t,I_p)=p(D_*;t,I_a)+(1+t)Q_1(I_a,I_p;t)\tag{$4.5.5$}
\end{equation}
where 
\begin{equation}
(1+t)Q_1(I_a,I_p;t)=\sum_{j\ge 0}(T_j(I_a,I_p)+T_{j-1}(I_a,I_p))t^j\tag{$4.5.6$}.
\end{equation}
But by Theorem 4.2.1, we have
\begin{align*}
p(D_*;t,I_a)&=\sum_{j\ge 0}\dim_{Q(P_s/\!I_a)}{H_j(Q(P_s/\!I_a)\otimes_{P_s}D_*)}t^j\notag\\
&=\sum_{j\ge 0}\dim_{\Bbb C}(H_j({\Bbb C}_a\otimes_{P_s}D_*))t^j\notag\\
&=\sum_{j\ge 0}\dim_{\Bbb C}(H_j(X;a^\phi\otimes E))t^j\tag{$4.5.7$}
\end{align*}
Combining (4.5.4), (4.5.5) and (4.5.7) and using the notation $p(X;t,a^\a\otimes E)$ for $p(D_*;t,I_a)$, we get 
\begin{align*}
&\sum^{n+1}_{j=1}p({\Bbb Z}_p\otimes J^* E\otimes C_*(\bar{N}_j,\bar{N}_{j-1});t)\\
&=p(X;t,a^\a\otimes E)+(1+t)Q_1(I_a,I_p;t)+(1+t)Q_2(I_p;t)\tag{$4.5.8$}
\end{align*}
\indent
Note that $(\bar{N}_j,\bar{N}_{j-1})$ is the index pair for the isolated invariant set $\bar{A}_j$ and hence is the index 
pair for the isolated invariant set $A_j\in \mathcal{IN}_\pi(v)$. (4.5.8) induces the equality (4.5.1)\qed\\
  \\
{\sl {\bf Corollary 4.5.3 (Euler-Poincar\'e formula)}$\;\;$ Under the hypothesis of Theorem 4.5.1,  
\begin{align*}
\sum_{A\in\mathcal{IN}_\pi(v)} p(h(A)\,;\,-1\,,\,{\Bbb Z}_p\otimes J^*_A E)=p(X\,,\,-1\,,\,a^\phi\otimes E)\tag{$4.5.9$}
\end{align*}
In particular, if $E$ is a trivial line bundle, then
\begin{align*}
\sum_{A\in\mathcal{IN}_\pi(v)}p(h(A)\,;\,-1\,,\,{\Bbb Z}_p)=\chi(X)\tag{$4.5.10$}
\end{align*}
for any prime number $p$. Here $\chi(X)$ is the Euler chracteristic number of the compact polyhedron.}\\
 \\
Proof$\;\;$ We only need to prove the second conclusion. If $E$ is a trivial line bundle, then
\begin{align*}
p(X\,;\,-1\,,\,a^\a)&=\sum_{j\ge 0}\dim_{\Bbb C}(H_j(X;a^\a))(-1)^j\\
&=\sum_{j\ge 0}\dim_{\Bbb C}(H_j({\Bbb C}_a\otimes_{P_s}D_*))(-1)^j\\
&=\sum_{j\ge 0}\hbox{rank}_{P_s}(P_s\otimes_{{\Bbb Z}[\pi_1(X)_+]}C_j(\tilde{X}_+))(-1)^j\\
&=\sum_{j\ge 0}\hbox{rank}_{\Bbb Z}C_j(X)(-1)^j=\chi(X)
\end{align*}
Here we use the fact $D_j$ is a free module over the ring $P_s$.\qed\\
\indent
If $v$ is a generalized $\a$-flow with respect to the nonwandering set $A=\{A_n,\cdots,\\ A_1\}$, then by Theorem 3.5.3 $v$ has a $\pi$-Morse decomposition with rank $\pi=\hbox{rank}\;\a$, and $A=\mathcal{N}_\pi(v)$. Therefore the following theorem is the direct consequence of Theorem 4.5.1.\\
  \\
{\sl {\bf Theorem 4.5.4}$\;\;$ Let $X$ be a compact polyhedron with a metric $d$. Let $v$ be a generalized 
$\a$-flow with respect to the nonwandering set $A=\{A_n,\cdots, A_1\}$ and assume that each $A_i$ is an 
isolated invariant set in $X$. Let $\tilde{E}$ be a local system of free abelian groups and let 
$E={\Bbb C}\otimes\tilde{E}$. If $a\in({\Bbb C}^*)^s$ is not in the arithmetic variety associated to 
the prime ideal $I=1+\langle t_1,\cdots,t_s\rangle\sbs P_s$, then there is a prime number $p$ relative to $a$ such that 
\begin{align*}
&\sum_{A_i\in A}p(h(A_i);t,{\Bbb Z}_p\otimes J_{A_i}^* E)=p(X;t,a^\a\otimes E)\notag\\
&+(1+t)Q_1(I_a,I_p;t)+(1+t)Q_2(I_p;t)\tag{$4.5.11$}
\end{align*}
where $J_{A_i}$ is the inclusion map from the isolated neighborhood of an isolated invariant set to $X$, and $Q_1(I_a,I_p;t)$ and $Q_2(I_p;t)$ are both polynomials with nonnegative integer coefficients.}\\
  \\
{\bf 4.6 Applications to special flows carrying a cocycle}\\
  \\
\indent
In this section, we will start from the general formula (4.5.1) and use the Conley index introduced in section 4.3 to get some Novikov-Morse type inequalities for some special important flows.\\
\indent
If $v$ is a generalized $\a$-flow with respect to the isolated nonwandering set $A=\{A_n,\cdots,A_1\}$ with $\a$ being a trivial cocyle, then $v$ has actually a ``global'' gradient-like structure in view of Theorem 3.4.1. Therefore the formula (4.5.1) is a Morse type inequality for an isolated invariant set which has the Morse decomposition $A=\{A_n,\cdots,A_1\}$. This result was given in the paper [10].\\
\indent
Now we consider the case that $\a$ is a non-trivial cocycle. \\
\indent
The important example is the $\o$-flow generated by a closed Morse 1-form $\o$ (see Example 3.4.2). One can also consider the Bott type Novikov inequalities. The following two theorems were given in [13] for $s=1$.\\
  \\
{\sl {\bf Theorem 4.6.1}$\;\;$ Let $X$ be an oriented closed smooth manifold and $\o$ a rank $s$ closed Morse 1-form. Let $\tilde{E}$ be a local system of free abelian group and $E={\Bbb C}\otimes\tilde{E}$. Assume that $a\in({\Bbb C}^*)^s$ is not in the arithmetic variety associated to the prime ideal $I=1+\langle T_1,\cdots,T_s\rangle\sbs P_s$. Then the number $c_j(\o)$ of zeros of $\o$ having index $j$ satisfies
\begin{align*}
c_j(\o)\ge& \frac{\dim_{\Bbb C}H_j(M,a^\o\otimes E)}{\dim E}\\
\sum^j_{i=0}(-1)^i c_{j-i}(\o)\ge& \sum^j_{i=0}(-1)^i\frac{\dim_{\Bbb C}H_{j-i}(M,a^\o\otimes E)}{\dim E}\tag{$4.6.1$}
\end{align*}
for $j=0,1,\cdots,m$.\\
}
 \\
{\sl {\bf Theorem 4.6.2}$\;\;$ Let $\o$ be a closed 1-form with Bott type nondegenerate zero sets with 
rank $s$. Let the quantities $\tilde{E},E$ and $a$ be the same as in Theorem 4.6.1. Then there is a prime
 $p$ such that 
\begin{align*}
&\sum_Z\sum_{j\ge 0}\dim_{{\Bbb Z}_p}H_j(Z,{\Bbb Z}_p\otimes \tilde{E}|_Z\otimes o(Z))t^{j+ind(Z)}\\
&=p(X,t;a^\o\otimes E)+(1+t)Q(E,t)\tag{$4.6.2$}
\end{align*}
where $o(Z)$ is the orientation bundle on $Z$ and the sum in the left hand side is taken for any zero manifold of $\o$ and $Q(E,t)$ is a polynomial with nonnegative coefficients.\\
}
 \\
Proof of Theorem 4.6.1 and 4.6.2$\;\;$ From Example 3.4.2 and 3.4.7, we know the flow generated by 
the dual vector field of $\o$ is a flow carrying the cocycle $\o$. Since the fixed points of these flow are just the (Bott type) nondegenerate zero points(sets), by section 4.3, we have the following result.\\
\indent
(1) For a zero point $x_0$, 
\begin{align*}
H_l(h(x_0)&,{\Bbb Z}_p\otimes J^*E)\cong H_l(*,S^j;{\Bbb Z}_p\otimes J^*E)\\
&\cong \Big\{
\begin{array}{ll}
\otimes^{\dim E}_1{\Bbb C}\otimes{\Bbb Z}_p & \;l=j\\
0 & l\neq j
\end{array}
\end{align*}
\indent
(2) For a Bott type zero set $Z$, it has 
$$
H_*(h(Z),\widetilde{o(Z)}\otimes E|_Z)\cong H_{*-k}(Z,\widetilde{o(Z)}\otimes E|_Z)
$$
Now applying formula (4.5.1), we get the conclusions.\qed\\
\indent
In view of Example 3.4.3, the formulas (4.6.1) and (4.6.2) not only hold for the closed 1-forms, they also hold for a family of flows which is a small perturbation of the flow generated by $\o$ because the Conley indices of fixed point sets are invariant under the perturbation. Hence we have\\
 \\ 
{\sl {\bf Corollary 4.6.3}$\;\;$ Let $v_\ape$ be a family of flows with parameter $\ape$ in Example 3.4.3. Under the hypothesis of Theorem 4.6.1 and 4.6.2, the formula (4.6.1) holds if $\o$ is a closed Morse 1-form and the formula (4.6.2) holds if $\o$ is a closed 1-form having Bott type zero sets.\\
}
 \\
\indent
If we let $E$ be the trivial line bundle and let each entry $a_i$ of $a=(a_1,\cdots,a_s)\in ({\Bbb C}^*)^s$ be a transcendental number, then from Theorem 4.6.1 we have the following corollary\\
  \\
{\sl {\bf Corollary 4.6.4 (Classical Novikov inequality)}$\;\;$ Let $X$ be an oriented closed smooth manifold and $\o$ a rank 1 Morse closed 1-form. Then the numbers $c_j(\o)$ of zeros of $\o$ having index $j$ satisfy
\begin{align*}
c_j(\o)\ge& b_j([\o]) \\
\sum^j_{i=0}(-1)^i c_{j-i}(\o)\ge&\sum^j_{i=0}(-1)^i b_{j-i}([\o])\tag{$4.6.3$}
\end{align*}
for $j=0,1,\cdots,m$.\\
}\\   
 \\
\indent
In the appendix, we will give a proof of the refined Novikov inequalities that  include the information of the torsion part.\\
\indent
In Example 3.4.5, we have considered a flow carrying a cohomology class. The existence of such flow gives a vanishing theorem for the Novikov numbers.\\
  \\
{\sl {\bf Theorem 4.6.5 (Vanishing theorem)}$\;\;$ Let $X$ be a compact polyhedron with a metric $d$. If there exists a flow carrying a cohomology class $[\a]$ on $X$, then  
$$
b_i([\a])=0,\;\;\forall i=0,1,\cdots,m
$$
}\\
 \\
{\sl {\bf Theorem 4.6.6}$\;\;$ Let $v$ be a $\a$-Morse-Smale flow on a closed oriented manifold $X$ and let its nonwandering set be $A=\{A_n,\cdots,A_1\}$ consisting of the hyperbolic fixed points and hyperbolic periodic orbits. Let $c_j$ be the number of hyperbolic fixed points with index $j$, $a_j$ be the number of the hyperbolic periodic orbits with index $j$, and $\mu_j=c_j+a_j+a_{j+1}$. Then 
\begin{align*}
\mu_j\ge& b_j([\a])\\
\sum^j_{i=0}(-1)^{j-i}\mu_i\ge&\sum^j_{i=0}(-1)^{j-i}b_j([\a])\tag{$4.6.4$}
\end{align*}
for $j=0,1,\cdots,m$.}\\
  \\
Proof$\;\;$ The conclusion is obtained by applying Theorem 4.5.4 with $E$ being the orientation bundle of $X$ and the calculation of the homology of the Conley indices of the hyperbolic fixed points and the orientable or unorientable hyperbolic periodic orbits.\\
\indent
For the calculation of the Conley indices, we can use the results in section 4.3. \qed\\
 \\
{\bf\large  Acknowledgement}\\
 \\
\indent
We thank Michael Farber for explaining his work on the Novikov inequalities for us and for stimulating discussions. The first author also wants to thank Guofang Wang for his discussions with him.\\

\begin{center}
{\bf\large Appendix}\\
\end{center}
\vskip 10pt
\indent
In this section we will provide a proof of the classical refined Novikov inequalities which includes also the information from the torsion parts of the homology groups. The starting point is the equality (4.5.1). If the $\a$-flow we consider is generated by the dual vector field of the closed Morse 1-form, then the connecting homomorphism is easily determined.\\
\indent
Let $X$ be a CW-complex with subcomplex $Y$ and let $X^k$ be the $k$-skeleton of $X$, and define
$$
X^k_Y=X^k\cup Y.
$$
Then we have a filtration relative to $Y$,
$$
Y=X^{-1}_Y\sbs X^0_Y\sbs \cdots\sbs X^m_Y=X
$$
This filtration induces the filtration of the complex $C_*(X)$:
$$
C_*(X^{-1}_Y)\sbs C_*(X_Y^0)\sbs\cdots\sbs C_*(X^m_Y)=C_*(X)
$$
\indent
Let $G$ be a field and set the subcomplex $C_j=G\otimes C_*(X^{m-j}_Y),\;j=0,1,\cdots,m,m+1$. Then we obtain the filtration
$$
0\sbs C_{m+1}\sbs C_m\sbs\cdots\sbs C_0.
$$
For this filtration of complexes, we get from proposition 4.3.3 the following formula:
\begin{align*}
\sum^{m+1}_{j=1}p(C_j,C_{j-1};t)=p(C_{m+1},C;t)+(1+t)\sum^m_{j=1}q(C_{m+1},C_j,C_{j-1};t),\tag{$1$}  
\end{align*}
where
\begin{align*}
&p(C_j,C_{j-1};t)=\sum_{l\ge 0}b_l(C_j,C_{j-1})t^l\\
&q(C_{m+1},C_j,C_{j-1};t)=\sum_{l\ge 0}d_l(C_{m+1},C_j,C_{j-1})t^l\\
&d_l(C_{m+1},C_j,C_{j-1})=\hbox{rank}(\hbox{im}\d_l)
\end{align*}
and $\d_l:H_l(C_j,C_{j-1})\longrightarrow H_{l-1}(C_{m+1},C_j)$ is the connectiong homomorphism of the triple of complexes $(C_{m+1},C_j,C_{j-1})$.\\
\indent
Because of the choice of the filtration, we have 
\begin{align*}
H_l(C_j,C_{j-1})\cong &H_l(G\otimes C_*(X^{m-j}_Y,X^{m-j+1}_Y)\\
\cong&\Big\{
\begin{array}{ll}
\otimes^{k_{m-j+1}}_1 G & \;l=m-j+1\\
0 & \;l\neq m-j+1
\end{array}\tag{$2$}
\end{align*}
where $k_{m-j+1}$ is the number of $m-j+1$-cells in $X$ but not in $Y$.\\
\indent
We also have 
$$
H_{l-1}(C_{m+1},C_j)\cong H_{l-1}(G\otimes G_*(Y,X^{m-j}).
$$
The connecting homomorphism $\d_l$ can be identified with the homomorphism going from $H_l(G\otimes C_*(X^{l-1}_Y,X^l_Y))$ to $H_{l-1}(G\otimes C_*(Y,X^{l-1})$. Consider the commutative diagram
$$
\diagram
H_l(G\otimes C_*(X_Y^{l-1}, X_Y^l)) \dto^{\d_l}\drto^{d_l}\\
H_{l-1}(G\otimes C_*(Y,X_Y^{l-1})) \rto^{i_*}& H_{l-1}(G\otimes C_*(X_Y^{l-2}, X^{l-1}_Y))
\enddiagram
$$
Here $i:(Y,X^{l-1}_Y)\longrightarrow (X^{l-2}_Y,X^{l-1}_Y)$ is the inclusion map of the space pair. Therefore $d_l=i_*\cdot\d_l$ and $d_l$ is just the boundary operator in the cellular complex $G\otimes W_*(X,Y)$ with coefficients in $G$. Thus 
\begin{align*}
\hbox{rank}_G(\hbox{im}d_l)=\hbox{rank}_G(\hbox{im}(i_*\cdot\d_l))\le \hbox{rank}_G(\hbox{im}\d_l)\tag{$3$}
\end{align*}
and the inequality holds if and only if $i_*$ is an injection.\\
\indent
Therefore using (1)-(3), we have 
\begin{align*}
\sum^m_{j=0}k_j t^j\ge p(G\otimes W_*(Y,X)\,;\,t)+(1+t)\sum^m_{j=1}\hbox{rank}_G(\hbox{im}d_j) t^j\tag{$4$}
\end{align*}
where $d_j:G\otimes W_j(Y,X)\longrightarrow G\otimes W_{j-1}(Y,X)$ is the boundary operator. \\
\indent
Let $X=\bar{X}_0,\;Y=\bar{R}_I\times\ape\bar{I}$ and $G={\Bbb Z}_p$, (4) becomes
\begin{align*}
\sum^m_{j=0}k_j t^j\ge p({\Bbb Z}_p\otimes W_*(Y,X)\,;\,t)+(1+t)\sum^m_{j=1}\hbox{rank}_{{\Bbb Z}_p}(\hbox{im}d_j) t^j\tag{$5$}
\end{align*}
\indent
Now we continue the progress from (4.5.4), to get 
\begin{align*}
\sum^m_{j=0}k_j t^j\ge& p(X;t,a^\a)+(1+t)\sum^m_{j=1}\hbox{rank}_{{\Bbb Z}_p}(\hbox{im}d_j) t^j\\
+&\sum_{j\ge 0}(T_j(I_a,I_p)+T_{j-1}(I_a,I_p))t^j\tag{$6$}
\end{align*}  
Note that 
\begin{align*}
T_j(I_a,I_p)=\hbox{rank}_{\Bbb C}(\hbox{im}d_{j+1}(I_a))-\hbox{rank}_{{\Bbb Z}_p}(\hbox{im}d_{j+1}(I_p)),
\end{align*}
where $d_j(I_a):Q(P_s/\!I_a)\otimes_{P_s}D_j\to Q(P_s/\!I_a)\otimes_{P_s}D_{j-1}$ and $d_j(I_p):Q(P_s/\!I_p)\otimes_{P_s}D_j\to Q(P_s/\!I_p)\otimes_{P_s}D_{j-1}$ are the boundary operators.\\
\indent
Therefore, from (6), we have
\begin{align*}
\sum^m_{j=0}k_j t^j\ge& p(X;t,a^\a)+\sum_{j\ge 0}(\hbox{rank}_{\Bbb C}(\hbox{im}d_j(I_a))+\hbox{rank}_{\Bbb C}(\hbox{im}d_{j-1}(I_a))) t^j\tag{$7$}
\end{align*}
\indent
Define the minimal number of the generators of the torsion part of $H_j(X,a^\o)$ by $q_j([\o])$.\\
\indent
Now we can prove the following classical Novikov inequality\\
 \\
{\sl{\bf Theorem (Novikov inequality)}$\;\;$ Let $\o$ be a closed Morse 1-form on a closed oriented manifold $X$. Let $c_j(\o)$ be the number of zero points of $\o$ with index $j$, then 
\begin{align*}
c_j(\o)\ge b_j([\o])+q_j([\o])+q_{j-1}([\o]),
\end{align*}
for $j=0,1,\cdots,m$.}\\
 \\
Proof$\;\;$ Lift the 1-form $\o$ to the covering space $\bar{X}$ such that $\o$ becomes a Morse function $\bar{f}$. By considering the handle decomposition of the fundamental domain $\bar{X}_0$ with respect to $\bar{f}$, it is easy to see that $k_j=c_j$ in (7). Since $\hbox{rank}(\hbox{im}d_j(I_a))\ge q_j([\o])$, (7) then induces the result.\qed\\

\end{document}